\documentclass[12pt]{amsart}
\usepackage{amsmath, amssymb}
\newtheorem{theorem}{Theorem}
\newtheorem{example}{Example}
\newtheorem{proposition}{Proposition}
\newtheorem{corollary}{Corollary}
\newtheorem{remark}{Remark}
\newtheorem{definition}{Definition}
\newtheorem{lemma}{Lemma}
\newcommand{\sm}{\setminus}
\newcommand{\K}{{\mathcal K}}
\newcommand{\Z}{{\mathbb Z}}
\newcommand{\R}{{\mathbb R}}

\newcommand{\CP}{{\mathbb {CP}}}

\author{V.~A.~Vassiliev}

\email{vva@mi.ras.ru}

\thanks{Supported in part by RFBR (project 98-01-00555a)
and NWO-RFBR grant (project 047-008-005)}

\title{On combinatorial formulas for cohomology of spaces of knots}

\date{Revised version published in 2001}

\begin{document}

\begin{abstract}
We develop homological techniques for finding explicit combinatorial
expressions of finite-type cohomology classes of spaces of knots in $\R^n, n
\ge 3,$ generalizing Polyak--Viro formulas \cite{PV} for invariants (i.e.
0-dimensional cohomology classes) of knots in $\R^3$.

As the first applications we give such formulas for the (reduced mod 2) {\em
generalized Teiblum--Turchin cocycle} of order 3 (which is the simplest
cohomology class of {\em long knots} $\R^1 \hookrightarrow \R^n$ not reducible
to knot invariants or their natural stabilizations), and for all integral
cohomology classes of orders 1 and 2 of spaces of {\em compact knots} $S^1
\hookrightarrow \R^n$. As a corollary, we prove the nontriviality of all these
cohomology classes in spaces of knots in $\R^3.$
\end{abstract}

\maketitle

\section{Introduction}
\label{int}

There is a wide family of cohomology classes of spaces of knots $S^1
\hookrightarrow \R^n$ ($n \ge 3$), called {\em finite-type cohomology classes};
see \cite{V1}, \cite{fasis}, \cite{bjo}. For $n>3$ they cover all of the
cohomology group of the space of knots in $\R^n$, for $n=3$ their 0-dimensional
part are the finite-type knot invariants.

These classes are defined as linking numbers (in the space of all smooth maps
$S^1 \to \R^n$) with appropriate cycles (of infinite dimension but finite
codimension) in the {\em discriminant space} $\Sigma$ (cf. \cite{arnold}); in
our case this space consists of maps which are not smooth embeddings. The group
of all such classes is filtered by their {\em orders} induced by some
filtration of (some resolution of) the discriminant: roughly speaking, the
order of a cohomology class indicates how much complicated strata of $\Sigma$
participate in the definition of its dual variety.

In \cite{PV}, M.~Polyak and O.~Viro have proposed some combinatorial formulas
for the finite-type invariants of knots in $\R^3$. Later, M.~Goussarov has
proved that any finite-type invariant can be represented by a formula of this
type, see \cite{GPV}.

We describe some calculus for finding (and proving) combinatorial formulas for
arbitrary finite type cohomology classes, in particular show what the answers
can look like. {\em Any such combinatorial formula is nothing else than some
semialgebraic chain in the space of maps $S^1 \to \R^n,$ whose boundary lies in
$\Sigma$ and our cohomology class is equal to the linking number with this
boundary.} We introduce several natural families of semialgebraic subvarieties
of the space of such maps, of which the desired chains are built. These
varieties are defined by easy differential geometrical conditions; they arise
naturally in the direct calculation of the main spectral sequence converging to
the (finite type) cohomology group of the space of knots. It is not surprising
that some elements of this calculus repeat pictures from \cite{PV}, \cite{GPV},
and also from the A.~B.~Merkov's works on invariants of plane curves \cite{M},
\cite{MM}.

We accomplish these calculations explicitly for several cohomology classes of
low orders. Before describing them three remarks more.

{\bf 1. Long and compact knots.} We shall distinguish two kinds of knot spaces.
The {\em compact knots} in $\R^n$ are any smooth embeddings $S^1 \to \R^n,$
while the {\em long knots} are the smooth embeddings $\R^1 \to \R^n$ coinciding
with a standard linear embedding outside some compact subset in $\R^1$. The
invariants of knots of both types in $\R^3$ naturally coincide, but generally
the cohomology ring of the space of compact knots is more complicated: it is
built of the similar ring for long knots (playing the role of a "coefficient"
ring) and homology groups of the space $S^1$ and certain its configuration
spaces.

{\bf 2. Stabilization.} If numbers $n$ and $m$ are of the same parity, then the
theories of (finite type) cohomology groups of spaces of knots in $\R^n$ and
$\R^m$ are very similar. Namely, the first terms of spectral sequences
calculating both groups and generated by the natural filtration of resolved
discriminants coincide up to shifts of indices:
\begin{equation}
\label{stabil} E_1^{p,q-pn}(\R^n) \simeq E_1^{p,q-pm}(\R^m).
\end{equation}
Moreover, for spectral sequences calculating $\Z_2$-cohomology groups this
identity is true also if $n$ and $m$ are of different parities. M.~Kontsevich
has proved (but not published) that in the case of complex coefficients our
spectral sequence degenerates at the first term: $E_\infty^{p,q} \equiv
E_1^{p,q}$, therefore also the limit groups of finite type cohomology classes
are very similar. (I conjecture that in the case of long knots the similar
degeneration holds also for any coefficients.)

{\bf 3.} This paper is very much a work in the differential geometry of spatial
curves and their projections to different subspaces, although almost all
results of this kind are hidden in the formulas for boundary operators in our
homological calculations.

\subsection{Results for long knots.}
\label{reslong}

Accordingly to \cite{V1}, \cite{tetu}, \cite{bjo}, \cite{fasis}, all cohomology
classes of orders $\le 3$ of the space of long knots in $\R^n,$ $n \ge 3,$ are
as follows.

\begin{proposition}
\label{longres} There are no cohomology classes of order 1. The classes of
order 2 are only in dimension $2n-6$ and form a group isomorphic to $\Z$ $($for
$n=3$ it is generated by the simplest knot invariant$)$. In order 3 additional
classes can be in exactly two dimensions more: $3n-9$ and $3n-8$. In dimension
$3n-9$ they form a group isomorphic to $\Z$ $($for $n=3$ it is generated by the
next simple knot invariant$)$. In dimension $3n-8$ the same is true if $n >3,$
and for $n=3$ the similar $($1-dimensional$)$ cohomology group is cyclic
$($maybe of order 1 or $\infty )$.
\end{proposition}

It was conjectured in \cite{fasis}, \cite{bjo} that the latter group for $n=3$
also is isomorphic to $\Z$; we shall prove it in the present work.

For any $n$ we call the generator of this $(3n-8)$-dimensional cohomology group
the {\em Turchin--Teiblum cocycle}. In the case of odd $n$ its existence was
discovered by D.~M.~Teiblum and V.~E.~Turchin about 1995 (\cite{tetu}). Its
(quite different) superanalog for even $n$ was found in \cite{fasis},
\cite{bjo}. However, all these works contain only an implicit proof of the
existence of such a class: namely, the calculation of the third column of our
spectral sequence (which is responsible for the third order cohomology classes
and is isomorphic to $\Z$ for exactly two values of $q$), and the remark that
all further differentials acting from or to this column are trivial by some
dimensional reasons.

In \S \ref{proof1} we prove the following combinatorial expression for this
class reduced mod 2.

Let us choose a direction "up" in $\R^n$, and say that a point $x \in \R^n$ is
{\em above} the point $y$ if the vector $\overrightarrow{(yx)}$ has the chosen
direction. Let $\R^{n-1}$ be the quotient space of $\R^n$ by this direction,
and ${\bf p}: \R^n \to \R^{n-1}$ the corresponding projection. We choose a
direction "to the right" in $\R^{n-1}$ and say that the point $x \in \R^n$ {\em
is to the right} of the point $y$ if the vector $\overrightarrow{({\bf
p}(y),{\bf p}(x))} \in \R^{n-1}$ has this chosen direction.

\begin{theorem}
\label{main} For any $n \ge 3,$ the value of the reduced mod 2 Teiblum--Turchin
class on any generic $(3n-8)$-dimensional singular cycle in the space of long
knots in $\R^n$ is equal to the parity of the number of points of this cycle
corresponding to such knots $f:\R^1 \to \R^n$ that one of three holds:

a) there are five points $a<b<c<d<e$ in $\R^1$ such that $f(a)$ is above
$f(d)$, and $f(e)$ is above $f(c)$ and $f(b)$;

b) there are four points $a<b<c<d$ in $\R^1$ such that $f(a)$ is above $f(c)$,
$f(b)$ is below $f(d)$, and the projection of the derivative $f'(b)$ to
$\R^{n-1}$ is directed to the right;

c) there are three points $a<b<c$ in $\R^1 $such that $f(a)$ is above $f(b)$
but below $f(c)$, and the "exterior" angle in $\R^{n-1}$ formed by projections
of $f'(a)$ and $f'(b)$ contains the direction "to the right" $($i.e. this
direction is equal to a linear combination of these projections, and at least
one of coefficients in this combination is non-positive$)$.
\end{theorem}

We prove this theorem in \S \ref{proof1}. In the next works I am planning to
accomplish all the same calculations taking respect on the orientations, and
thus to obtain similar results with integer coefficients.

\begin{corollary}
\label{nontriv} The group of order 3 one-dimensional cohomology classes of the
space of long knots in $\R^3$ is free cyclic and generated by the
$($integral$)$ Teiblum--Turchin class.
\end{corollary}

More precisely, let us consider the connected sum of two equal (long) trefoil
knots in $\R^3$ and a path in the space of knots connecting this knot with
itself as in the proof of the commutativity of the knot semigroup: we shrink
the first summand, move it "through" the second, and then blow up again.

\begin{proposition}
\label{realiz} This closed path in the space of long knots has an odd number of
intersection points with the union of three varieties indicated in items a, b
and c of Theorem \ref{main}.
\end{proposition}

The proof will be given in \S~\ref{proreal}.
\medskip

On the other hand, for any $n$ the Teiblum--Turchin cocycle is a well-defined
integral cohomology class, and Corollary \ref{nontriv} is proved.

\subsection{Answers for compact knots}
\label{comp}

Nontrivial cohomology classes in the space of {\em compact knots} $S^1
\hookrightarrow \R^n$ appear already in filtrations 1 and 2. We assume that a
cyclic coordinate in $S^1$, i.e. an identification $S^1 \simeq \R^1/2\pi \Z$,
is fixed.

\begin{proposition}[see \cite{arman}, \cite{bjo}]
\label{ordone} For any $n \ge 3$ the group of $\Z_2$-cohomology classes of
order 1 of the space of compact knots in $\R^n$ is nontrivial only in
dimensions $n-2$ and $n-1$, and is isomorphic to $\Z_2$ in these dimensions.
Moreover, for $($only$)$ even $n$ similar integral cohomology groups in these
dimensions are isomorphic to $\Z.$ The generator of the $(n-2)$-dimensional
group is Alexander dual to the set ${\mathcal L}$ of discriminant maps $S^1 \to
\R^n$ gluing together some two {\em opposite} points of $S^1$, and the
generator of the $(n-1)$-dimensional group is dual to the set of maps gluing
some {\em chosen} opposite points, say $0$ and $\pi$.
\end{proposition}

\begin{proposition}[see \cite{fasis}, \cite{bjo}]
\label{ordtwo} Additional classes of order 2 exist in exactly two dimensions:
$2n-6$ and $2n-3$. In dimension $2n-6$ they for any $n$ form a group isomorphic
to $\Z$ $($for $n=3$ it is generated by the simplest knot invariant$)$. The
group in dimension $2n-3$ is isomorphic to $\Z$ for $n>3$ and cyclic for $n=3$;
its generator is Alexander dual to the cycle in the discriminant, whose {\em
principal part} $($see Definition \ref{fitype} in \S \ref{method} below$)$ in
the double self-intersection of $\Sigma$ is swept out by such maps $f:S^1 \to
\R^n$ that for some $\alpha \in S^1$ we have $f(\alpha) = f(\alpha+\pi),$
$f(\alpha+\pi/2) =f(\alpha+3\pi/2).$
\end{proposition}

Below we prove in particular that for $n=3$ the last group also is free cyclic,
see Corollary \ref{realiz3}. Now we give explicit combinatorial formulas for
all classes mentioned in Propositions \ref{ordone} and \ref{ordtwo}.

\begin{theorem}
\label{comain} For any $n \ge 3$, the values of any of these four basic
cohomology classes on any generic cycle of corresponding dimension in the space
${\mathcal K}_n \setminus \Sigma$ of compact knots $S^1 \hookrightarrow \R^n$
is equal to the number of points of this cycle, corresponding to knots
satisfying the following conditions $($and in the case of integer coefficients
taken with appropriate signs$)$.

A. For the $(n-1)$-dimensional class of order $1$: projections of $f(0)$ and
$f(\pi)$ into the plane $\R^{n-1}$ coincide, and $f(0)$ is above $f(\pi)$.

B. For the $(n-2)$-dimensional class of order 1, one of the following two
conditions:

a) there is a point $\alpha \in [0,\pi)$ such that the projections of
$f(\alpha)$ and $f(\alpha+\pi)$ to $\R^{n-1}$ coincide, and moreover
$f(\alpha)$ is above $f(\alpha+\pi)$;

b) the projection of the point $f(0)$ to $\R^{n-1}$ lies "to the right" from
the projection of $f(\pi)$.

C. For the $(2n-3)$-dimensional class of order $2$, one of following two
conditions:

a) there is a point $\alpha \in [0,\pi/2)$ such that projections of $f(\alpha)$
and $f(\alpha+\pi)$ to $\R^{n-1}$ coincide, projections of $f(\alpha+\pi/2)$
and $f(\alpha+3\pi/2)$ to $\R^{n-1}$ coincide, and additionally $f(\alpha+\pi)$
is above $f(\alpha)$ and $f(\alpha+\pi/2)$ is above $f(\alpha+3\pi/2)$;

b) projections of $f(0)$ and $f(\pi)$ to $\R^{n-1}$ coincide, $f(\pi)$ is above
$f(0)$, and the projection of $f(\pi/2)$ to $\R^{n-1}$ lies "to the right" from
the projection of $f(3\pi/2)$.

D. For the $(2n-6)$-dimensional class of order $2$, one of two conditions:

a) there are four distinct points $\alpha, \beta, \gamma, \delta \in S^1$
$($whose cyclic coordinates satisfy $0 \le \alpha < \beta < \gamma < \delta <
2\pi)$ such that projections of $f(\alpha)$ and $f(\gamma)$ to $\R^{n-1}$
coincide, projections of $f(\beta)$ and $f(\delta)$ to $\R^{n-1}$ coincide, and
additionally $f(\gamma)$ is above $f(\alpha)$ and $f(\beta)$ is above
$f(\delta)$.

b) If $n=3$ then second condition is void $($and we have only the first one
coinciding with the combinatorial formula from \mbox{\cite{PV}}$)$, but for
$n>3$ we have additional condition: there are three distinct points $ \beta,
\gamma, \delta$ $($whose cyclic coordinates satisfy $0 < \beta < \gamma <
\delta < 2\pi)$ such that projections of $f(\gamma)$ and $f(0)$ to $\R^{n-1}$
coincide, $f(\gamma)$ is above $f(0)$, and the projection of $f(\delta)$ to
$\R^{n-1}$ lies "to the right" of the projection of $f(\beta)$.
\end{theorem}

Proofs see in \S \ \ref{proof2}.
\medskip

\begin{corollary}
\label{realiz2} For any $n\ge 3$, the basic class of order $2$ and dimension
$2n-3$ takes value $\pm 1$ on the submanifold of the space of knots, consisting
of all naturally parametrized great circles of the unit sphere in $\R^n$.
\end{corollary}

Indeed, the variety a) of statement C does not intersect this submanifold, and
variety b) has with it exactly one intersection point. \quad $\square$
\medskip

In the case of even $n$ the fact that this variety in the space of knots is not
homologous to zero was proved in \cite{cotta} by very different methods.
\medskip

\begin{corollary}
\label{realiz3} The group of $(2n-3)$-dimensional cohomology classes of order 2
is free cyclic for $n=3$ as well. \quad $\square$
\end{corollary}

I am indebted to A.~B.~Merkov very much for many interesting conversations.

\section{Methodology and nature of combinatorial expressions.}
\label{method}

In fact, our main purpose is to develop a general method of finding
combinatorial formulas of this type.

Any such formula is just a relative cycle in the space of knots (modulo the
discriminant $\Sigma$) whose boundary in $\Sigma$ is Alexander dual to our
cohomology class. The problem is to construct such a variety explicitly and as
simply as possible.

Our method of doing it consists in the conscientious calculation of our
spectral sequence. In this subsection we outline the definition of this
sequence and this calculation. This spectral sequence for spaces of knots is
very analogous to that calculating the cohomology of complements of plane
arrangements (see \cite{congr}); let us demonstrate their main common features
on the latter more simple example.

\subsection{Simplicial resolutions for plane arrangements}

Let $L \subset \R^N$ be an {\em affine plane arrangement,} i.e. the union of
finitely many affine planes $L_i$ of any dimensions, $L = \bigcup_{i=1}^k L_i.$
The cohomology group of its complement is Alexander dual to the homology group
of $L$: $H^j(\R^N \setminus L) \simeq \bar H_{N-j-1}(L);$ here $\bar H_*$
denotes the {\em Borel--Moore homology group}, i.e. the homology group of the
one-point compactification reduced modulo the added point. To calculate the
latter group it is convenient to use the {\em simplicial resolution} of $L$
(which is a continuous version of the combinatorial formula of inclusions and
exclusions).

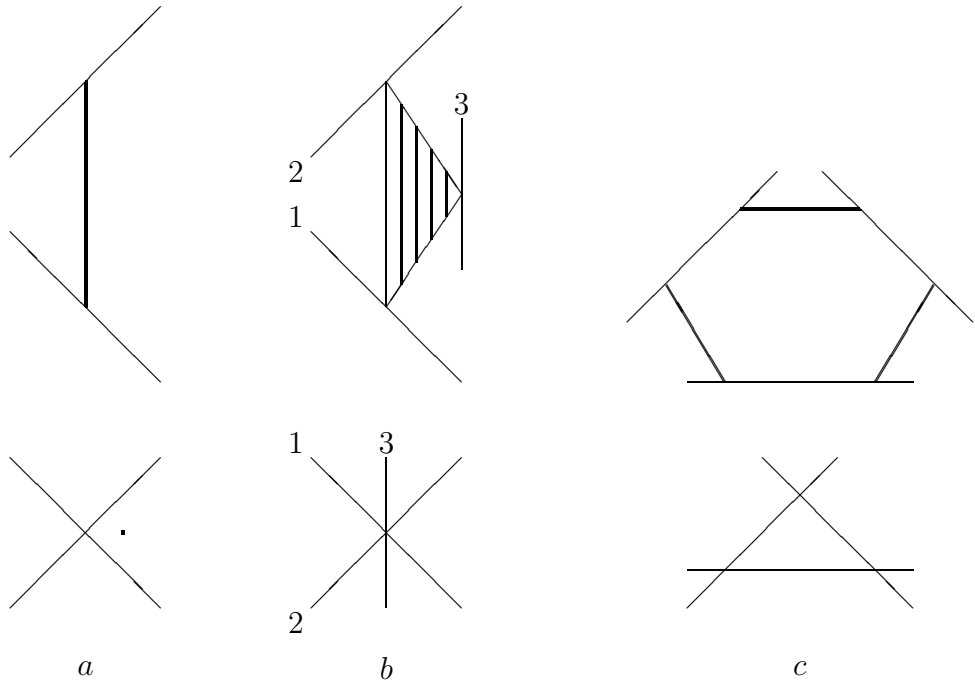
\begin{figure}
\unitlength 1.00mm \linethickness{0.4pt}
\begin{picture}(128.00,90.00)
\put(0.00,10.00){\line(1,1){20.00}} \put(0.00,30.00){\line(1,-1){20.00}}
\put(20.00,40.00){\line(-1,1){20.00}} \put(0.00,70.00){\line(1,1){20.00}}
\put(10.10,80.10){\line(0,-1){30.20}} \put(9.90,79.90){\line(0,-1){29.80}}
\put(60.00,10.00){\line(-1,1){20.00}} \put(60.00,30.00){\line(-1,-1){20.00}}
\put(50.00,10.00){\line(0,1){20.00}} \put(60.00,40.00){\line(-1,1){20.00}}
\put(40.00,70.00){\line(1,1){20.00}} \put(60.00,55.00){\line(0,1){20.00}}
\put(50.00,50.00){\line(0,1){30.00}} \put(50.00,80.00){\line(2,-3){10.00}}
\put(60.00,65.00){\line(-2,-3){10.00}}
\put(60.00,77.00){\makebox(0,0)[cc]{$3$}}
\put(38.00,68.00){\makebox(0,0)[cc]{$2$}}
\put(38.00,62.00){\makebox(0,0)[cc]{$1$}}
\put(38.00,32.00){\makebox(0,0)[cc]{$1$}}
\put(38.00,8.00){\makebox(0,0)[cc]{$2$}}
\put(50.00,32.00){\makebox(0,0)[cc]{$3$}} \put(90.00,10.00){\line(1,1){20.00}}
\put(100.00,30.00){\line(1,-1){20.00}} \put(120.00,15.00){\line(-1,0){30.00}}
\put(82.00,48.00){\line(1,1){20.00}} \put(108.00,68.00){\line(1,-1){20.00}}
\put(120.00,40.00){\line(-1,0){30.00}} \put(114.90,40.00){\line(3,5){7.80}}
\put(115.10,40.00){\line(3,5){7.70}} \put(95.10,40.00){\line(-3,5){7.80}}
\put(94.90,40.00){\line(-3,5){7.70}} \put(97.10,63.10){\line(1,0){15.80}}
\put(96.90,62.90){\line(1,0){16.20}}
\put(14.75,19.75){\rule{0.50\unitlength}{0.50\unitlength}}
\put(10.00,2.00){\makebox(0,0)[cc]{$a$}}
\put(50.00,2.00){\makebox(0,0)[cc]{$b$}}
\put(105.00,2.00){\makebox(0,0)[cc]{$c$}} \put(52.00,53.00){\line(0,1){24.00}}
\put(54.00,74.00){\line(0,-1){18.00}} \put(56.00,59.00){\line(0,1){12.00}}
\put(58.00,68.00){\line(0,-1){6.00}}
\end{picture}
\caption{Examples of line arrangements} \label{arr}
\end{figure}

For some three line arrangements in $\R^2$ (shown in the lower part of
Fig.~\ref{arr}) the corresponding simplicial resolutions are given above them
in the same picture. These resolutions are constructed as follows.

First, we embed the set of indices $\{1, \ldots, k\}$ into a space $\R^T$ of
dimension $T \ge k-1$ in such a way that their convex hull is a
$(k-1)$-dimensional simplex. The resolution will be constructed as a subset in
$\R^T \times \R^N.$ For any point $x \in L$ denote by $\tilde \Delta(x)$ the
convex hull in $\R^T$ of images of such indices $i$ that $L_i \ni x$, i.e. the
simplex with vertices at images of all these indices. Denote by $\Delta(x)$ the
simplex $\tilde \Delta(x) \times \{x\} \subset \R^T \times \R^N.$ Denote by
$L'$ the union of all simplices $\Delta(x)$, $x \in L.$ The obvious projection
$L' \to L$ (sending any $\Delta(x)$ to $x$) is a homotopy equivalence, as well
as its extension to the map of one-point compactifications $\bar L' \to \bar
L.$ In particular $\bar H_*(L') \equiv \bar H_*(L).$

On the other hand, $L'$ has a very useful filtration. For any set of indices $I
\subset \{1, \ldots, k\}$, denote by $L_I$ the plane $\cap_{i \in I} L_i.  $
Let $L'_I \subset L'$ be the {\em proper preimage} of $L_I$, i.e. the closure
of the union of complete preimages of all {\em generic} points of $L_I$ (i.e.
of points not from even smaller strata $L_J \subset L_I,$ $L_J \ne L_I$). There
is obvious homeomorphism $L'_I \simeq \tilde \Delta(I) \times L_I,$ where
$\tilde \Delta(I) \subset \R^T$ is the simplex whose vertices correspond to all
indices $i$ such that $L_i \supset L_I.$ (It is equal to $\tilde \Delta(x)$
where $x$ is any generic point of $L_I.$)

By definition, $L' = \bigcup L'_I,$ where the union is taken over all {\em
geometrically different} planes $L_I$. We define the term $F_p$ of the desired
filtration of $L'$ as the similar union of prisms $L'_I$ over all planes $L_I$
{\em of codimension $\le p$}. Then we extend it to a filtration on the
one-point compactification $\overline{L'}$ of $L'$ setting $F_0=$ \{the added
point\}.

This filtration defines a spectral sequence calculating the group $\bar H_*(L')
\simeq \bar H_*(L)$: by definition its term $E_{p,q}^1$ is equal to $\bar
H_{p+q}(F_p \setminus F_{p-1}) \equiv
H_{p+q}(\overline{F_p}/\overline{F_{p-1}}).$ This space $F_p \setminus F_{p-1}$
splits into a disjoint union (over all spaces $L_I$ of codimension exactly $p$)
of spaces $\check L'_I \stackrel{def}{=} \check \Delta(I) \times L'_I$, where
$\check \Delta(I)$ is the simplex $\tilde \Delta(I)$ from which several faces
are removed: namely such faces $\tilde \Delta(J)$, $J \subset I,$ that the
plane $L_J$ is strictly greater than $L_I$. For instance, for the
configurations shown in pictures a), b), c) of Fig.~\ref{arr} the planes $L_I$
of codimension 2 are: the point $(1,2)$, the point $(1,2,3)$, and three points
$(1,2),$ $(1,3)$, $(2,3)$ respectively. The proper preimages of them are: one
segment, one triangle (shadowed vertically in the picture), and three segments.
In all these cases the corresponding spaces $\check L_I$ coincide with $\check
\Delta_I$, namely they are: an open interval, a triangle without vertices, and
three open intervals, respectively.

In general, any face of the simplex $\tilde \Delta(I)$ is characterized by its
vertices, i.e. some indices $i \in \{1, \ldots, k\}$. A face of $\tilde
\Delta(I)$ is called {\em marginal}\label{marginal} if the intersection of
planes $L_i$ labeled by its vertices is strictly greater than $L_I.$ $\check
\Delta(I)$ is equal to $\tilde \Delta(I)$ with all marginal faces removed. By
the K\"unneth formula, $E^1_{p,q}= \bigoplus \bar H_{p+q-(N-p)}(\check
\Delta(I)),$ summation over all planes $L_I$ of codimension $p$.

The geometrical sense of the corresponding filtration  in the Alexander dual
group $H^*(\R^N \setminus L)$ is as follows: any element of this group has
filtration $p$ if and only if it is equal to a linear combination of finitely
many elements $\gamma_j$, any of which can be represented by the intersection
index with some semilinear\footnote{= semialgebraic distinguished by only
linear equations and inequalities} subvariety $V_j \subset \R^N,$ $\partial V_j
\subset L,$ invariant under the group $\R^{N-p_j}$ of translations in all
directions parallel to some $(N-p_j)$-dimensional plane $L_{I_j}$ with $p_j \le
p$.

\begin{proposition}[see \cite{congr}]
Our filtration of the space $\overline{L'}$ always homotopically splits, i.e.
we have the homotopy equivalence
\begin{equation}
\bar L' \sim \bar F_1 \vee (\bar F_2/\bar F_1) \vee \ldots \vee (\bar F_N/ \bar
F_{N-1}). \label{splitt}
\end{equation}
In particular, the spectral sequence degenerates in the first term: $E^\infty
\equiv E^1$, and we have $\bar H_{p+q}(\bar L') = \oplus_{p=1}^N E_{p,q}^1$.
\quad $\square$
\end{proposition}

An equivalent statement was proved in \cite{ZZ}.
\medskip

This theorem reduces the structure of cohomology groups of $\R^N \setminus L$
to dimensions of all spaces $L_I$. However, it does not allow us to calculate
the value of an arbitrary element of the group $E_{p,q}^1$ on any cycle in
$\R^N \setminus L$. For instance, in the case of the arrangement shown in
Fig.~\ref{arr}a, the entire group $E_{2,*}^1$ appears from the unique crossing
point $L_{\{1,2\}}$. This group is nontrivial only for $*=-1$, is isomorphic to
$\Z$ and is generated by the homology class of the segment $\Delta(1,2)$ modulo
its endpoints (lying in $F_1$). The splitting formula (\ref{splitt}) means that
we {\em can} extend this relative cycle of $\bar F_2 \ (\mbox{mod} \ \bar F_1)$
(or, equivalently, a closed locally finite cycle in $F_2 \setminus F_1$) to a
cycle in entire $\bar L'$ (respectively, in entire $L'$). However, to be able
to define the value of this point or of this segment on any $0$-dimensional
cycle (i.e. on a point) in $\R^2 \setminus L$ we need to {\em choose} such an
extension explicitly. Then we project it to $L$ and get a cycle there. Finally,
we need to choose a relative cycle in $\R^2 \ (\mbox{mod} \ L)$ whose boundary
coincides with this cycle. Then we call this relative cycle "a combinatorial
formula": its value on a point in $\R^2 \sm L$ is equal to the multiplicity of
this cycle in the neighborhood of this point.

If we have a more complicated plane arrangement, then the most convenient way
to extend an element of $E^1_{p,q}$ to a closed cycle in $L'$ is to do it step
by step over our filtration. Our starting element $\gamma \in E_{p,q}^1$ is
represented by a cycle with closed supports in $F_p \setminus F_{p-1}$. We take
its {\em first boundary} $d_1(\gamma)$, which is a cycle in $F_{p-1} \
(\mbox{mod} \ F_{p-2})$. Then we {\em span} it, i.e. construct a closed chain
$\tilde \gamma_1 $ in $ F_{p-1} \setminus F_{p-2}$ such that $\partial \tilde
\gamma_1 = d_1 \gamma $ there. Then we take the boundary of $\gamma+\tilde
\gamma_1$ in the space $F_{p-2} \setminus F_{p-3}$ and span it there by a chain
$\tilde \gamma_2,$ etc. The degeneration formula (\ref{splitt}) ensures that
all this sequence of choices can be accomplished. See \cite{ZZ}, \cite{M} for
some precise algorithms of doing it in the case of plane arrangements.

\subsection{All the same for knots}

The case of knots (say, of long knots) is very similar to that of plane
arrangements. A list of parallel notions is given in Table \ref{glos} (whose
left part was explained in the previous subsection, and the right-hand part
will be explained in the present one).

\begin{table}
\begin{tabular}{|l|l|}
\hline
Space $\R^N$ & Space ${\mathcal K}_n$ of smooth maps $\R^1 \to \R^n$ \\
& with a fixed behavior at $\infty$ \\
\hline Union of planes $L=\cup L_i \subset \R^N$ & Discriminant subset
$\Sigma \subset {\mathcal K}_n$ \\
\hline
Set of indices $\{1, \ldots, k\}$ & Chord space $\overline{B(\R^1,2)}$ \\
\hline
A plane $L_i$ & A subspace $L(a,b)$, $a,b \in \R^1$ \\
\hline Disjoint union of hyperplanes $L_i$ & Tautological resolution $F_1
\sigma$
of $\Sigma$ \\
\hline Simplicial resolution $L'$ of $L$ & Simplicial resolution
$\sigma$ of $\Sigma$ \\
\hline Subsets $I \subset \{1, \ldots, k\}$ & Combinatorial types of chord
configurations $J$ \\
with codim$L_I =p$ & with codim$L(J) = pn$ \\
\hline
A prism $L'_I$ & A $J$-block in $\sigma$ \\
\hline K\"unneth isomorphism for &
Thom isomorphism for the fibration of \\
homology of $\check L'_I = \check \Delta(I) \times L_I$ &
pure $J$-blocks by spaces $L(J')$ \\
\hline
Homotopy splitting (\ref{splitt}) & Kontsevich's degeneration theorem \\
\hline
\end{tabular}
\caption{} \label{glos}
\end{table}

So, instead of $\R^N$ we consider the affine space ${\mathcal K}_n$ of all
smooth maps $\R^1 \to \R^n$ coinciding with a fixed linear embedding "at
infinity", and instead of $L$ the discriminant space $\Sigma \subset {\mathcal
K}_n$ of all such maps which are not smooth embeddings.

Of course, the space ${\mathcal K}_n$ is infinite-dimensional, and formally we
cannot use the Alexander duality in it: the (finite-dimensional) cohomology
classes of the space of knots ${\mathcal K}_n \setminus \Sigma$ should
correspond to "infinite-dimensional cycles" in $\Sigma$, whose definition
requires some effort. The strict definition of such cycles corresponding to
finite-type cohomology classes was proposed in \cite{V1} and is as follows. We
consider a sequence of finite-dimensional approximating subspaces ${\mathcal
K}_n^j$ in ${\mathcal K}_n$, calculate (some) cohomology classes of ${\mathcal
K}^j_n \setminus \Sigma$ dual to certain cycles in $\Sigma$, and then prove a
stabilization theorem when $j \to \infty.$ It follows from the Weierstrass
approximation theorem that these stable cocycles are well-defined cohomology
classes in ${\mathcal K}_n \setminus \Sigma$. A rigorous reader can either read
\cite{V1} or \cite{fasis} for all justifications or to think of the spaces
${\mathcal K}_n$ as of such approximating spaces of very high but finite
dimension. Let us denote this virtual dimension of ${\mathcal K}_n$ by
$\omega.$

Again, $\Sigma$ is the union of a family of subspaces of very simple nature.
For any pair of points $(a,b)$ in $\R^1$, denote by $L(a,b)$ the space of all
maps $f \in {\mathcal K}_n$ such that
\begin{equation}
\label{cond} f(a)=f(b) \mbox{ (if }a\ne b\mbox{ ) or }f'(a)=0\mbox{ (if }a=b).
\end{equation}
Such spaces form a 2-parametric family parametrized by all points $(a,b)$ of
the space $\overline{B(\R^1,2)}$ of all unordered collections of two points in
$\R^1.$ Since \cite{V1} such pairs are depicted by arcs connecting the points
$a,b$ (called {\em chords} in \cite{BN}), so the space $\overline{B(\R^1,2)}$
will be called here the {\em chord space}. Its degenerated points
(corresponding to pairs $a=b$) are depicted by an asterisk at the point $a$.

The {\em tautological resolution} $F_1\sigma$ of $\Sigma$ is constructed as a
subspace of the direct product $\overline{B(\R^1,2)} \times {\mathcal K}_n $:
this is the space of pairs $((a,b),f)$ satisfying (\ref{cond}). It is the space
of an $(\omega-n)$-dimensional vector bundle over $\overline{B(\R^1,2)}.$
Therefore by the Thom isomorphism we have $\bar H_*(F_1\sigma)\simeq \bar
H_{*-(\omega-n)}(\overline{B(\R^1,2)}) \equiv 0:$ indeed,
$\overline{B(\R^1,2)}$ is homeomorphic to the closed half-plane. There is
obvious projection $F_1 \sigma \to \Sigma$; it is a map onto, and close to
generic points of $\Sigma$ is a homeomorphism.

Further, we insert simplices spanning preimages of non-generic points of
$\Sigma$. As previously, we embed the space $\overline{B(\R^1,2)}$ generically
and algebraically into a space $\R^T$ of a huge dimension ($T \approx \omega^3
$). Then for any point $f \in \Sigma$ we mark all the points $(a,b) \in
\overline{B(\R^1,2)}$ such that $L(a,b)\ni f,$ and denote by $\tilde \Delta(f)$
the convex hull of images of all these points in $\R^T.$

Of course, there exist points $f \in \Sigma$ having infinitely many preimages.
However they form a subset of infinite codimension in ${\mathcal K}_n$, and we
can ignore them by considering only finite-dimensional approximations
${\mathcal K}_n^j$ in general position with the stratification of $\Sigma$.
Then all the sets $\tilde \Delta(f),$ $f \in {\mathcal K}_n^j,$ still will be
the simplices with vertices at the images of all corresponding points $(a,b)$
of the chord space. The simplicial resolution $\sigma \subset \R^T \times
{\mathcal K}_n$ is defined as the union of all simplices $\Delta(f) \equiv
\tilde \Delta(f) \times \{f\}$.

Again, $\sigma$ has a useful increasing filtration. Let $I\subset
\overline{B(\R^1,2)}$ be a finite set of chords $(a,b)$. The intersection of
corresponding planes $L(a,b)$ is a subspace $L(I) \subset {\mathcal K}_n,$
whose codimension is a multiple of $n$. The proportionality coefficient
$\mbox{codim}L(I)/n$ is called the {\em complexity} of $I$. Consider all the
points $(a,b) \in \overline{B(\R^1,2)}$ such that $L(a,b) \supset L(I),$ and
denote by $\tilde \Delta(I) \subset \R^T$ the convex hull of images of all
these points. (It is equal to the space $\tilde \Delta(f)$ where $f$ is a {\em
generic} point of the space $L(I)$.) Set $L'(I)=\tilde \Delta(I) \times L(I)
\subset \R^T \times {\mathcal K}_n.$ Finally, define the term $F_p(\sigma)$ of
the filtration as the union of all simplices $\Delta(I)$ over all $I$ of
complexity $\le p.$

\begin{definition}
\label{fitype} {\rm A cohomology class of the space of knots ${\mathcal K}_n
\setminus \Sigma$ is a {\it finite type class of order} $p$ if it can be
defined as the linking number with the direct image in $\Sigma$ of a cycle
(with closed support) lying in the term $F_p$ of the standard filtration of
$\sigma$. For any such class of order $p$ and dimension $d$, its {\em principal
part} is the class of the corresponding cycle in the group $\bar
H_{\omega-d-1}(F_p \setminus F_{p-1})$.}
\end{definition}

The important property of this filtration is as follows: any its term $F_p
\setminus F_{p-1}$ is the space of an $(\omega-pn)$-dimensional affine bundle
over some finite-dimensional semialgebraic base: the projection of this bundle
is induced by the obvious projection $\R^T \times {\mathcal K}_n \to \R^T.$ In
particular, the Thom isomorphism reduces the Borel--Moore homology group of
this term to the homology group of locally finite chains of this base (in the
case of odd $n$ with coefficients in some system of groups locally isomorphic
to $\Z$, which is constant only for $p=1$).

These finite-dimensional bases, and hence also entire spaces $F_p \setminus
F_{p-1}$  of our filtration, admit an easy description, in particular their
one-point compactifications have a natural structure of $CW$-complexes. First
let us describe all the spaces $L(I)$ of complexity exactly $p$.

\begin{definition}[see \cite{V1}]{\rm
Let $A$ is a unordered finite collection of naturals $A=(a_1, \ldots, a_{\#
A}),$  $a_j \ge 2,$ and $b$  any nonnegative integer. Then an $(A,b)$-{\em
configuration} in $\R^1$ is any collection of distinct $a_1+ \cdots + a_{\# A}$
points in $\R^1$ separated into groups of cardinalities $a_1, \ldots, a_{\#
A}$, plus a collection of $b$ distinct points in  $\R^1$  (some of which can
coincide with the points of the $A$-part). A map $f:\R^1 \to \R^n$ {\em
respects} an $(A,b)$-configuration $J$ if it maps all points of any of groups
of cardinality $a_j,$ $j=1, \ldots, \# A,$ into one point (these points for
different groups may coincide), and $f'=0$ at all points of the $b$-part of the
configuration. The space of all maps $f$ respecting a fixed
$(A,b)$-configuration $J$ is denoted by $L(J)$. Two $(A,b)$-configurations are
{\em equivalent} if they can be transformed one into the other  by an
orientation-preserving homeomorphism of $\R^1$. The {\em complexity} of an
$(A,b)$-configuration is the number $\sum_{j=1}^{\# A}(a_j-1)+b.$}
\end{definition}

Obviously the codimension in ${\mathcal K}_n$ of any space $L(J)$ is equal to
$n$ times the complexity of $J$. The space of all $(A,b)$-configurations of
complexity $1$ is the chord space $\overline{B(\R^1,2)}$. Any intersection of
finitely many planes $L(a,b)$, $(a,b) \in\overline{B(\R^1,2)},$ is a plane of
form $L(J)$ for some $(A,b)$-configuration $J$. The corresponding simplex
$\tilde \Delta(J)$ in $\R^T$ has exactly $\sum_{i=1}^{\# A}\binom{a_i}{2} +b$
vertices. The $J$-{\em block} $\square(J)$ in $\sigma$ is the union of all
pairs $(x,f) \subset \R^T \times {\mathcal K}_n,$ such that $x$ belongs to the
simplex $\tilde \Delta(J')$ for some $(A,b)$-configuration $J'$ equivalent to
$J$, and $f$ respects this configuration $J'$. It belongs to the term $F_p$ of
our filtration, where $p$ is the complexity of $J$.

The {\em pure $J$-block} $\Check {\square}(J)$ is equal to $\square(J)
\setminus F_{p-1}$. It is fibered over the space of $(A,b)$-configurations $J'$
equivalent to $J$, with fiber equal to $\check \Delta(J) \times L(J)$, where
$\check \Delta(J)$ is the union of several (non-marginal in some sense) faces
of $\tilde \Delta(J)$. The base of this fiber bundle is an open cell, thus the
bundle is trivial, and we have a canonical decomposition of $\Check
{\square}(J)$ into open cells corresponding to all such non0marginal faces.

The canonical notation of any such cell is a {\em generalized chord diagram},
i.e. a finite collection of arcs connecting some points of $\R^1$ and of
asterisks marking some points, say as in the picture
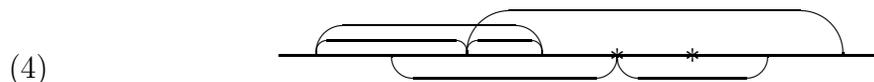
\begin{equation}
\label{samcell} \unitlength=1.00mm \special{em:linewidth 0.4pt}
\linethickness{0.4pt}
\begin{picture}(80.00,9.00)
\put(15.00,3.00){\oval(20.00,4.00)[t]} \put(30.00,3.00){\oval(10.00,4.00)[t]}
\put(20.00,3.00){\oval(30.00,8.00)[t]} \put(30.00,3.00){\oval(30.00,6.00)[b]}
\put(55.00,3.00){\oval(20.00,6.00)[b]} \put(55.00,3.00){\makebox(0,0)[cc]{$*$}}
\put(45.00,3.00){\makebox(0,0)[cc]{$*$}}
\put(50.00,3.00){\oval(50.00,12.00)[t]} \put(80.00,3.00){\line(-1,0){80.00}}
\end{picture}
\end{equation}
presenting one of cells of a certain equivalence class of
$((4,3),2)$-configurations. Namely, for any such cell related with a class of
equivalent $(A,b)$-configurations, we fix some configuration $J \subset \R^1$
of this class, mark by asterisks all points of its $b$-part ("singular points")
and draw a chord between any its two points $a,b$ such that the point $(a,b)
\in \overline{B(\R^1,2)}$ is a vertex of the face of $\tilde \Delta(J)$
corresponding to this cell.

The space $F_p \setminus F_{p-1}$ is the union of such pure blocks $\Check
{\square}(J)$ over (finitely many) equivalence classes of all
$(A,b)$-configurations of complexity exactly $p$. So we get also the
decomposition of this space into finitely many open cells. This decomposition
can be extended to the structure of a $CW$-complex on the one-point
compactification of $F_p \setminus F_{p-1}$. Its structure and incidence
coefficients are explicitly described in \cite{V1}, \cite{fasis}, which gives
also an algorithm for calculating the term $E^1$  of the spectral sequence
generated by this filtration and converging to the group of all finite-type
cohomology classes of the space of knots. In particular, if $n=3$ then all
$J$-blocks of complexity $p$ which (by dimensional reasons) can be valuable for
the calculation of knot invariants, are only the blocks with $(A,b)$ equal to
$((2,\ldots,2),0)$ (chord diagrams), $((2,\ldots,2),1)$ ({\em one-term
relations}, see \cite{BL}, \cite{BN}), and $((3,2,\ldots,2),0)$ (any such block
corresponds to the totality of {\em 4-term relations} arising from the
neighborhood of a triple point: there are three such relations, any two of
which are independent).

\begin{example}
\label{exf1} {\rm The term $F_1$ consists of exactly two cells, one of which is
the boundary of the other:
\begin{equation}
\label{dif01} \unitlength 1mm \linethickness{0.4pt}
\begin{picture}(72.00,4.00)
\put(3.00,0.00){\makebox(0,0)[cc]{$\partial$}}
\put(6.00,0.00){\line(1,0){30.00}} \put(39.00,0.00){\makebox(0,0)[cc]{$=$}}
\put(42.00,0.00){\line(1,0){30.00}} \put(21.00,0.00){\oval(20.00,8.00)[t]}
\put(57.00,0.00){\makebox(0,0)[cc]{{\large $*$}}}
\end{picture} \ ,
\end{equation}
thus there are no cohomology classes of order 1 of the space of long knots.
(Here and in the next example we assume some natural orientations of such
cells, see \cite{V1}, \cite{fasis}.)}
\end{example}

\begin{example}
\label{exf2} {\rm The term $F_2 \setminus F_1$ consists of the following cells:
four cells of maximal dimension
\begin{equation}
\unitlength 1.00mm \special{em:linewidth 0.4pt} \linethickness{0.4pt}
\begin{picture}(16.00,5.00)
\put(0.00,2.00){\line(1,0){16.00}} \put(6.50,2.00){\oval(9.00,4.00)[b]}
\put(9.50,2.00){\oval(9.00,4.00)[t]}
\end{picture}
 \ , \
\unitlength 1.00mm \linethickness{0.4pt}
\begin{picture}(39.00,4.00)
\put(0.00,2.00){\line(1,0){16.00}} \put(3.50,2.00){\oval(5.00,4.00)[t]}
\put(12.50,2.00){\oval(5.00,4.00)[t]} \put(19.00,0.00){\makebox(0,0)[cc]{,}}
\put(23.00,2.00){\line(1,0){16.00}} \put(31.00,2.00){\oval(14.00,4.00)[t]}
\put(31.00,2.00){\oval(6.00,4.00)[b]}
\end{picture}
 \ , \
\unitlength 1.00mm \special{em:linewidth 0.4pt} \linethickness{0.4pt}
\begin{picture}(18.00,4.00)
\put(0.00,0.00){\line(1,0){16.00}} \put(8.00,0.00){\oval(14.00,8.00)[t]}
\put(4.50,0.00){\oval(7.00,4.00)[t]} \put(11.50,0.00){\oval(7.00,4.00)[t]}
\end{picture}
\label{maxcel}
\end{equation}
(only the first and the last of which will be interesting for us), three cells
forming the boundary of any of these two cells,
\begin{equation}
\label{ztripno} \unitlength 1.00mm \special{em:linewidth 0.4pt}
\linethickness{0.4pt}
\begin{picture}(20.00,4.00)
\put(0.00,0.00){\line(1,0){20.00}} \put(10.00,0.00){\oval(16.00,8.00)[t]}
\put(14.00,0.00){\oval(8.00,4.00)[t]}
\end{picture}
 \ , \quad
\unitlength 1.00mm \special{em:linewidth 0.4pt} \linethickness{0.4pt}
\begin{picture}(20.00,4.00)
\put(0.00,0.00){\line(1,0){20.00}} \put(14.00,0.00){\oval(8.00,4.00)[t]}
\put(6.00,0.00){\oval(8.00,4.00)[t]}
\end{picture}
 \ , \quad
\mbox{and} \quad \unitlength 1.00mm \special{em:linewidth 0.4pt}
\linethickness{0.4pt}
\begin{picture}(20.00,4.00)
\put(0.00,0.00){\line(1,0){20.00}} \put(10.00,0.00){\oval(16.00,8.00)[t]}
\put(6.00,0.00){\oval(8.00,4.00)[t]}
\end{picture} \ ,
\end{equation}
and also six cells defined by maps with singular points (i.e. in whose notation
the asterisks $*$ participate):
\begin{equation}
\label{marg} \unitlength 1.00mm \special{em:linewidth 0.4pt}
\linethickness{0.4pt}
\begin{picture}(119.00,3.00)
\put(1.00,1.00){\line(1,0){15.00}} \put(21.00,1.00){\line(1,0){15.00}}
\put(41.00,1.00){\line(1,0){15.00}} \put(61.00,1.00){\line(1,0){15.00}}
\put(81.00,1.00){\line(1,0){15.00}} \put(101.00,1.00){\line(1,0){15.00}}
\put(10.50,1.00){\oval(7.00,4.00)[t]} \put(4.00,1.00){\makebox(0,0)[cc]{$*$}}
\put(28.50,1.00){\oval(11.00,4.00)[t]} \put(29.00,1.00){\makebox(0,0)[cc]{$*$}}
\put(46.50,1.00){\oval(7.00,4.00)[t]} \put(54.00,1.00){\makebox(0,0)[cc]{$*$}}
\put(68.50,1.00){\oval(9.00,4.00)[t]} \put(64.00,1.00){\makebox(0,0)[cc]{$*$}}
\put(88.50,1.00){\oval(9.00,4.00)[t]} \put(93.00,1.00){\makebox(0,0)[cc]{$*$}}
\put(104.00,1.00){\makebox(0,0)[cc]{$*$}}
\put(113.00,1.00){\makebox(0,0)[cc]{$*$}}
\put(18.00,1.00){\makebox(0,0)[cc]{,}} \put(38.00,1.00){\makebox(0,0)[cc]{,}}
\put(58.00,1.00){\makebox(0,0)[cc]{,}} \put(78.00,1.00){\makebox(0,0)[cc]{,}}
\put(98.00,1.00){\makebox(0,0)[cc]{,}} \put(119.00,1.00){\makebox(0,0)[cc]{.}}
\end{picture}
\end{equation}

The boundary operator in this term $F_2 \setminus F_1$ (i.e. the vertical
differential $d^0$ of the spectral sequence) acts as follows:
$$
\unitlength 1.00mm \linethickness{0.4pt}
\begin{picture}(125.00,53.00)
\put(3.00,34.00){\makebox(0,0)[cc]{$\partial$}}
\put(7.00,34.00){\line(1,0){22.00}} \put(33.00,34.00){\makebox(0,0)[cc]{$=$}}
\put(45.50,34.00){\line(1,0){21.00}} \put(125.00,34.00){\line(-1,0){22.00}}
\put(96.00,34.00){\line(-1,0){22.00}}
\put(3.00,49.00){\makebox(0,0)[cc]{$\partial$}}
\put(7.00,49.00){\line(1,0){22.00}} \put(33.00,49.00){\makebox(0,0)[cc]{$=$}}
\put(40.00,49.00){\line(1,0){22.00}} \put(125.00,49.00){\line(-1,0){22.00}}
\put(94.00,49.00){\line(-1,0){22.00}}
\put(3.00,4.00){\makebox(0,0)[cc]{$\partial$}}
\put(7.00,4.00){\line(1,0){22.00}} \put(33.00,4.00){\makebox(0,0)[cc]{$=$}}
\put(40.00,4.00){\line(1,0){22.00}} \put(125.00,4.00){\line(-1,0){22.00}}
\put(94.00,4.00){\line(-1,0){22.00}}
\put(3.00,19.00){\makebox(0,0)[cc]{$\partial$}}
\put(7.00,19.00){\line(1,0){22.00}} \put(33.00,19.00){\makebox(0,0)[cc]{$=$}}
\put(45.50,19.00){\line(1,0){21.00}} \put(125.00,19.00){\line(-1,0){22.00}}
\put(96.00,19.00){\line(-1,0){22.00}} \put(15.00,49.00){\oval(12.00,8.00)[t]}
\put(21.00,49.00){\oval(12.00,8.00)[b]}
\put(37.00,49.00){\makebox(0,0)[cc]{$-$}}
\put(50.50,49.00){\oval(15.00,8.00)[t]} \put(47.00,49.00){\oval(8.00,8.00)[b]}
\put(67.00,49.00){\makebox(0,0)[cc]{$+$}}
\put(78.50,49.00){\oval(9.00,8.00)[t]} \put(87.50,49.00){\oval(9.00,8.00)[t]}
\put(99.00,49.00){\makebox(0,0)[cc]{$-$}}
\put(114.00,49.00){\oval(18.00,8.00)[t]}
\put(118.50,49.00){\oval(9.00,8.00)[b]} \put(12.00,34.00){\oval(6.00,8.00)[t]}
\put(24.00,34.00){\oval(6.00,8.00)[t]}
\put(40.00,34.00){\makebox(0,0)[cc]{$(-1)^n$}}
\put(48.00,34.00){\makebox(0,0)[cc]{{\large $*$}}}
\put(60.00,34.00){\oval(8.00,8.00)[t]}
\put(70.00,34.00){\makebox(0,0)[cc]{$+$}}
\put(80.50,34.00){\oval(9.00,8.00)[t]} \put(89.50,34.00){\oval(9.00,8.00)[t]}
\put(100.00,34.00){\makebox(0,0)[cc]{$-$}}
\put(109.50,34.00){\oval(9.00,8.00)[t]}
\put(123.00,34.00){\makebox(0,0)[cc]{{\large $*$}}}
\put(18.00,19.00){\oval(18.00,8.00)[t]} \put(18.00,19.00){\oval(6.00,8.00)[b]}
\put(40.00,19.00){\makebox(0,0)[cc]{$(-1)^n$}}
\put(55.50,19.00){\oval(15.00,8.00)[t]} \put(52.00,19.00){\oval(8.00,8.00)[b]}
\put(70.00,19.00){\makebox(0,0)[cc]{$+$}}
\put(85.00,19.00){\oval(18.00,8.00)[t]}
\put(85.00,19.00){\makebox(0,0)[cc]{{\large $*$}}}
\put(100.00,19.00){\makebox(0,0)[cc]{$-$}}
\put(114.00,19.00){\oval(18.00,8.00)[t]}
\put(118.50,19.00){\oval(9.00,8.00)[b]} \put(114.00,4.00){\oval(18.00,8.00)[t]}
\put(118.50,4.00){\oval(9.00,8.00)[b]} \put(99.00,4.00){\makebox(0,0)[cc]{$+$}}
\put(78.50,4.00){\oval(9.00,8.00)[t]} \put(87.50,4.00){\oval(9.00,8.00)[t]}
\put(67.00,4.00){\makebox(0,0)[cc]{$-$}} \put(50.50,4.00){\oval(15.00,8.00)[t]}
\put(47.00,4.00){\oval(8.00,8.00)[b]} \put(18.00,4.00){\oval(18.00,8.00)[t]}
\put(13.50,4.00){\oval(9.00,4.00)[t]} \put(22.50,4.00){\oval(9.00,4.00)[t]}
\end{picture}
$$
$$
\unitlength 1.00mm \linethickness{0.4pt}
\begin{picture}(105.00,38.00)
\put(3.00,34.00){\makebox(0,0)[cc]{$\partial$}}
\put(6.00,34.00){\line(1,0){25.00}} \put(34.00,34.00){\makebox(0,0)[cc]{$=$}}
\put(47.00,34.00){\line(1,0){25.00}} \put(59.50,34.00){\oval(19.00,8.00)[t]}
\put(3.00,19.00){\makebox(0,0)[cc]{$\partial$}}
\put(3.00,4.00){\makebox(0,0)[cc]{$\partial$}}
\put(6.00,19.00){\line(1,0){25.00}} \put(6.00,4.00){\line(1,0){25.00}}
\put(34.00,19.00){\makebox(0,0)[cc]{$=$}}
\put(34.00,4.00){\makebox(0,0)[cc]{$=$}} \put(47.00,19.00){\line(1,0){25.00}}
\put(40.00,4.00){\line(1,0){25.00}} \put(59.50,19.00){\oval(19.00,8.00)[t]}
\put(52.50,4.00){\oval(19.00,8.00)[t]} \put(18.50,34.00){\oval(19.00,8.00)[t]}
\put(14.00,34.00){\oval(10.00,8.00)[b]}
\put(42.00,34.00){\makebox(0,0)[cc]{$(-1)^n$}}
\put(50.00,34.00){\makebox(0,0)[cc]{{\large $*$}}}
\put(92.50,19.00){\oval(19.00,8.00)[t]} \put(80.00,19.00){\line(1,0){25.00}}
\put(76.00,19.00){\makebox(0,0)[cc]{$+$}}
\put(42.00,19.00){\makebox(0,0)[cc]{$(-1)^n$}}
\put(50.00,19.00){\makebox(0,0)[cc]{{\large $*$}}}
\put(102.00,19.00){\makebox(0,0)[cc]{{\large $*$}}}
\put(62.00,4.00){\makebox(0,0)[cc]{{\large $*$}}}
\put(18.50,4.00){\oval(19.00,8.00)[t]} \put(23.50,4.00){\oval(9.00,8.00)[b]}
\put(23.50,19.00){\oval(9.00,8.00)[t]} \put(14.00,19.00){\oval(10.00,8.00)[t]}
\end{picture}
$$
$$
\unitlength 1.00mm \linethickness{0.4pt}
\begin{picture}(112.00,38.00)
\put(3.00,34.00){\makebox(0,0)[cc]{$\partial$}}
\put(6.00,34.00){\line(1,0){25.00}} \put(34.00,34.00){\makebox(0,0)[cc]{$=$}}
\put(42.00,34.00){\line(1,0){25.00}} \put(73.00,34.00){\line(1,0){25.00}}
\put(85.50,34.00){\oval(19.00,8.00)[t]}
\put(95.00,34.00){\makebox(0,0)[cc]{{\large $*$}}}
\put(3.00,19.00){\makebox(0,0)[cc]{$\partial$}}
\put(3.00,4.00){\makebox(0,0)[cc]{$\partial$}}
\put(6.00,19.00){\line(1,0){25.00}} \put(6.00,4.00){\line(1,0){25.00}}
\put(34.00,19.00){\makebox(0,0)[cc]{$=$}}
\put(34.00,4.00){\makebox(0,0)[cc]{$=$}} \put(42.00,19.00){\line(1,0){25.00}}
\put(42.00,4.00){\line(1,0){25.00}} \put(73.00,19.00){\line(1,0){25.00}}
\put(87.00,4.00){\line(1,0){25.00}} \put(85.50,19.00){\oval(19.00,8.00)[t]}
\put(95.00,19.00){\makebox(0,0)[cc]{{\large $*$}}}
\put(90.00,4.00){\makebox(0,0)[cc]{{\large $*$}}}
\put(109.00,4.00){\makebox(0,0)[cc]{{\large $*$}}}
\put(45.00,34.00){\makebox(0,0)[cc]{{\large $*$}}}
\put(64.00,34.00){\makebox(0,0)[cc]{{\large $*$}}}
\put(14.00,34.00){\oval(10.00,8.00)[t]} \put(18.50,19.00){\oval(19.00,8.00)[t]}
\put(28.00,34.00){\makebox(0,0)[cc]{{\large $*$}}}
\put(19.00,19.00){\makebox(0,0)[cc]{{\large $*$}}}
\put(9.00,4.00){\makebox(0,0)[cc]{{\large $*$}}}
\put(45.00,4.00){\makebox(0,0)[cc]{{\large $*$}}}
\put(54.50,4.00){\oval(19.00,8.00)[t]}
\put(45.00,19.00){\makebox(0,0)[cc]{{\large $*$}}}
\put(54.50,19.00){\oval(19.00,8.00)[t]}
\put(80.00,4.00){\makebox(0,0)[cc]{$(-1)^{n-1}$}}
\put(70.00,4.00){\makebox(0,0)[cc]{$+$}}
\put(70.00,19.00){\makebox(0,0)[cc]{$+$}}
\put(70.00,34.00){\makebox(0,0)[cc]{$+$}}
\put(39.00,34.00){\makebox(0,0)[cc]{$-$}}
\put(39.00,19.00){\makebox(0,0)[cc]{$-$}}
\put(39.00,4.00){\makebox(0,0)[cc]{$-$}} \put(23.50,4.00){\oval(9.00,8.00)[t]}
\end{picture} .
$$
In particular there is exactly one nontrivial group $\bar H_i(F_2 \setminus
F_1) \equiv E^1_{2,i-2}$, namely such group with $i=\omega-(2n-5)$ is
isomorphic to $\Z$ and is generated by the sum of the first and the last cells
in (\ref{maxcel}).}
\end{example}

Thus we obtain a proof of the first statement of Proposition \ref{longres}. The
group $\bar H_*(F_3 \setminus F_2)$ of (possible) principal parts of third
order cohomology classes was calculated in \cite{tetu} for odd $n$ (another
proof, not relying on the computer's honesty, see in \cite{bjo}) and in
\cite{bjo} for even $n$. In both cases, there are exactly two nontrivial groups
$\bar H_i(F_3 \setminus F_2) \simeq E^1_{3,i-3}$, namely with $i=\omega
-(3n-8)$ and $\omega-(3n-7)$; they both are isomorphic to $\Z$. By the
dimensional reasons both these groups for any $n \ge 3$ coincide with
corresponding groups $E^\infty_{3,i-3}$, and their generators extend to
well-defined cohomology classes of spaces of long knots in $\R^n$. By similar
considerations (see e.g. \cite{bjo}) for $n \ge 4$ these generators are
nontrivial and free. If $n=3$ then for the first of these classes the same
follows from the fact that it coincides with the second simple knot invariant
(calculated in \cite{V1}) whose nontriviality is well known. The other class is
exactly the Teiblum--Turchin class studied below; the fact that it also is
nontrivial and free for $n=3$ will follow from the proof of Corollary
\ref{nontriv}.

\begin{remark} {\rm
It is often convenient to replace formally our homological spectral sequence
calculating $\bar H_*(\sigma)$ by the  "Alexander dual" cohomological spectral
sequence
$$
E_r^{p,q} \equiv E^r_{-p,\omega-q-1}.
$$
It lies in the second quadrant in the wedge $\{(p,q)| p \le 0, q +pn \ge 0\}$
and converges to some subgroup of the group $H^*(\K_n \setminus \Sigma)$ (if
$n>3$ then to entire this group).}
\end{remark}

\begin{remark} {\rm
There are beautiful algebraic structures on the above-described spectral
sequence, and hence on its limit filtered group $H^*(\K_n \setminus \Sigma)$
and the corresponding adjoint graded group, see \cite{Turchin}.}
\end{remark}

All of this theory can be extended almost literally to the cohomology of the
space of compact knots $S^1 \hookrightarrow \R^n$.\footnote{As well as of
compact links, i.e. of embeddings of a disjoint union of several circles into
$\R^n$; we shall not discuss here the latter theory} However, in this case the
chord space $\overline{B(S^1,2)}$ is not topologically trivial (it is a closed
M\"obius band); also the spaces of equivalent $(A,b)$-configurations are not
the cells. To get the cell decomposition of all spaces $F_p \setminus F_{p-1}$
we need to mark one point in $S^1$ ("the origin") and call two configurations
equivalent if they are transformed one into the other by a homeomorphism of
$S^1$ preserving the origin and the orientation.
\medskip

The direct calculation of the spectral sequence and obtaining the {\em
combinatorial formulas}  for the finite-type cohomology classes of spaces of
knots is formally the same process as in the case of plane arrangements.
However, the exact choice of the spanning chains in all the consecutive terms
of the filtration and in entire ${\mathcal K}_n$ depends very much of the
features of the knot space.

\begin{remark}{\rm
There is another, sometimes more convenient construction of the resolution of
discriminant sets, namely the {\em conical resolutions} based on the notion of
a continuous order complex of a topologized partially ordered set, see e.g.
\cite{kotor}. In particular it allows us to resolve the points of $\Sigma$ with
infinitely many preimages in the tautological resolution space. However, for
the calculations in the present work it will be enough to use the "naive"
simplicial resolution described above.}
\end{remark}

\subsection{Finite type knot invariants and Polyak-Viro combinatorial
formulas}

Suppose that $n=3$ and we are interested in the knot invariants, i.e. the
0-dimensional cohomology classes of $\K_n \setminus \Sigma$. For any such class
of finite filtration $p$, its principal part in $F_p \setminus F_{p-1}$ is a
linear combination of cells depicted by $p$-chord diagrams (i.e. collections of
$p$ chords with distinct endpoints) and $\tilde p$-{\em configurations}, i.e.
collections of $p-2$ chords with different endpoints and one triple of points
joined by three chords. E.g., among all diagrams in
(\ref{samcell})--(\ref{marg}) only the left picture in (\ref{dif01}) and three
left pictures in (\ref{maxcel}) are chord diagrams, and only the last picture
in (\ref{maxcel}) is a $\tilde 2$-configuration. The coefficients with which
all these cells can enter the linear combination satisfy the homological
condition. In particular the coefficients at $\tilde p$-configurations are
determined by these at $p$-chord diagrams, and any admissible linear
combination is characterized uniquely only by the collection of latter
coefficients, which is called a {\em weight system}.

{\em The elementary characterization of these invariants} is as follows (see
e.g. \S 0.2 in \cite{V1}). Let us consider any immersion $\R^1 \to \R^3$ with
exactly $k$ transverse self-intersection points. We can resolve any of these
points in two locally distinct ways to get a knot without intersections. One of
these two local resolutions can be invariantly defined as a positive, and the
other as the negative one. The $k$-th index of a knot invariant at our singular
immersion is equal to the alternated sum of its values at all $2^k$ knots
obtained by all different possible resolutions of double points: the value at a
knot is counted with sign 1 or $-1$ depending on the parity of the number of
positive local resolutions. A knot invariant is of filtration $p$ if and only
if all its indices at all immersions with $k>p$ self-intersections are equal to
0. The same definition can be applied to define the filtration of invariants of
compact knots $S^1 \to \R^3$. On the other hand, it is easy to see that there
is a natural one-to-one correspondence between connected components of spaces
of long and compact knots, in particular the theories of their invariants
naturally coincide.

Some combinatorial formulas for the simplest finite-type knot invariants --- of
orders 2 and 3 --- were found in \cite{lannes}. Another, more convenient
formulas were introduced by M.~Polyak and O.~Viro in \cite{PV}. These formulas
for long knots look as the linear combinations of chord diagrams with oriented
chords. E.g. the formula \unitlength 1mm \linethickness{0.4pt}
\begin{picture}(22.00,6.00)
\put(0.00,3.00){\line(1,0){22.00}} \put(2.00,3.00){\vector(-4,-3){0.2}}
\bezier{68}(14.00,3.00)(8.00,9.00)(2.00,3.00)
\put(20.00,3.00){\vector(4,3){0.2}}
\bezier{68}(8.00,3.00)(14.00,-3.00)(20.00,3.00)
\end{picture} \ should be read as follows. Consider a generic
long knot $f : \R^1 \to \R^3$. A {\em representation} of the above picture in
this knot is any collection of points $a<b<c<d \subset \R^1$ such that $f(a)$
lies below $f(c)$ and $f(d)$ lies below $f(b)$. The value of this picture on
our knot is equal to the number of its representations (counted with
appropriate signs). An immediate calculation shows that this number is a knot
invariant of order 2.

In the case of compact knots, there are Polyak-Viro formulas of two types:
absolute and punctured ones. They also look as oriented chord diagrams, but
with endpoints in the oriented circle $S^1$ instead of $\R^1$; moreover, a
punctured Polyak-Viro diagram contains a point in $S^1$ not coinciding with the
endpoints of chords. E.g. a representation of the diagram \unitlength 1mm
\linethickness{0.4pt}
\begin{picture}(8.00,8.00)
\put(4.00,4.00){\circle{8.00}} \put(7.00,1.00){\vector(-1,1){6.00}}
\put(1.00,1.00){\vector(1,1){6.00}}
\end{picture}
 in a compact knot is any collection of four points
in $S^1$ with cyclic order $a<b<c<d<a$, satisfying the same conditions as
previously. A representation of the punctured diagram \unitlength 1.00mm
\linethickness{0.4pt}
\begin{picture}(8.00,9.00)
\put(4.00,4.00){\circle{8.00}} \put(7.00,1.00){\vector(-1,1){6.00}}
\put(1.00,1.00){\vector(1,1){6.00}} \put(4.00,8.00){\circle*{1.00}}
\end{picture}
 is such a collection of points in the parametrized
circle, whose cyclic coordinates satisfy a more strong condition
$0<a<b<c<d<2\pi;$ the origin $0=2\pi \in S^1$ corresponds to the marked point
in the diagram. It is easy to see that the number of representations of the
last punctured diagram (counted with natural signs) is a knot invariant, and
the similar number for the absolute diagram is not.

However, some of finite type knot invariants can be realized by absolute
diagrams: in particular the unique invariant of order 3 can be given by the
diagram \unitlength 1mm \linethickness{0.4pt}
\begin{picture}(29.00,8.00)
\put(9.00,4.00){\circle{8.00}} \put(25.00,4.00){\circle{8.00}}
\put(25.00,0.00){\vector(0,1){8.00}} \put(28.30,6.20){\vector(-3,-2){6.60}}
\put(21.70,6.20){\vector(3,-2){6.60}} \put(9.00,0.00){\vector(0,1){8.00}}
\put(5.60,6.00){\vector(1,0){6.80}} \put(12.40,2.00){\vector(-1,0){6.80}}
\put(2.00,4.00){\makebox(0,0)[cc]{$\frac{1}{2}$}}
\put(17.00,4.00){\makebox(0,0)[cc]{$+\frac{1}{3}$}}
\end{picture}
 , see \cite{PV}.

M.~Goussarov has proved that any finite-type knot invariant of long knots can
be represented by a formula of Polyak--Viro type, see \cite{GPV}.

Similar (but more complicated) formulas appear naturally in the calculation of
higher-dimensional cohomology classes of spaces of knots, see the next
sections.

\begin{remark}{\rm
The homological calculations discussed below provide numerous possibilities to
make a mistake: to miss some component of the boundary, to calculate wrongly
some orientation, etc. Fortunately, we always can check our calculations. If we
have calculated some boundary operator and suspect that it is not correct, just
calculate the boundary of this boundary, and try to understand why it is not
equal to zero! My experience says that no mistakes survive this examination.}
\end{remark}

\section{Proof of Theorem 1}
\label{proof1}

\subsection{Principal part of the cocycle}
\label{f3}

In the original calculation \cite{tetu}, the principal part of the
Teiblum-Turchin cocycle in the term $F_3 \setminus F_2$ of the natural
filtration of the resolved discriminant was found as a linear combination of
some 8 cells of the canonical cell decomposition, see e.g. \cite{bjo},
\cite{fasis}.

This expression can be simplified, especially if $n$ is even.

\begin{proposition}
For any $n\ge 3$, the group of order 3 cohomology classes of dimension $3n-8$
of the space of long knots $\R^1 \to \R^n$ is cyclic; for $n\ge 3$ it is free
Abelian.

If $n$ is even, then this group is generated by the sum of only two cells:
\begin{equation}
\label{prinpart} \special{em:linewidth 0.4pt} \unitlength 1.00mm
\linethickness{0.4pt}
\begin{picture}(68.00,7.00)
\put(68.00,3.00){\line(-1,0){23.00}} \put(54.00,3.00){\oval(10.00,4.00)[t]}
\put(61.50,3.00){\oval(5.00,4.00)[t]} \put(56.50,3.00){\oval(15.00,8.00)[t]}
\put(59.00,3.00){\oval(10.00,6.00)[b]} \put(41.00,3.00){\makebox(0,0)[cc]{$+$}}
\put(19.50,3.00){\oval(15.00,6.00)[b]} \put(27.00,3.00){\oval(10.00,4.00)[t]}
\put(24.50,3.00){\oval(15.00,8.00)[t]} \put(36.00,3.00){\line(-1,0){28.00}}
\put(2.00,3.00){\makebox(0,0)[cc]{$TT =$}}
\end{picture} .
\end{equation}

For odd $n$ it is generated by the linear combination

\begin{equation}
\label{prinpar} \unitlength=1.00mm \special{em:linewidth 0.4pt}
\linethickness{0.4pt}
\begin{picture}(125.00,7.00)
\put(18.00,3.00){\line(-1,0){18.00}} \put(6.00,3.00){\oval(10.00,6.00)[b]}
\put(8.50,3.00){\oval(5.00,4.00)[t]} \put(13.50,3.00){\oval(5.00,4.00)[t]}
\put(11.00,3.00){\oval(10.00,8.00)[t]} \put(21.00,3.00){\makebox(0,0)[cc]{$+$}}
\put(32.50,3.00){\oval(5.00,4.00)[t]} \put(37.50,3.00){\oval(5.00,4.00)[t]}
\put(35.00,3.00){\oval(10.00,8.00)[t]} \put(32.50,3.00){\oval(15.00,6.00)[b]}
\put(42.00,3.00){\line(-1,0){19.00}} \put(44.00,3.00){\makebox(0,0)[cc]{$+$}}
\put(70.00,3.00){\line(-1,0){24.00}} \put(55.50,3.00){\oval(15.00,6.00)[b]}
\put(63.00,3.00){\oval(10.00,4.00)[t]} \put(60.50,3.00){\oval(15.00,8.00)[t]}
\put(73.00,3.00){\makebox(0,0)[cc]{$+$}} \put(97.00,3.00){\line(-1,0){22.00}}
\put(86.00,3.00){\oval(10.00,6.00)[b]} \put(90.50,3.00){\oval(9.00,4.00)[t]}
\put(86.00,3.00){\oval(18.00,8.00)[t]}
\put(101.00,3.00){\makebox(0,0)[cc]{$-$}} \put(125.00,3.00){\line(-1,0){21.00}}
\put(112.50,3.00){\oval(15.00,8.00)[t]} \put(115.00,3.00){\oval(10.00,4.00)[t]}
\put(119.50,3.00){\oval(9.00,6.00)[b]}
\end{picture} .
\end{equation}
\end{proposition}

The first statement of this proposition for odd $n$ was essentially proved by
Teiblum and Turchin \cite{tetu}; the justification of entire statement see in
\S 6 of \cite{bjo} or \S V.8.8 of \cite{fasis}.

Since we consider our class mod 2, the stabilization formula (\ref{stabil})
allows us to use any of expressions (\ref{prinpart}), (\ref{prinpar}) in the
case of any $n$. We shall use the shorter "even" version (\ref{prinpart}).

All further calculations in this section are mod 2 only.

\subsection{On the pictures}

The system of notation in this work is an extension of that used in
\S~\ref{method} for the cells of the natural simplicial resolution of the
discriminant. Any of our pictures consists of a horizontal segment (the {\em
Wilson loop} symbolizing the line $\R^1$), several asterisks placed on it, and
several arcs ("chords") connecting some its points (these data determine such a
cell), plus some additional furniture consisting of broken lines ({\em
zigzags}) and subscripts, which distinguish certain subvarieties in these
cells.

For instance, the picture \unitlength 1.00mm \special{em:linewidth 0.4pt}
\linethickness{0.4pt}
\begin{picture}(21.00,5.00)
\put(8.50,2.00){\oval(13.00,4.00)[b]} \put(12.50,2.00){\oval(13.00,6.00)[t]}
\put(0.00,2.00){\line(1,0){21.00}} \put(19.00,2.00){\line(-2,1){4.00}}
\put(15.00,4.00){\line(-2,-1){4.00}}
\end{picture}
 \
means, first of all, that we are in the cell \unitlength 1.00mm
\special{em:linewidth 0.4pt} \linethickness{0.4pt}
\begin{picture}(16.00,5.00)
\put(0.00,2.00){\line(1,0){16.00}} \put(6.50,2.00){\oval(9.00,4.00)[b]}
\put(9.50,2.00){\oval(9.00,4.00)[t]}
\end{picture}
 \ of the term $F_2 \setminus F_1$.
This cell can be considered as the space of all triples $(\alpha,t,f)$ where
$\alpha$ is some quadruple of points $a<b<c<d$ in $\R^1$, $f$ is a smooth map
$\R^1 \to \R^n$ such that $f(a)=f(c),$ $f(b)=f(d)$, and $t$ is a point of the
segment $\tilde \Delta(J),$ $J=((a,c);(b,d))$, participating in the
construction of the resolution: its endpoints correspond formally to the pairs
of points $(a,c)$ and $(b,d)$ glued together by $f$. The additional zigzag in
the picture \unitlength 1.00mm \special{em:linewidth 0.4pt}
\linethickness{0.4pt}
\begin{picture}(21.00,5.00)
\put(8.50,2.00){\oval(13.00,4.00)[b]} \put(12.50,2.00){\oval(13.00,6.00)[t]}
\put(0.00,2.00){\line(1,0){21.00}} \put(19.00,2.00){\line(-2,1){4.00}}
\put(15.00,4.00){\line(-2,-1){4.00}}
\end{picture}
 \
distinguishes the subvariety in this cell, consisting of such triples
$(\alpha,t,f)$ that there exists one point $\lambda \in \R^1$ more, $b <
\lambda < c$, such that $f(\lambda)=f(d)$. By definition, this subvariety is
identical with the one encoded by the picture \unitlength 1.00mm
\special{em:linewidth 0.4pt} \linethickness{0.4pt}
\begin{picture}(20.00,5.00)
\put(8.50,2.00){\oval(13.00,4.00)[b]} \put(11.50,2.00){\oval(13.00,6.00)[t]}
\put(0.00,2.00){\line(1,0){20.00}} \put(5.00,2.00){\line(2,1){4.00}}
\put(9.00,4.00){\line(2,-1){4.00}}
\end{picture} .

The picture \unitlength 1.00mm \special{em:linewidth 0.4pt}
\linethickness{0.4pt}
\begin{picture}(21.00,5.00)
\put(8.50,2.00){\oval(13.00,4.00)[b]} \put(12.50,2.00){\oval(13.00,6.00)[t]}
\put(0.00,2.00){\line(1,0){21.00}} \put(19.00,2.00){\line(-2,1){4.00}}
\put(15.00,4.00){\vector(-2,-1){4.00}}
\end{picture}
(respectively, \unitlength 1.00mm \special{em:linewidth 0.4pt}
\linethickness{0.4pt}
\begin{picture}(21.00,5.00)
\put(8.50,2.00){\oval(13.00,4.00)[b]} \put(12.50,2.00){\oval(13.00,6.00)[t]}
\put(0.00,2.00){\line(1,0){21.00}} \put(15.00,4.00){\vector(2,-1){4.00}}
\put(15.00,4.00){\line(-2,-1){4.00}}
\end{picture} )
will denote almost the same, but with the condition $f(\lambda)=f(d)$ replaced
by the condition that $f(\lambda)$ has the same projection to $\R^{n-1}$ as
$f(d)=f(b)$ and lies {\em below} (respectively, {\em above}) $f(d)$ in the line
of all points with the same projection.

The subscript of type \unitlength 1.00mm \special{em:linewidth 0.4pt}
\linethickness{0.4pt}
\begin{picture}(8.00,5.00)
\put(4.00,0.00){\vector(1,2){2.50}} \put(4.00,0.00){\vector(-1,2){2.50}}
\put(0.33,2.00){\makebox(0,0)[cc]{\small 1}}
\put(7.67,2.00){\makebox(0,0)[cc]{\small 2}}
\end{picture}
 \
under a picture denotes the condition that the "vertical" direction in $\R^n$
lies in the angle between the tangent directions $f'(a_1)$ and $f'(a_2)$, where
$a_1$ and $a_2$ are the first and the second from the left points of $\R^1$
participating actively in the picture. Similarly, the subscript \unitlength
1.00mm \special{em:linewidth 0.4pt} \linethickness{0.4pt}
\begin{picture}(5.00,6.00)
\put(0.00,3.00){\vector(3,1){5.00}} \put(0.00,3.00){\vector(3,-1){5.00}}
\put(2.50,0.00){\makebox(0,0)[cc]{\small 1}}
\put(2.50,6.00){\makebox(0,0)[cc]{\small 2}}
\end{picture}
(respectively, \special{em:linewidth 0.4pt} \unitlength 1.00mm
\linethickness{0.4pt}
\begin{picture}(8.00,5.00)
\put(4.00,5.00){\vector(1,-2){2.50}} \put(4.00,5.00){\vector(-1,-2){2.50}}
\put(7.67,2.00){\makebox(0,0)[cc]{\small 2}}
\put(0.33,2.00){\makebox(0,0)[cc]{\small 1}}
\end{picture}
 )
says that the vertical direction lies in the angle between the vectors
$-f'(a_1)$ and $f'(a_2)$ (respectively, between the vectors $-f'(a_1)$ and
$-f'(a_2)$). The subscript \unitlength 1.00mm \special{em:linewidth 0.4pt}
\linethickness{0.4pt}
\begin{picture}(12.00,1.00)
\put(6.00,-1.00){\vector(1,0){6.00}} \put(6.00,-1.00){\vector(-1,0){6.00}}
\put(6.00,-1.00){\circle*{1.00}} \put(1.00,1.00){\makebox(0,0)[cc]{\small 1}}
\put(11.00,1.00){\makebox(0,0)[cc]{\small 2}}
\end{picture}
means that the tangent directions $f'(a_1)$ and $f'(a_2)$ are opposite in
$\R^n$.  The notation of all these types appears only if the condition
$f(a_1)=f(a_2)$ is satisfied (and can be seen from the chords and zigzags on
the picture).

The subscript \unitlength 1.00mm \special{em:linewidth 0.4pt}
\linethickness{0.4pt}
\begin{picture}(6.00,5.00)
\put(0.00,1.00){\vector(2,-1){4.00}} \put(0.00,1.00){\vector(2,1){4.00}}
\put(0.00,1.00){\vector(1,0){6.00}} \put(2.00,-2.00){\makebox(0,0)[cc]{$2$}}
\put(2.00,4.00){\makebox(0,0)[cc]{$1$}} \put(0.00,-0.50){\line(0,1){3.00}}
\end{picture}
means that the distinguished direction "to the right" in ${\bf R}^{n-1}$ lies
between the projections of such tangents $f'(a_1),$ $f'(a_2)$ to $\R^{n-1}$
(i.e. this direction is the linear combination of these projections with
nonnegative coefficients). The subscript \unitlength 1.00mm
\special{em:linewidth 0.4pt} \linethickness{0.4pt}
\begin{picture}(11.50,3.00)
\put(6.00,1.00){\makebox(0,0)[cc]{$1 \longleftrightarrow 2$}}
\put(6.00,-0.50){\line(0,1){3.00}}
\end{picture}
\/ means that the projections of tangents $f'(a_1), f'(a_2)$ to $\R^{n-1}$ have
opposite directions. The notation of last two types can appear only if the
projections of corresponding points $f(a_1), f(a_2)$ coincide in $\R^{n-1}$.

The sum of varieties distinguished by conditions of types \unitlength 1.00mm
\special{em:linewidth 0.4pt} \linethickness{0.4pt}
\begin{picture}(8.00,5.00)
\put(4.00,0.00){\vector(1,2){2.50}} \put(4.00,0.00){\vector(-1,2){2.50}}
\put(0.33,2.00){\makebox(0,0)[cc]{\small 1}}
\put(7.67,2.00){\makebox(0,0)[cc]{\small 2}}
\end{picture}
 \
and \special{em:linewidth 0.4pt} \unitlength 1.00mm \linethickness{0.4pt}
\begin{picture}(8.00,5.00)
\put(4.00,5.00){\vector(1,-2){2.50}} \put(4.00,5.00){\vector(-1,-2){2.50}}
\put(7.67,2.00){\makebox(0,0)[cc]{\small 2}}
\put(0.33,2.00){\makebox(0,0)[cc]{\small 1}}
\end{picture}
 \
in one and the same cell is equal to the variety of type \unitlength 1.00mm
\special{em:linewidth 0.4pt} \linethickness{0.4pt}
\begin{picture}(11.50,3.00)
\put(6.00,1.00){\makebox(0,0)[cc]{$1 \longleftrightarrow 2$}}
\put(6.00,-0.50){\line(0,1){3.00}}
\end{picture}
 \ ;
"of type" here means that some other two numbers instead of 1 and 2 can stay in
all three pictures. This identity is not symmetric: indeed, the variety of type
\unitlength 1.00mm \special{em:linewidth 0.4pt} \linethickness{0.4pt}
\begin{picture}(11.50,3.00)
\put(6.00,1.00){\makebox(0,0)[cc]{$1 \longleftrightarrow 2$}}
\put(6.00,-0.50){\line(0,1){3.00}}
\end{picture}
 \
can be well defined even when the former two varieties have no sense.

The subscript $2 \updownarrow$ (respectively, $2 \uparrow$, respectively, $2
\downarrow$) means that the tangent vector $f'(a_2)$ is vertical (respectively,
vertical directed up, respectively, vertical directed down). The subscript of
type $2 \mapsto$ means that the projection of the tangent $f'(a_2)$ to
$\R^{n-1}$ is directed "to the right".

Abbreviation $f_1$ replaces the composition ${\bf p} \circ f: \R^1 \to
\R^{n-1}$. Finally, a collection of vectors in a framebox means that these
vectors are linearly dependent. For instance the subscript \unitlength 1.00mm
\special{em:linewidth 0.4pt} \linethickness{0.4pt}
\begin{picture}(28.00,5.00)
\put(0.00,-1.00){\framebox(28.00,5.00)[cc]{$f_1'(1),f_1''(2),\mapsto$}}
\end{picture}
means that some three vectors in $\R^{n-1}$, namely the projection of
$f'(a_1)$, the projection of $f''(a_2)$, and the direction "to the right", span
a subspace of dimension $\le 2$. Several more specific abbreviations will be
explained later, close to their first use.

The boundary of the variety distinguished in any cell by the condition
\unitlength 1.00mm \special{em:linewidth 0.4pt} \linethickness{0.4pt}
\begin{picture}(8.00,5.00)
\put(4.00,0.00){\vector(1,2){2.50}} \put(4.00,0.00){\vector(-1,2){2.50}}
\put(0.33,2.00){\makebox(0,0)[cc]{\small 1}}
\put(7.67,2.00){\makebox(0,0)[cc]{\small 2}}
\end{picture}
 \ (respectively,
\special{em:linewidth 0.4pt} \unitlength 1.00mm \linethickness{0.4pt}
\begin{picture}(8.00,5.00)
\put(4.00,5.00){\vector(1,-2){2.50}} \put(4.00,5.00){\vector(-1,-2){2.50}}
\put(7.67,2.00){\makebox(0,0)[cc]{\small 2}}
\put(0.33,2.00){\makebox(0,0)[cc]{\small 1}}
\end{picture}
 \ )
is equal to the sum of varieties distinguished by conditions \unitlength 1.00mm
\special{em:linewidth 0.4pt} \linethickness{0.4pt}
\begin{picture}(12.00,1.00)
\put(6.00,-1.00){\vector(1,0){6.00}} \put(6.00,-1.00){\vector(-1,0){6.00}}
\put(6.00,-1.00){\circle*{1.00}} \put(1.00,1.00){\makebox(0,0)[cc]{\small 1}}
\put(11.00,1.00){\makebox(0,0)[cc]{\small 2}}
\end{picture}
 \ ,
$1 \uparrow$ and $2 \uparrow$ (respectively, \unitlength 1.00mm
\special{em:linewidth 0.4pt} \linethickness{0.4pt}
\begin{picture}(12.00,1.00)
\put(6.00,-1.00){\vector(1,0){6.00}} \put(6.00,-1.00){\vector(-1,0){6.00}}
\put(6.00,-1.00){\circle*{1.00}} \put(1.00,1.00){\makebox(0,0)[cc]{\small 1}}
\put(11.00,1.00){\makebox(0,0)[cc]{\small 2}}
\end{picture}
 \ ,
$1 \downarrow$ and $2 \downarrow$) plus maybe something in the boundary of the
cell. Similarly, the boundary of the variety distinguished by the condition
\unitlength 1.00mm \special{em:linewidth 0.4pt} \linethickness{0.4pt}
\begin{picture}(6.00,5.00)
\put(0.00,1.00){\vector(2,-1){4.00}} \put(0.00,1.00){\vector(2,1){4.00}}
\put(0.00,1.00){\vector(1,0){6.00}} \put(2.00,-2.00){\makebox(0,0)[cc]{$2$}}
\put(2.00,4.00){\makebox(0,0)[cc]{$1$}} \put(0.00,-0.50){\line(0,1){3.00}}
\end{picture}
is equal (modulo the boundary of the cell) to the sum of varieties
distinguished in the same cell by conditions \unitlength 1.00mm
\special{em:linewidth 0.4pt} \linethickness{0.4pt}
\begin{picture}(11.50,3.00)
\put(6.00,1.00){\makebox(0,0)[cc]{$1 \longleftrightarrow 2$}}
\put(6.00,-0.50){\line(0,1){3.00}}
\end{picture}
 \ ,
$1 \mapsto$, and $2 \mapsto$. The boundary of the condition $2 \mapsto$ is
equal to $2 \updownarrow$ plus something in smaller cells.

\subsection{The first differential}
\label{f2}

Formula (\ref{prinpart}) defines a relative cycle in the term $F_3$ of our
filtration modulo $F_2.$ In this subsection we calculate its boundary in the
term $F_2 \setminus F_1,$ and span it by some chain with closed supports in
this term (i.e. we represent it as the boundary of such a chain).

\begin{proposition}
The boundary of the cycle $($\ref{prinpart}$)$ in $F_2\setminus F_1$ is equal
to the chain
\begin{equation}
\label{difone} \unitlength=1.00mm \special{em:linewidth 0.4pt}
\linethickness{0.4pt}
\begin{picture}(120.00,16.00)
\put(0.00,12.00){\line(1,0){20.00}} \put(25.00,12.00){\line(1,0){20.00}}
\put(50.00,12.00){\line(1,0){20.00}} \put(75.00,12.00){\line(1,0){20.00}}
\put(100.00,12.00){\line(1,0){20.00}} \put(110.00,12.00){\oval(16.00,8.00)[t]}
\put(106.00,12.00){\oval(8.00,4.00)[t]} \put(114.00,12.00){\oval(8.00,4.00)[t]}
\put(118.00,12.00){\line(-2,-1){6.00}} \put(112.00,9.00){\line(-2,1){6.00}}
\put(85.00,12.00){\oval(16.00,8.00)[t]} \put(81.00,12.00){\oval(8.00,4.00)[t]}
\put(89.00,12.00){\oval(8.00,4.00)[t]} \put(84.00,8.00){\vector(1,0){6.00}}
\put(84.00,8.00){\vector(-1,0){6.00}} \put(80.00,10.00){\makebox(0,0)[cc]{$1$}}
\put(88.00,10.00){\makebox(0,0)[cc]{$2$}} \put(84.00,8.00){\circle*{1.00}}
\put(58.50,12.00){\oval(13.00,8.00)[b]} \put(64.00,12.00){\oval(8.00,8.00)[t]}
\put(33.00,12.00){\oval(12.00,8.00)[b]} \put(36.50,12.00){\oval(13.00,8.00)[t]}
\put(43.00,12.00){\line(-2,1){6.00}} \put(37.00,15.00){\line(-4,-3){4.00}}
\put(11.50,12.00){\oval(13.00,8.00)[t]} \put(13.50,12.00){\oval(9.00,4.00)[t]}
\put(9.00,8.33){\line(-2,1){7.00}} \put(10.00,2.00){\makebox(0,0)[cc]{$A$}}
\put(35.00,2.00){\makebox(0,0)[cc]{$B$}}
\put(60.00,2.00){\makebox(0,0)[cc]{$C$}}
\put(85.00,2.00){\makebox(0,0)[cc]{$D$}}
\put(110.00,2.00){\makebox(0,0)[cc]{$E$}}
\put(22.50,12.00){\makebox(0,0)[cc]{$+$}}
\put(47.50,12.00){\makebox(0,0)[cc]{$+$}}
\put(72.50,12.00){\makebox(0,0)[cc]{$+$}}
\put(97.50,12.00){\makebox(0,0)[cc]{$+$}} \put(9.00,8.00){\line(3,2){6.00}}
\put(60.00,12.00){\line(-1,2){2.00}} \put(58.00,16.00){\line(-1,-2){2.00}}
\end{picture}
 \ .
\end{equation}
$($Namely, the boundary of the first term of $($\ref{prinpart}$)$ consists of
four chains A, B, C and D in $($\ref{difone}$)$, and the boundary of the second
is equal to the fifth chain E.$)$
\end{proposition}

The unique nontrivial term of this formula is the 4th one: it appears when the
first point of the 5-configuration participating in the first term of
(\ref{prinpart}) tends to the second, and simultaneously the fourth point tends
to the third. \quad $\square$
\medskip

{\bf Exercise:} to check that the chain (\ref{difone}) actually is a cycle in
$F_2 \setminus F_1$.
\medskip

Now, let us span this cycle by a chain in the term $F_2 \sm F_1$. The cellular
structure of this term was described in Example \ref{exf2} of section
\ref{method}.

First we span the components D and E of (\ref{difone}) inside the cell
\unitlength 1.00mm \special{em:linewidth 0.4pt} \linethickness{0.4pt}
\begin{picture}(18.00,4.00)
\put(0.00,0.00){\line(1,0){16.00}} \put(8.00,0.00){\oval(14.00,8.00)[t]}
\put(4.50,0.00){\oval(7.00,4.00)[t]} \put(11.50,0.00){\oval(7.00,4.00)[t]}
\end{picture}
, i.e. we construct the homology between their sum and some chain in the
boundary of this cell. It is natural to span a chain with condition of type
\unitlength 1.00mm \special{em:linewidth 0.4pt} \linethickness{0.4pt}
\begin{picture}(12.00,1.00)
\put(6.00,-1.00){\vector(1,0){6.00}} \put(6.00,-1.00){\vector(-1,0){6.00}}
\put(6.00,-1.00){\circle*{1.00}} \put(1.00,1.00){\makebox(0,0)[cc]{\small 1}}
\put(11.00,1.00){\makebox(0,0)[cc]{\small 2}}
\end{picture}
 \ by a similar chain with condition of type
\unitlength 1.00mm \special{em:linewidth 0.4pt} \linethickness{0.4pt}
\begin{picture}(8.00,5.00)
\put(4.00,0.00){\vector(1,2){2.50}} \put(4.00,0.00){\vector(-1,2){2.50}}
\put(0.33,2.00){\makebox(0,0)[cc]{\small 1}}
\put(7.67,2.00){\makebox(0,0)[cc]{\small 2}}
\end{picture}
 ,
and a chain having zigzag without arrows by a similar chain with an arrow added
at one of endpoints of the zigzag. The chains obtained in this way from the
ones encoded by parts D and E of (\ref{difone}) are indicated in the left parts
of the next two equations (\ref{boundd}) and (\ref{bounde}) respectively.

In the right-hand parts of these formulas, as well as in all forthcoming
expressions for boundary operators in this work, we first count the components
of the boundary defined by the degenerations of the subvarieties in the
corresponding cells, distinguished by arrowed zigzags and subscripts. Then we
count the components defined by the limit positions of these varieties when the
cell itself degenerates because of the collision of some points forming its
underlying $J$-configuration in $\R^1$. The latter degenerations appear in the
lexicographic order: first by the number of colliding pairs of points in
$\R^1$, and then by their positions in $\R^1$.

\begin{equation}
\label{boundd} \unitlength=1.00mm \special{em:linewidth 0.4pt}
\linethickness{0.4pt}
\begin{picture}(114.00,28.00)
\put(3.00,24.00){\makebox(0,0)[cc]{$\partial$}}
\put(6.00,24.00){\line(1,0){20.00}} \put(16.00,24.00){\oval(16.00,8.00)[t]}
\put(12.00,24.00){\oval(8.00,4.00)[t]} \put(20.00,24.00){\oval(8.00,4.00)[t]}
\put(14.00,17.00){\vector(1,1){5.00}} \put(14.00,17.00){\vector(-1,1){5.00}}
\put(20.00,19.00){\makebox(0,0)[cc]{$1$}}
\put(8.00,19.00){\makebox(0,0)[cc]{$2$}}
\put(30.00,24.00){\makebox(0,0)[cc]{=}} \put(36.00,24.00){\line(1,0){20.00}}
\put(46.00,24.00){\oval(16.00,8.00)[t]} \put(42.00,24.00){\oval(8.00,4.00)[t]}
\put(50.00,24.00){\oval(8.00,4.00)[t]} \put(63.00,24.00){\line(1,0){20.00}}
\put(73.00,24.00){\oval(16.00,8.00)[t]} \put(69.00,24.00){\oval(8.00,4.00)[t]}
\put(77.00,24.00){\oval(8.00,4.00)[t]} \put(90.00,24.00){\line(1,0){20.00}}
\put(100.00,24.00){\oval(16.00,8.00)[t]} \put(96.00,24.00){\oval(8.00,4.00)[t]}
\put(104.00,24.00){\oval(8.00,4.00)[t]} \put(87.00,24.00){\makebox(0,0)[cc]{+}}
\put(60.00,24.00){\makebox(0,0)[cc]{+}} \put(41.00,20.00){\makebox(0,0)[cc]{$1
\uparrow$}} \put(70.00,20.00){\makebox(0,0)[cc]{$2 \uparrow$}}
\put(98.00,20.00){\vector(1,0){6.00}} \put(98.00,20.00){\vector(-1,0){6.00}}
\put(98.00,20.00){\circle*{1.00}} \put(102.00,22.00){\makebox(0,0)[cc]{$1$}}
\put(94.00,22.00){\makebox(0,0)[cc]{$2$}}
\put(114.00,24.00){\makebox(0,0)[cc]{$+$}} \put(36.00,7.00){\line(1,0){20.00}}
\put(46.00,7.00){\oval(16.00,8.00)[t]} \put(50.00,7.00){\oval(8.00,4.00)[t]}
\put(44.00,0.00){\vector(1,1){5.00}} \put(44.00,0.00){\vector(-1,1){5.00}}
\put(50.00,2.00){\makebox(0,0)[cc]{$1$}}
\put(38.00,2.00){\makebox(0,0)[cc]{$2$}} \put(63.00,7.00){\line(1,0){20.00}}
\put(69.00,7.00){\oval(8.00,4.00)[t]} \put(77.00,7.00){\oval(8.00,4.00)[t]}
\put(71.00,0.00){\vector(1,1){5.00}} \put(71.00,0.00){\vector(-1,1){5.00}}
\put(77.00,2.00){\makebox(0,0)[cc]{$1$}}
\put(65.00,2.00){\makebox(0,0)[cc]{$2$}} \put(90.00,7.00){\line(1,0){20.00}}
\put(100.00,7.00){\oval(16.00,8.00)[t]} \put(96.00,7.00){\oval(8.00,4.00)[t]}
\put(98.00,0.00){\vector(1,1){5.00}} \put(98.00,0.00){\vector(-1,1){5.00}}
\put(104.00,2.00){\makebox(0,0)[cc]{$1$}}
\put(92.00,2.00){\makebox(0,0)[cc]{$2$}}
\put(87.00,7.00){\makebox(0,0)[cc]{$+$}}
\put(60.00,7.00){\makebox(0,0)[cc]{$+$}}
\put(32.00,7.00){\makebox(0,0)[cc]{$+$}}
\end{picture}
 ,
\end{equation}
\begin{equation}
\label{bounde} \unitlength=1.00mm \special{em:linewidth 0.4pt}
\linethickness{0.4pt}
\begin{picture}(108.00,23.00)
\put(3.00,19.00){\makebox(0,0)[cc]{$\partial$}}
\put(6.00,19.00){\line(1,0){20.00}} \put(16.00,19.00){\oval(16.00,8.00)[t]}
\put(12.00,19.00){\oval(8.00,4.00)[t]} \put(20.00,19.00){\oval(8.00,4.00)[t]}
\put(24.00,19.00){\line(-2,-1){6.00}} \put(18.00,16.00){\vector(-2,1){6.00}}
\put(32.00,19.00){\line(1,0){20.00}} \put(42.00,19.00){\oval(16.00,8.00)[t]}
\put(38.00,19.00){\oval(8.00,4.00)[t]} \put(46.00,19.00){\oval(8.00,4.00)[t]}
\put(58.00,19.00){\line(1,0){20.00}} \put(68.00,19.00){\oval(16.00,8.00)[t]}
\put(64.00,19.00){\oval(8.00,4.00)[t]} \put(72.00,19.00){\oval(8.00,4.00)[t]}
\put(84.00,19.00){\line(1,0){20.00}} \put(94.00,19.00){\oval(16.00,8.00)[t]}
\put(90.00,19.00){\oval(8.00,4.00)[t]} \put(98.00,19.00){\oval(8.00,4.00)[t]}
\put(81.00,19.00){\makebox(0,0)[cc]{$+$}}
\put(55.00,19.00){\makebox(0,0)[cc]{$+$}}
\put(29.00,19.00){\makebox(0,0)[cc]{$=$}}
\put(36.00,14.00){\makebox(0,0)[cc]{$1 \downarrow$}}
\put(66.00,14.00){\makebox(0,0)[cc]{$2 \uparrow$}}
\put(102.00,19.00){\line(-2,-1){6.00}} \put(96.00,16.00){\line(-2,1){6.00}}
\put(108.00,19.00){\makebox(0,0)[cc]{$+$}} \put(32.00,3.00){\line(1,0){20.00}}
\put(42.00,3.00){\oval(16.00,8.00)[t]} \put(46.00,3.00){\oval(8.00,4.00)[t]}
\put(50.00,3.00){\line(-2,-1){6.00}} \put(44.00,0.00){\vector(-2,1){6.00}}
\put(58.00,3.00){\line(1,0){20.00}} \put(64.00,3.00){\oval(8.00,4.00)[t]}
\put(72.00,3.00){\oval(8.00,4.00)[t]} \put(76.00,3.00){\line(-2,-1){6.00}}
\put(70.00,0.00){\vector(-2,1){6.00}} \put(84.00,3.00){\line(1,0){20.00}}
\put(94.00,3.00){\oval(16.00,8.00)[t]} \put(90.00,3.00){\oval(8.00,4.00)[t]}
\put(102.00,3.00){\line(-2,-1){6.00}} \put(96.00,0.00){\vector(-2,1){6.00}}
\put(81.00,3.00){\makebox(0,0)[cc]{$+$}}
\put(55.00,3.00){\makebox(0,0)[cc]{$+$}}
\put(29.00,3.00){\makebox(0,0)[cc]{$+$}}
\end{picture}
 .
\end{equation}

\begin{proposition}
The equalities $($\ref{boundd}$)$, $($\ref{bounde}$)$ are correct, i.e. the
algebraic boundaries $($mod 2$)$ in $F_2 \setminus F_1$ of the varieties
indicated in their left parts are equal to the sums of varieties indicated in
their right-hand parts. \quad $\square$
\end{proposition}

In (\ref{bounde}) first two summands are degenerations of the variety defined
by the zigzag when its arrowed endpoint tends to one of boundaries of the
corresponding segment, and the third summand belongs to its boundary as the
equality of type $\phi(x)=\phi(y)$ defines a component of the boundary of the
set defined by the inequality $\phi(x) \ge \phi(y)$.

The last three summands in both (\ref{boundd}) and (\ref{bounde}) belong to the
boundary (\ref{ztripno}) of the cell \unitlength 1.00mm \special{em:linewidth
0.4pt} \linethickness{0.4pt}
\begin{picture}(18.00,4.00)
\put(0.00,0.00){\line(1,0){16.00}} \put(8.00,0.00){\oval(14.00,8.00)[t]}
\put(4.50,0.00){\oval(7.00,4.00)[t]} \put(11.50,0.00){\oval(7.00,4.00)[t]}
\end{picture}
 .

The sum of all varieties indicated in right-hand parts of (\ref{boundd}),
(\ref{bounde}) consists of part D + E of (\ref{difone}), some chain in the
boundary of the cell \unitlength 1.00mm \special{em:linewidth 0.4pt}
\linethickness{0.4pt}
\begin{picture}(18.00,4.00)
\put(0.00,0.00){\line(1,0){16.00}} \put(8.00,0.00){\oval(14.00,8.00)[t]}
\put(4.50,0.00){\oval(7.00,4.00)[t]} \put(11.50,0.00){\oval(7.00,4.00)[t]}
\end{picture}
 ,
and the first chain in the right-hand part of the equation
\begin{equation}
\label{nev} \unitlength=1.00mm \special{em:linewidth 0.4pt}
\linethickness{0.4pt}
\begin{picture}(130.00,8.00)
\put(6.00,4.00){\line(1,0){20.00}} \put(16.00,4.00){\oval(16.00,8.00)[t]}
\put(12.00,4.00){\oval(8.00,4.00)[t]} \put(20.00,4.00){\oval(8.00,4.00)[t]}
\put(32.00,4.00){\line(1,0){20.00}} \put(42.00,4.00){\oval(16.00,8.00)[t]}
\put(38.00,4.00){\oval(8.00,4.00)[t]} \put(46.00,4.00){\oval(8.00,4.00)[t]}
\put(37.00,1.00){\makebox(0,0)[cc]{$1 \updownarrow$}}
\put(29.00,4.00){\makebox(0,0)[cc]{$=$}}
\put(3.00,4.00){\makebox(0,0)[cc]{$\partial$}}
\put(58.00,4.00){\line(1,0){20.00}} \put(68.00,4.00){\oval(16.00,8.00)[t]}
\put(72.00,4.00){\oval(8.00,4.00)[t]} \put(84.00,4.00){\line(1,0){20.00}}
\put(90.00,4.00){\oval(8.00,4.00)[t]} \put(98.00,4.00){\oval(8.00,4.00)[t]}
\put(110.00,4.00){\line(1,0){20.00}} \put(120.00,4.00){\oval(16.00,8.00)[t]}
\put(116.00,4.00){\oval(8.00,4.00)[t]} \put(81.00,4.00){\makebox(0,0)[cc]{$+$}}
\put(107.00,4.00){\makebox(0,0)[cc]{$+$}}
\put(55.00,4.00){\makebox(0,0)[cc]{$+$}} \put(116.00,1.00){\makebox(0,0)[cc]{$1
\mapsto$}} \put(90.00,1.00){\makebox(0,0)[cc]{$1 \mapsto$}}
\put(64.00,1.00){\makebox(0,0)[cc]{$1 \mapsto$}}
\put(12.00,1.00){\makebox(0,0)[cc]{$1 \mapsto$}}
\end{picture}
 \ .
\end{equation}

In other words, the sum of three chains in the left parts of equations
(\ref{boundd}), (\ref{bounde}) and (\ref{nev}) realizes homology between the
sum of chains $D$ and $E$ and some chain in the boundary of the cell
\unitlength 1.00mm \special{em:linewidth 0.4pt} \linethickness{0.4pt}
\begin{picture}(18.00,4.00)
\put(0.00,0.00){\line(1,0){16.00}} \put(8.00,0.00){\oval(14.00,8.00)[t]}
\put(4.50,0.00){\oval(7.00,4.00)[t]} \put(11.50,0.00){\oval(7.00,4.00)[t]}
\end{picture}
 \/.

Now we span the summands B and C of (\ref{difone}) inside the open cell
\unitlength 1.00mm \special{em:linewidth 0.4pt} \linethickness{0.4pt}
\begin{picture}(16.00,5.00)
\put(0.00,2.00){\line(1,0){16.00}} \put(6.50,2.00){\oval(9.00,4.00)[b]}
\put(9.50,2.00){\oval(9.00,4.00)[t]}
\end{picture}
 . We need to find varieties in this cell,
whose boundaries include these summands. The obvious candidates for this are
the chains shown in the left parts of equations (\ref{boundb}) and
(\ref{boundc}) respectively.

\begin{equation}
\unitlength 1.00mm \linethickness{0.4pt}
\begin{picture}(121.00,29.00)
\put(3.00,25.00){\makebox(0,0)[cc]{$\partial$}}
\put(7.00,25.00){\line(1,0){23.00}} \put(17.00,25.00){\oval(16.00,8.00)[b]}
\put(20.00,25.00){\oval(16.00,8.00)[t]}
\put(33.00,25.00){\makebox(0,0)[cc]{$=$}} \put(36.00,25.00){\line(1,0){23.00}}
\put(46.00,25.00){\oval(16.00,8.00)[b]} \put(49.00,25.00){\oval(16.00,8.00)[t]}
\put(36.00,8.00){\line(1,0){23.00}} \put(46.00,8.00){\oval(16.00,8.00)[b]}
\put(49.00,8.00){\oval(16.00,8.00)[t]} \put(94.00,25.00){\line(1,0){23.00}}
\put(104.00,25.00){\oval(16.00,8.00)[b]}
\put(107.00,25.00){\oval(16.00,8.00)[t]}
\put(120.00,25.00){\makebox(0,0)[cc]{$+$}}
\put(91.00,25.00){\makebox(0,0)[cc]{$+$}}
\put(33.00,8.00){\makebox(0,0)[cc]{$+$}} \put(41.00,8.00){\line(3,-1){9.00}}
\put(50.00,5.00){\vector(4,3){4.00}} \put(114.00,19.00){\makebox(0,0)[cc]{$2
\downarrow$}} \put(62.00,25.00){\makebox(0,0)[cc]{$+$}}
\put(65.00,25.00){\line(1,0){23.00}} \put(76.50,25.00){\oval(19.00,8.00)[t]}
\put(74.00,25.00){\oval(14.00,8.00)[b]}
\put(92.00,8.00){\makebox(0,0)[cc]{$+$}} \put(95.00,8.00){\line(1,0){23.00}}
\put(106.50,8.00){\oval(19.00,8.00)[t]} \put(104.00,8.00){\oval(14.00,8.00)[b]}
\put(62.00,8.00){\makebox(0,0)[cc]{$+$}} \put(66.00,8.00){\line(1,0){23.00}}
\put(77.50,8.00){\oval(19.00,8.00)[b]} \put(80.00,8.00){\oval(14.00,8.00)[t]}
\put(116.00,0.00){\vector(-3,4){3.67}} \put(116.00,0.00){\vector(3,4){3.67}}
\put(121.00,3.00){\makebox(0,0)[cc]{$1$}}
\put(111.00,3.00){\makebox(0,0)[cc]{$2$}} \put(12.00,25.00){\line(4,-3){4.00}}
\put(16.00,22.00){\vector(4,3){4.00}} \put(41.00,25.00){\line(4,-3){4.00}}
\put(45.00,22.00){\line(4,3){4.00}} \put(81.00,25.00){\line(-4,3){4.00}}
\put(77.00,28.00){\vector(-4,-3){4.00}} \put(73.00,8.00){\line(4,-3){4.00}}
\put(77.00,5.00){\vector(4,3){4.00}}
\end{picture}
\label{boundb}
\end{equation}

\begin{equation}
\label{boundc} \unitlength 1.00mm \linethickness{0.4pt}
\begin{picture}(113.00,28.00)
\put(3.00,24.00){\makebox(0,0)[cc]{$\partial$}}
\put(6.00,24.00){\line(1,0){21.00}} \put(14.00,24.00){\oval(14.00,8.00)[b]}
\put(20.50,24.00){\oval(11.00,8.00)[t]}
\put(30.00,24.00){\makebox(0,0)[cc]{$=$}} \put(33.00,24.00){\line(1,0){22.00}}
\put(41.50,24.00){\oval(15.00,8.00)[b]} \put(48.50,24.00){\oval(11.00,8.00)[t]}
\put(58.00,24.00){\makebox(0,0)[cc]{$+$}} \put(61.00,24.00){\line(1,0){21.00}}
\put(69.50,24.00){\oval(13.00,8.00)[b]} \put(73.00,24.00){\oval(14.00,8.00)[t]}
\put(66.00,24.00){\line(2,-1){6.00}} \put(72.00,21.00){\vector(4,3){4.00}}
\put(85.00,24.00){\makebox(0,0)[cc]{$+$}} \put(88.00,24.00){\line(1,0){21.00}}
\put(97.00,24.00){\oval(14.00,8.00)[b]}
\put(100.50,24.00){\oval(13.00,8.00)[t]}
\put(107.00,19.00){\makebox(0,0)[cc]{$2 \uparrow$}}
\put(113.00,24.00){\makebox(0,0)[cc]{$+$}} \put(7.00,9.00){\line(1,0){21.00}}
\put(13.50,9.00){\oval(9.00,8.00)[b]} \put(22.00,9.00){\oval(8.00,8.00)[t]}
\put(18.00,9.00){\line(-1,2){2.00}} \put(16.00,13.00){\vector(-1,-2){2.00}}
\put(58.00,9.00){\makebox(0,0)[cc]{$+$}} \put(61.00,9.00){\line(1,0){21.00}}
\put(67.50,9.00){\oval(9.00,8.00)[b]} \put(76.00,9.00){\oval(8.00,8.00)[t]}
\put(77.00,7.00){\vector(-2,-3){4.00}} \put(77.00,7.00){\vector(2,-3){4.00}}
\put(82.00,4.00){\makebox(0,0)[cc]{$1$}}
\put(72.00,4.00){\makebox(0,0)[cc]{$2$}}
\put(30.00,9.00){\makebox(0,0)[cc]{$+$}} \put(33.00,9.00){\line(1,0){21.00}}
\put(43.50,9.00){\oval(17.00,8.00)[b]} \put(47.50,9.00){\oval(9.00,8.00)[t]}
\put(43.00,9.00){\line(-1,2){2.00}} \put(41.00,13.00){\vector(-1,-2){2.00}}
\put(85.00,9.00){\makebox(0,0)[cc]{$+$}} \put(88.00,9.00){\line(1,0){21.00}}
\put(98.50,9.00){\oval(17.00,8.00)[b]} \put(102.50,9.00){\oval(9.00,8.00)[t]}
\put(15.00,24.00){\line(-3,4){3.00}} \put(12.00,28.00){\vector(-3,-4){3.00}}
\put(43.00,24.00){\line(-3,4){3.00}} \put(40.00,28.00){\line(-3,-4){3.00}}
\put(4.00,9.00){\makebox(0,0)[cc]{$+$}} \put(106.00,4.00){\vector(3,-1){6.00}}
\put(106.00,4.00){\vector(3,1){6.00}}
\put(109.00,7.00){\makebox(0,0)[cc]{{\small 3}}}
\put(109.00,1.00){\makebox(0,0)[cc]{{\small 1}}}
\end{picture}
\end{equation}

Again, all summands in lower rows of these equalities belong to the boundary of
the cell \unitlength 1.00mm \special{em:linewidth 0.4pt} \linethickness{0.4pt}
\begin{picture}(16.00,5.00)
\put(0.00,2.00){\line(1,0){16.00}} \put(6.50,2.00){\oval(9.00,4.00)[b]}
\put(9.50,2.00){\oval(9.00,4.00)[t]}
\end{picture}
 .

\begin{proposition}
\label{propbc} The equalities $($\ref{boundb}$)$, $($\ref{boundc}$)$ are
correct, i.e. the algebraic $($mod 2$)$ boundaries in $F_2 \setminus F_1$ of
the varieties indicated in their left parts are equal to the sums of varieties
indicated in their right-hand parts. \quad $\square$
\end{proposition}

In particular we get that the boundary of the sum of these two left-side
varieties is equal to the sum of varieties denoted in (\ref{difone}) by B and
C, plus some chain in the boundary of the cell \unitlength 1.00mm
\special{em:linewidth 0.4pt} \linethickness{0.4pt}
\begin{picture}(16.00,5.00)
\put(0.00,2.00){\line(1,0){16.00}} \put(6.50,2.00){\oval(9.00,4.00)[b]}
\put(9.50,2.00){\oval(9.00,4.00)[t]}
\end{picture}
 , plus the variety distinguished in this cell by the additional
condition $2 \updownarrow$. The last variety is a part of the boundary of the
similar set distinguished by the condition $2 \mapsto$. Entire boundary of this
set in $F_2 \setminus F_1$ is expressed by the formula
\begin{equation}
\label{nev2} \unitlength=1.00mm \special{em:linewidth 0.4pt}
\linethickness{0.4pt}
\begin{picture}(127.00,11.00)
\put(16.00,0.00){\makebox(0,0)[cc]{$2 \mapsto$}}
\put(38.00,0.00){\makebox(0,0)[cc]{$2 \updownarrow$}}
\put(56.00,0.00){\makebox(0,0)[cc]{$1 \mapsto$}}
\put(74.00,0.00){\makebox(0,0)[cc]{$2 \mapsto$}}
\put(93.00,0.00){\makebox(0,0)[cc]{$2 \mapsto$}}
\put(1.00,7.00){\makebox(0,0)[cc]{$\partial$}}
\put(4.00,7.00){\line(1,0){14.00}} \put(9.50,7.00){\oval(9.00,8.00)[b]}
\put(12.50,7.00){\oval(9.00,8.00)[t]} \put(22.00,7.00){\makebox(0,0)[cc]{$=$}}
\put(25.00,7.00){\line(1,0){14.00}} \put(30.50,7.00){\oval(9.00,8.00)[b]}
\put(33.50,7.00){\oval(9.00,8.00)[t]} \put(42.00,7.00){\makebox(0,0)[cc]{$+$}}
\put(45.00,7.00){\line(1,0){15.00}} \put(50.50,7.00){\oval(9.00,8.00)[b]}
\put(52.50,7.00){\oval(13.00,8.00)[t]} \put(63.00,7.00){\makebox(0,0)[cc]{$+$}}
\put(66.00,7.00){\line(1,0){14.00}} \put(70.00,7.00){\oval(6.00,8.00)[b]}
\put(76.00,7.00){\oval(6.00,8.00)[t]} \put(83.00,7.00){\makebox(0,0)[cc]{$+$}}
\put(86.00,7.00){\line(1,0){14.00}} \put(93.00,7.00){\oval(12.00,8.00)[b]}
\put(95.00,7.00){\oval(8.00,8.00)[t]} \put(103.00,7.00){\makebox(0,0)[cc]{$+$}}
\put(107.00,7.00){\line(1,0){20.00}} \put(115.00,7.00){\oval(9.00,8.00)[t]}
\put(110.50,7.00){\makebox(0,0)[cc]{\large *}}
\put(118.00,4.00){\makebox(0,0)[cc]{\small $-f_{1}'''(1)/f_{1}''(1) \sim$}}
\put(118.00,0.00){\makebox(0,0)[cc]{\small $\sim \ \mapsto / f_{1}''(1)$}}
\end{picture}
 \ .
\end{equation}

The subscript under the last term of (\ref{nev2}) means, that the projections
of second and third derivatives of $f$ at the point $a_1$ into $\R^{n-1}$ lie
in the same 2-plane as the direction "to the right", and two frames in this
2-plane obtained by adding to the projection of $f''(a_1)$ either the
projection of $f'''(a_1)$ or the direction "to the right" have {\em opposite}
orientations. This term occurs when both endpoints $a_1, a_3$ of the "lower"
arc in the left picture of (\ref{nev2}) tend from different sides to the first
endpoint $a_2$ of the "upper" arc.

Finally we get that the cycle $d^1(TT)$ shown in (\ref{difone}) is homologous
in $F_2 \setminus F_1$ to a chain lying in the union of cells of nonmaximal
dimensions listed in (\ref{ztripno}), (\ref{marg}); this homology is provided
by the sum of six varieties indicated in the left parts of equalities
(\ref{boundd}), (\ref{bounde}), (\ref{nev}), (\ref{boundb}), (\ref{boundc}),
and (\ref{nev2}).

Namely, this cycle homologous to $d^1(TT)$ is as follows. In the cell
\unitlength 1.00mm \special{em:linewidth 0.4pt} \linethickness{0.4pt}
\begin{picture}(20.00,4.00)
\put(0.00,0.00){\line(1,0){20.00}} \put(10.00,0.00){\oval(16.00,8.00)[t]}
\put(6.00,0.00){\oval(8.00,4.00)[t]}
\end{picture}
it is zero, in the cell \unitlength 1.00mm \special{em:linewidth 0.4pt}
\linethickness{0.4pt}
\begin{picture}(20.00,4.00)
\put(0.00,0.00){\line(1,0){20.00}} \put(14.00,0.00){\oval(8.00,4.00)[t]}
\put(6.00,0.00){\oval(8.00,4.00)[t]}
\end{picture}
 \ it is equal to the chain
\begin{equation}
\label{homztrip2} \unitlength=1.00mm \special{em:linewidth 0.4pt}
\linethickness{0.4pt}
\begin{picture}(98.00,13.00)
\put(0.00,9.00){\line(1,0){20.00}} \put(6.00,9.00){\oval(8.00,8.00)[b]}
\put(14.00,9.00){\oval(8.00,8.00)[t]} \put(23.00,9.00){\makebox(0,0)[cc]{$+$}}
\put(26.00,9.00){\line(1,0){20.00}} \put(32.00,9.00){\oval(8.00,8.00)[b]}
\put(40.00,9.00){\oval(8.00,8.00)[t]} \put(49.00,9.00){\makebox(0,0)[cc]{$+$}}
\put(52.00,9.00){\line(1,0){20.00}} \put(58.00,9.00){\oval(8.00,8.00)[b]}
\put(66.00,9.00){\oval(8.00,8.00)[t]} \put(75.00,9.00){\makebox(0,0)[cc]{$+$}}
\put(78.00,9.00){\line(1,0){20.00}} \put(84.00,9.00){\oval(8.00,8.00)[b]}
\put(92.00,9.00){\oval(8.00,8.00)[t]} \put(91.00,1.00){\makebox(0,0)[cc]{$2
\mapsto$}} \put(65.00,6.00){\vector(-2,-3){4.00}}
\put(65.00,6.00){\vector(2,-3){4.00}} \put(70.00,3.00){\makebox(0,0)[cc]{1}}
\put(60.00,3.00){\makebox(0,0)[cc]{2}} \put(38.00,1.00){\makebox(0,0)[cc]{$1
\mapsto$}} \put(16.00,0.00){\vector(-2,3){4.00}}
\put(16.00,0.00){\vector(2,3){4.00}} \put(11.00,3.00){\makebox(0,0)[cc]{1}}
\put(21.00,3.00){\makebox(0,0)[cc]{2}}
\end{picture}
 \ ,
\end{equation}
in the cell \unitlength 1.00mm \special{em:linewidth 0.4pt}
\linethickness{0.4pt}
\begin{picture}(20.00,4.00)
\put(0.00,0.00){\line(1,0){20.00}} \put(10.00,0.00){\oval(16.00,8.00)[t]}
\put(14.00,0.00){\oval(8.00,4.00)[t]}
\end{picture}
 \ it is equal to the chain
\begin{equation}
\label{homztrip3} \unitlength 1.00mm \special{em:linewidth 0.4pt}
\linethickness{0.4pt}
\begin{picture}(101.00,28.00)
\put(0.00,24.00){\line(1,0){20.00}} \put(10.00,24.00){\oval(16.00,8.00)[t]}
\put(14.00,24.00){\oval(8.00,4.00)[t]}
\put(23.00,24.00){\makebox(0,0)[cc]{$+$}} \put(26.00,24.00){\line(1,0){20.00}}
\put(36.00,24.00){\oval(16.00,8.00)[t]} \put(40.00,24.00){\oval(8.00,4.00)[t]}
\put(49.00,24.00){\makebox(0,0)[cc]{$+$}} \put(52.00,24.00){\line(1,0){20.00}}
\put(62.00,24.00){\oval(16.00,8.00)[t]} \put(66.00,24.00){\oval(8.00,4.00)[t]}
\put(75.00,24.00){\makebox(0,0)[cc]{$+$}} \put(78.00,24.00){\line(1,0){20.00}}
\put(88.00,24.00){\oval(16.00,8.00)[t]} \put(92.00,24.00){\oval(8.00,4.00)[t]}
\put(101.00,24.00){\makebox(0,0)[cc]{$+$}} \put(52.00,4.00){\line(1,0){20.00}}
\put(62.00,4.00){\oval(16.00,8.00)[t]} \put(66.00,4.00){\oval(8.00,4.00)[t]}
\put(26.00,4.00){\line(1,0){20.00}} \put(37.00,4.00){\oval(14.00,8.00)[t]}
\put(40.00,4.00){\oval(8.00,4.00)[t]} \put(49.00,4.00){\makebox(0,0)[cc]{$+$}}
\put(10.00,15.00){\vector(-1,2){3.00}} \put(10.00,15.00){\vector(1,2){3.00}}
\put(15.00,18.00){\makebox(0,0)[cc]{2}} \put(5.00,18.00){\makebox(0,0)[cc]{1}}
\put(35.00,18.00){\makebox(0,0)[cc]{$1 \mapsto$}}
\put(88.00,18.00){\makebox(0,0)[cc]{$2 \mapsto$}}
\put(70.00,4.00){\line(-1,-2){2.00}} \put(68.00,0.00){\vector(-1,2){2.00}}
\put(40.00,4.00){\line(-5,-3){7.00}} \put(33.00,-0.33){\line(0,0){0.00}}
\put(33.00,0.00){\line(-3,2){6.00}} \put(23.00,4.00){\makebox(0,0)[cc]{$+$}}
\put(59.00,18.00){\vector(3,1){6.00}} \put(59.00,18.00){\vector(3,-1){6.00}}
\put(61.00,21.00){\makebox(0,0)[cc]{{\small 3}}}
\put(61.00,15.00){\makebox(0,0)[cc]{{\small 1}}}
\end{picture}
 \ ,
\end{equation}
in the 4th cell of (\ref{marg}) it is equal to
\begin{equation}
\label{mar2} \unitlength=1.00mm \special{em:linewidth 0.4pt}
\linethickness{0.4pt}
\begin{picture}(26.00,13.00)
\put(2.00,9.00){\line(1,0){24.00}} \put(13.00,9.00){\oval(13.00,8.00)[t]}
\put(6.50,9.00){\makebox(0,0)[cc]{\large *}}
\put(15.00,6.00){\makebox(0,0)[cc]{\small $-f_{1}'''(1)/f_{1}''(1) \sim$}}
\put(15.00,2.00){\makebox(0,0)[cc]{\small $\sim \ \mapsto / f_{1}''(1)$}}
\end{picture}
 \ \qquad ,
\end{equation}
and its intersections with all other cells are empty.

The sum of the first and the third terms in (\ref{homztrip2}) is equal to the
variety denoted by the subscript \unitlength 1.00mm \special{em:linewidth
0.4pt} \linethickness{0.4pt}
\begin{picture}(11.50,3.00)
\put(6.00,1.00){\makebox(0,0)[cc]{$1 \longleftrightarrow 2$}}
\put(6.00,-0.50){\line(0,1){3.00}}
\end{picture}
 \ .
To kill it (and something else) we consider the equality
\begin{equation}
\label{nev3} \unitlength 1.00mm \special{em:linewidth 0.4pt}
\linethickness{0.4pt}
\begin{picture}(125.00,17.00)
\put(5.00,14.00){\line(1,0){16.00}} \put(9.50,14.00){\oval(7.00,6.00)[t]}
\put(16.50,14.00){\oval(7.00,6.00)[t]} \put(9.00,5.00){\vector(2,1){6.00}}
\put(9.00,5.00){\vector(2,-1){6.00}} \put(9.00,5.00){\vector(1,0){10.00}}
\put(9.00,3.00){\line(0,1){4.00}} \put(12.00,9.00){\makebox(0,0)[cc]{\small
$1$}} \put(12.00,1.00){\makebox(0,0)[cc]{\small $2$}}
\put(2.00,14.00){\makebox(0,0)[cc]{$\partial$}}
\put(24.00,14.00){\makebox(0,0)[cc]{$=$}} \put(26.00,14.00){\line(1,0){16.00}}
\put(30.50,14.00){\oval(7.00,6.00)[t]} \put(37.50,14.00){\oval(7.00,6.00)[t]}
\put(47.00,14.00){\line(1,0){16.00}} \put(51.50,14.00){\oval(7.00,6.00)[t]}
\put(58.50,14.00){\oval(7.00,6.00)[t]} \put(68.00,14.00){\line(1,0){16.00}}
\put(72.50,14.00){\oval(7.00,6.00)[t]} \put(79.50,14.00){\oval(7.00,6.00)[t]}
\put(87.00,14.00){\makebox(0,0)[cc]{$+$}}
\put(66.00,14.00){\makebox(0,0)[cc]{$+$}}
\put(45.00,14.00){\makebox(0,0)[cc]{$+$}} \put(89.00,14.00){\line(1,0){16.00}}
\put(97.00,14.00){\oval(10.00,6.00)[t]}
\put(92.00,14.00){\makebox(0,0)[cc]{\large *}}
\put(108.00,14.00){\makebox(0,0)[cc]{$+$}}
\put(110.00,14.00){\line(1,0){15.00}} \put(117.50,14.00){\oval(9.00,6.00)[t]}
\put(97.00,8.00){\makebox(0,0)[cc]{\small $f_{1}'''(1)/f_{1}''(1) \sim$}}
\put(97.00,4.00){\makebox(0,0)[cc]{\small $\sim \ \mapsto /f_{1}''(1)$}}
\put(114.00,3.00){\line(0,1){4.00}} \put(114.00,5.00){\vector(1,0){10.00}}
\put(114.00,5.00){\vector(2,1){8.00}} \put(118.00,9.00){\makebox(0,0)[cc]{$1$}}
\put(122.00,1.00){\vector(-2,1){6.00}} \put(118.00,3.00){\vector(-2,1){4.00}}
\put(117.00,1.00){\makebox(0,0)[cc]{$2$}} \put(33.00,3.00){\line(0,1){4.00}}
\put(33.00,5.00){\vector(1,0){6.00}} \put(33.00,5.00){\vector(-1,0){6.00}}
\put(25.00,5.00){\makebox(0,0)[cc]{$1$}}
\put(41.00,5.00){\makebox(0,0)[cc]{$2$}} \put(53.00,5.00){\makebox(0,0)[cc]{$1
\mapsto$}} \put(74.00,5.00){\makebox(0,0)[cc]{$2 \mapsto$}}
\put(122.00,14.00){\makebox(0,0)[cc]{{\large *}}}
\end{picture}
 \ .
\end{equation}
The variety in its left part consists of such points of the cell \unitlength
1.00mm \special{em:linewidth 0.4pt} \linethickness{0.4pt}
\begin{picture}(18.00,4.00)
\put(0.00,0.00){\line(1,0){16.00}} \put(4.50,0.00){\oval(7.00,4.00)[t]}
\put(11.50,0.00){\oval(7.00,4.00)[t]}
\end{picture}
 \
that the direction "to the right" in $\R^{n-1}$ lies between the projections of
$f'(a_1)$ and $f'(a_2)$ to $\R^{n-1}$. The sum of three first terms in the
right-hand part of (\ref{nev3}) is equal to entire (\ref{homztrip2}). The
subscript under the fourth term in (\ref{nev3}) means almost the same as in
(\ref{nev2}) or (\ref{mar2}), but now the two frames compared there should
define {\em equal} orientations.

Finally, the last term in (\ref{nev3}) belongs to the 5th cell in (\ref{marg}).
This cell can be considered as the space of triples $(\alpha,t,f)$ where
$\alpha$ is a pair of points $(a<b)$ in $\R^1$, $f$ a map $\R^1 \to \R^n$ such
that $f(a)=f(b), f'(b)=0$, and $t$ is a point of a segment participating in the
construction of the simplicial resolution (its endpoints formally correspond to
the above two linear conditions). The subscript under the picture of this cell
in (\ref{nev3}) denotes a subvariety in the space of such triples, defined by
the following additional condition: the direction "to the right" in $\R^{n-1}$
belongs to the angle between projections of vectors $f'(a)$ and $-f''(b)$. Here
the number of arrows labeled by $2$ shows us the order of the derivative at the
second point $b$ participating in this condition, and the reversed direction of
these arrows indicates that we need to take this derivative with the opposite
sign.

Now we span the chain (\ref{homztrip3}) inside the cell \unitlength 1.00mm
\special{em:linewidth 0.4pt} \linethickness{0.4pt}
\begin{picture}(20.00,4.00)
\put(0.00,0.00){\line(1,0){20.00}} \put(10.00,0.00){\oval(16.00,8.00)[t]}
\put(14.00,0.00){\oval(8.00,4.00)[t]}
\end{picture}
 . First of all we kill the 5th picture in
(\ref{homztrip3}) by the variety shown in the left part of the next equality:
\begin{equation}
\label{span3} \unitlength 1.00mm \special{em:linewidth 0.4pt}
\linethickness{0.4pt}
\begin{picture}(109.00,26.00)
\put(3.00,22.00){\makebox(0,0)[cc]{$\partial$}}
\put(6.00,22.00){\line(1,0){21.00}} \put(18.50,22.00){\oval(13.00,8.00)[t]}
\put(21.50,22.00){\oval(7.00,4.00)[t]} \put(9.00,22.00){\line(2,-1){8.00}}
\put(17.00,18.00){\vector(1,1){4.00}} \put(33.00,22.00){\line(1,0){21.00}}
\put(45.50,22.00){\oval(13.00,8.00)[t]} \put(48.50,22.00){\oval(7.00,4.00)[t]}
\put(36.00,22.00){\line(2,-1){8.00}} \put(44.00,18.00){\line(1,1){4.00}}
\put(30.00,22.00){\makebox(0,0)[cc]{$=$}}
\put(57.00,22.00){\makebox(0,0)[cc]{$+$}} \put(60.00,22.00){\line(1,0){20.00}}
\put(71.00,22.00){\oval(14.00,8.00)[t]} \put(74.00,22.00){\oval(8.00,4.00)[t]}
\put(64.00,22.00){\line(3,-2){6.00}} \put(70.00,18.00){\vector(1,1){4.00}}
\put(83.00,22.00){\makebox(0,0)[cc]{$+$}} \put(86.00,22.00){\line(1,0){20.00}}
\put(97.00,22.00){\oval(14.00,8.00)[t]} \put(100.00,22.00){\oval(8.00,4.00)[t]}
\put(109.00,22.00){\makebox(0,0)[cc]{$+$}} \put(33.00,8.00){\line(1,0){20.00}}
\put(44.00,8.00){\oval(14.00,8.00)[t]} \put(47.00,8.00){\oval(8.00,4.00)[t]}
\put(56.00,8.00){\makebox(0,0)[cc]{$+$}} \put(60.00,8.00){\line(1,0){20.00}}
\put(71.00,8.00){\oval(14.00,8.00)[t]} \put(74.50,8.00){\oval(7.00,4.00)[t]}
\put(86.00,8.00){\line(1,0){20.00}} \put(97.00,8.00){\oval(14.00,8.00)[t]}
\put(100.50,8.00){\oval(7.00,4.00)[t]} \put(83.00,8.00){\makebox(0,0)[cc]{$+$}}
\put(92.00,18.00){\line(-5,4){5.00}} \put(92.00,18.00){\vector(1,1){4.00}}
\put(43.00,4.00){\line(-2,1){8.00}} \put(43.00,4.00){\vector(2,1){8.00}}
\put(71.00,6.00){\vector(-1,-2){3.00}} \put(71.00,6.00){\vector(1,-2){3.00}}
\put(76.00,3.00){\makebox(0,0)[cc]{1}} \put(66.00,3.00){\makebox(0,0)[cc]{2}}
\put(94.00,3.50){\vector(3,-1){6.00}} \put(94.00,3.50){\vector(3,1){6.00}}
\put(96.00,6.00){\makebox(0,0)[cc]{{\small 3}}}
\put(96.00,1.00){\makebox(0,0)[cc]{{\small 1}}}
\end{picture}
 ,
\end{equation}
thus reducing it to the sum of other five pictures in the right part of this
equality. The last picture in the upper row of (\ref{span3}) and the first
picture in the lower row denote one and the same set and annihilate. The first
term in (\ref{homztrip3}) together with the second from the end term in
(\ref{span3}) form a subvariety in the same cell defined by the condition of
the type \unitlength 1.00mm \special{em:linewidth 0.4pt} \linethickness{0.4pt}
\begin{picture}(11.50,3.00)
\put(6.00,1.00){\makebox(0,0)[cc]{$1 \longleftrightarrow 2$}}
\put(6.00,-0.50){\line(0,1){3.00}}
\end{picture}
 \ .
It is natural to kill it by the left part of the following equation:
\begin{equation}
\label{span4} \unitlength 1.00mm \linethickness{0.4pt}
\begin{picture}(125.00,15.00)
\put(3.00,11.00){\makebox(0,0)[cc]{$\partial$}}
\put(6.00,11.00){\line(1,0){19.00}} \put(15.50,11.00){\oval(15.00,8.00)[t]}
\put(19.00,11.00){\oval(8.00,4.00)[t]} \put(31.00,11.00){\line(1,0){19.00}}
\put(40.50,11.00){\oval(15.00,8.00)[t]} \put(44.00,11.00){\oval(8.00,4.00)[t]}
\put(56.00,11.00){\line(1,0){19.00}} \put(65.50,11.00){\oval(15.00,8.00)[t]}
\put(69.00,11.00){\oval(8.00,4.00)[t]} \put(81.00,11.00){\line(1,0){19.00}}
\put(90.50,11.00){\oval(15.00,8.00)[t]} \put(94.00,11.00){\oval(8.00,4.00)[t]}
\put(106.00,11.00){\line(1,0){19.00}}
\put(103.00,11.00){\makebox(0,0)[cc]{$+$}}
\put(78.00,11.00){\makebox(0,0)[cc]{$+$}}
\put(53.00,11.00){\makebox(0,0)[cc]{$+$}}
\put(28.00,11.00){\makebox(0,0)[cc]{=}} \put(12.00,4.00){\vector(2,1){6.00}}
\put(12.00,4.00){\vector(2,-1){6.00}} \put(12.00,4.00){\vector(1,0){10.00}}
\put(12.00,2.00){\line(0,1){4.00}} \put(14.00,7.00){\makebox(0,0)[cc]{1}}
\put(14.00,1.00){\makebox(0,0)[cc]{2}} \put(62.00,6.00){\makebox(0,0)[cc]{$1
\mapsto$}} \put(90.00,6.00){\makebox(0,0)[cc]{$2 \mapsto$}}
\put(40.00,6.00){\line(0,-1){4.00}} \put(46.00,4.00){\vector(1,0){0.2}}
\put(40.00,4.00){\line(1,0){6.00}} \put(34.00,4.00){\vector(-1,0){0.2}}
\put(40.00,4.00){\line(-1,0){6.00}} \put(31.00,4.00){\makebox(0,0)[cc]{$1$}}
\put(49.00,4.00){\makebox(0,0)[cc]{$2$}}
\put(115.50,11.00){\oval(13.00,8.00)[t]}
\put(122.00,11.00){\makebox(0,0)[cc]{$*$}} \put(112.00,3.00){\line(0,1){4.00}}
\put(112.00,5.00){\vector(1,0){10.00}} \put(112.00,5.00){\vector(2,1){8.00}}
\put(116.00,9.00){\makebox(0,0)[cc]{$1$}}
\put(120.00,1.00){\vector(-2,1){4.00}} \put(116.00,3.00){\vector(-2,1){4.00}}
\put(115.00,1.00){\makebox(0,0)[cc]{$2$}}
\end{picture}
 \ .
\end{equation}
Summing up all terms in right-hand parts of equations
(\ref{nev3})--(\ref{span4}) and subtracting the chains (\ref{homztrip2}),
(\ref{homztrip3}), (\ref{mar2}), we annihilate almost all of their summands
except for the term (\ref{mar2}) and the second from the right term of
(\ref{nev3}). The sum of these two terms is equal to the right-hand part of the
identity
\begin{equation}
\label{span5} \unitlength 1.00mm \special{em:linewidth 0.4pt}
\linethickness{0.4pt}
\begin{picture}(72.00,10.00)
\put(3.00,6.00){\makebox(0,0)[cc]{$\partial$}}
\put(6.00,6.00){\line(1,0){24.00}} \put(18.00,6.00){\oval(18.00,8.00)[t]}
\put(12.50,6.00){\oval(7.00,4.00)[t]} \put(39.00,6.00){\makebox(0,0)[cc]{$=$}}
\put(48.00,6.00){\line(1,0){24.00}} \put(59.50,6.00){\oval(17.00,8.00)[t]}
\put(51.00,6.00){\makebox(0,0)[cc]{$*$}}
\put(6.00,0.00){\framebox(25.00,4.00)[cc]{\small
$f_{1}'(1),f_{1}'(2),\mapsto$}}
\put(46.00,0.00){\framebox(28.00,4.00)[cc]{\small
$f_{1}''(1),f_{1}'''(1),\mapsto$}}
\end{picture}
 \ \ .
\end{equation}
The subscript under this right-hand part means that the projections of
$f''(a_1)$ and $f'''(a_1)$ to $\R^{n-1}$ and the direction "to the right"
should be linearly dependent; the subscript in the left part says the same
about projections of vectors $f'(a_1)$ and $f'(a_2)$. If $n=3$ then both these
subscripts mean nothing.

Summarizing, we get that for the desired chain spanning (\ref{difone}) in $F_2
\setminus F_1$ we can take the sum of varieties shown in left parts of
equalities (\ref{boundd})--(\ref{nev2}) and (\ref{nev3})--(\ref{span5}), i.e.
the chain
\begin{equation}
\label{spanlast} \unitlength=1.00mm \special{em:linewidth 0.4pt}
\linethickness{0.4pt}
\begin{picture}(124.00,30.00)
\put(0.00,26.00){\line(1,0){19.00}} \put(9.50,26.00){\oval(15.00,8.00)[t]}
\put(6.00,26.00){\oval(8.00,4.00)[t]} \put(13.50,26.00){\oval(7.00,4.00)[t]}
\put(8.00,18.00){\vector(1,2){3.00}} \put(8.00,18.00){\vector(-1,2){3.00}}
\put(3.00,21.00){\makebox(0,0)[cc]{2}} \put(13.00,21.00){\makebox(0,0)[cc]{1}}
\put(22.00,26.00){\makebox(0,0)[cc]{$+$}} \put(25.00,26.00){\line(1,0){19.00}}
\put(34.50,26.00){\oval(15.00,8.00)[t]} \put(31.00,26.00){\oval(8.00,4.00)[t]}
\put(38.50,26.00){\oval(7.00,4.00)[t]}
\put(47.00,26.00){\makebox(0,0)[cc]{$+$}} \put(50.00,26.00){\line(1,0){19.00}}
\put(59.50,26.00){\oval(15.00,8.00)[t]} \put(56.00,26.00){\oval(8.00,4.00)[t]}
\put(63.50,26.00){\oval(7.00,4.00)[t]}
\put(72.00,26.00){\makebox(0,0)[cc]{$+$}} \put(35.00,26.00){\line(-1,-2){2.00}}
\put(33.00,22.00){\vector(-1,2){2.00}} \put(57.00,21.00){\makebox(0,0)[cc]{$1
\mapsto$}} \put(75.00,26.00){\line(1,0){19.00}}
\put(82.50,26.00){\oval(11.00,8.00)[b]} \put(86.50,26.00){\oval(11.00,8.00)[t]}
\put(81.00,26.00){\line(-1,4){1.00}} \put(80.00,30.00){\vector(-1,-4){1.00}}
\put(97.00,26.00){\makebox(0,0)[cc]{$+$}} \put(100.00,26.00){\line(1,0){19.00}}
\put(122.00,26.00){\makebox(0,0)[cc]{$+$}} \put(3.00,10.00){\line(1,0){19.00}}
\put(1.00,10.00){\makebox(0,0)[cc]{$+$}}
\put(10.50,10.00){\oval(11.00,8.00)[b]} \put(14.50,10.00){\oval(11.00,8.00)[t]}
\put(12.00,2.00){\makebox(0,0)[cc]{$2 \mapsto$}}
\put(25.00,10.00){\makebox(0,0)[cc]{$+$}} \put(28.00,10.00){\line(1,0){19.00}}
\put(34.00,10.00){\oval(8.00,6.00)[t]} \put(42.00,10.00){\oval(8.00,6.00)[t]}
\put(35.00,4.00){\vector(2,1){6.00}} \put(35.00,4.00){\vector(2,-1){6.00}}
\put(35.00,4.00){\vector(1,0){10.00}} \put(35.00,2.00){\line(0,1){4.00}}
\put(37.00,0.00){\makebox(0,0)[cc]{2}} \put(37.00,8.00){\makebox(0,0)[cc]{1}}
\put(50.00,10.00){\makebox(0,0)[cc]{$+$}} \put(53.00,10.00){\line(1,0){19.00}}
\put(63.50,10.00){\oval(13.00,8.00)[t]} \put(66.50,10.00){\oval(7.00,6.00)[t]}
\put(55.00,10.00){\line(3,-2){6.00}} \put(61.00,6.00){\vector(3,2){6.00}}
\put(75.00,10.00){\makebox(0,0)[cc]{$+$}} \put(78.00,10.00){\line(1,0){19.00}}
\put(88.00,10.00){\oval(14.00,8.00)[t]} \put(91.50,10.00){\oval(7.00,6.00)[t]}
\put(85.00,4.00){\vector(2,1){6.00}} \put(85.00,4.00){\vector(2,-1){6.00}}
\put(85.00,4.00){\vector(1,0){10.00}} \put(85.00,2.00){\line(0,1){4.00}}
\put(87.00,0.00){\makebox(0,0)[cc]{2}} \put(87.00,8.00){\makebox(0,0)[cc]{1}}
\put(108.50,26.00){\oval(15.00,8.00)[b]}
\put(111.00,26.00){\oval(14.00,8.00)[t]} \put(104.00,26.00){\line(5,3){5.00}}
\put(109.00,29.00){\vector(4,-3){4.00}} \put(103.00,10.00){\line(1,0){21.00}}
\put(113.50,10.00){\oval(15.00,8.00)[t]}
\put(109.50,10.00){\oval(7.00,4.00)[t]}
\put(100.00,10.00){\makebox(0,0)[cc]{$+$}}
\put(101.00,2.00){\framebox(25.00,5.00)[cc]{\small
$f_{1}'(1),f_{1}'(2),\mapsto$}}
\end{picture}
 \ .
\end{equation}

\subsection{The second differential and its homology to zero}

Now let us consider the boundary of the chain (\ref{spanlast}) in the term
$F_1$ of the filtration. This term consists of two cells, one of which is
characterized by a single chord and the second by one asterisk; see
(\ref{dif01}). It is easy to see that the first three summands in
(\ref{spanlast}) do not have any homological boundary in these cells, and the
next seven have two components of the boundary each, and these pairs of
components are shown consecutively in the next formula (\ref{diftwo}):
\begin{equation}
\label{diftwo} \unitlength 1.00mm \linethickness{0.4pt}
\begin{picture}(120.00,77.00)
\put(0.00,72.00){\line(1,0){25.00}} \put(10.00,72.00){\oval(16.00,10.00)[b]}
\put(23.00,72.00){\line(-3,2){6.00}} \put(17.00,76.00){\line(-1,-1){4.00}}
\put(13.00,72.00){\line(-4,5){4.00}} \put(9.00,77.00){\vector(-3,-4){3.67}}
\put(28.00,72.00){\makebox(0,0)[cc]{$+$}} \put(31.00,72.00){\line(1,0){25.00}}
\put(47.50,72.00){\oval(13.00,10.00)[b]} \put(33.00,72.00){\line(2,1){10.00}}
\put(43.00,77.00){\line(6,-5){6.00}} \put(41.00,72.00){\line(-2,-5){2.00}}
\put(39.00,67.00){\vector(-2,3){3.33}}
\put(59.00,72.00){\makebox(0,0)[cc]{$+$}} \put(62.00,72.00){\line(1,0){25.00}}
\put(72.00,72.00){\oval(16.00,10.00)[b]} \put(85.00,72.00){\line(-2,1){10.00}}
\put(75.00,77.00){\line(-6,-5){6.00}} \put(69.00,72.00){\line(1,-2){2.00}}
\put(71.00,68.00){\vector(1,1){4.00}} \put(90.00,72.00){\makebox(0,0)[cc]{$+$}}
\put(93.00,72.00){\line(1,0){24.00}} \put(108.00,72.00){\oval(14.00,10.00)[b]}
\put(95.00,72.00){\line(2,1){10.00}} \put(105.00,77.00){\line(6,-5){6.00}}
\put(101.00,72.00){\line(1,-1){4.00}} \put(105.00,68.00){\vector(1,2){2.00}}
\put(120.00,72.00){\makebox(0,0)[cc]{$+$}} \put(6.00,54.00){\line(1,0){23.00}}
\put(14.50,54.00){\oval(15.00,10.00)[b]} \put(28.00,54.00){\line(-2,1){10.00}}
\put(18.00,59.00){\line(-6,-5){6.00}} \put(13.00,45.00){\makebox(0,0)[cc]{$2
\mapsto$}} \put(32.00,54.00){\makebox(0,0)[cc]{$+$}}
\put(35.00,54.00){\line(1,0){23.00}} \put(49.50,54.00){\oval(15.00,10.00)[b]}
\put(36.00,54.00){\line(2,1){10.00}} \put(46.00,59.00){\line(6,-5){6.00}}
\put(42.00,45.00){\makebox(0,0)[cc]{$2 \mapsto$}}
\put(61.00,54.00){\makebox(0,0)[cc]{$+$}} \put(65.00,54.00){\line(1,0){23.00}}
\put(71.50,54.00){\oval(11.00,10.00)[b]} \put(77.00,54.00){\line(1,-1){5.00}}
\put(82.00,49.00){\line(1,1){5.00}} \put(68.00,43.00){\vector(2,1){6.00}}
\put(68.00,43.00){\vector(2,-1){6.00}} \put(68.00,43.00){\vector(1,0){10.00}}
\put(68.00,41.00){\line(0,1){4.00}} \put(91.00,54.00){\makebox(0,0)[cc]{$+$}}
\put(94.00,54.00){\line(1,0){23.00}} \put(109.00,54.00){\oval(12.00,10.00)[b]}
\put(103.00,54.00){\line(-4,-5){4.00}} \put(99.00,49.00){\line(-4,5){4.00}}
\put(96.00,43.00){\vector(2,1){6.00}} \put(96.00,43.00){\vector(2,-1){6.00}}
\put(96.00,43.00){\vector(1,0){10.00}} \put(96.00,41.00){\line(0,1){4.00}}
\put(120.00,54.00){\makebox(0,0)[cc]{$+$}} \put(5.00,30.00){\line(1,0){25.00}}
\put(32.00,30.00){\makebox(0,0)[cc]{$+$}} \put(35.00,30.00){\line(1,0){23.00}}
\put(61.00,30.00){\makebox(0,0)[cc]{$+$}} \put(65.00,30.00){\line(1,0){23.00}}
\put(76.50,30.00){\oval(21.00,10.00)[b]} \put(87.00,30.00){\line(-6,5){6.00}}
\put(81.00,35.00){\line(-6,-5){6.00}} \put(91.00,30.00){\makebox(0,0)[cc]{$+$}}
\put(94.00,30.00){\line(1,0){23.00}} \put(110.00,30.00){\oval(10.00,10.00)[b]}
\put(115.00,30.00){\line(-2,1){12.00}} \put(103.00,36.00){\line(-4,-3){8.00}}
\put(77.00,19.00){\vector(2,1){6.00}} \put(77.00,19.00){\vector(2,-1){6.00}}
\put(77.00,19.00){\vector(1,0){10.00}} \put(77.00,17.00){\line(0,1){4.00}}
\put(96.00,19.00){\vector(2,1){6.00}} \put(96.00,19.00){\vector(2,-1){6.00}}
\put(96.00,19.00){\vector(1,0){10.00}} \put(96.00,17.00){\line(0,1){4.00}}
\put(70.00,46.00){\makebox(0,0)[cc]{1}} \put(70.00,39.00){\makebox(0,0)[cc]{2}}
\put(96.00,41.00){\line(0,1){4.00}} \put(98.00,46.00){\makebox(0,0)[cc]{1}}
\put(98.00,39.00){\makebox(0,0)[cc]{2}} \put(77.00,17.00){\line(0,1){4.00}}
\put(79.00,22.00){\makebox(0,0)[cc]{1}} \put(79.00,15.00){\makebox(0,0)[cc]{2}}
\put(96.00,17.00){\line(0,1){4.00}} \put(96.00,17.00){\line(0,1){4.00}}
\put(98.00,22.00){\makebox(0,0)[cc]{1}} \put(98.00,15.00){\makebox(0,0)[cc]{2}}
\put(6.00,30.00){\line(0,1){0.00}} \put(6.00,30.00){\line(1,-1){6.00}}
\put(12.00,24.00){\line(1,0){6.00}} \put(18.00,24.00){\vector(1,1){6.00}}
\put(20.00,30.00){\oval(18.00,10.00)[t]} \put(29.00,30.00){\line(-3,2){6.00}}
\put(23.00,34.00){\line(-3,-2){6.00}} \put(36.00,30.00){\line(1,-1){6.00}}
\put(42.00,24.00){\line(1,0){6.00}} \put(51.00,30.00){\oval(6.00,8.00)[t]}
\put(48.00,30.00){\line(-1,1){4.00}} \put(44.00,34.00){\line(-3,-4){3.00}}
\put(120.00,30.00){\makebox(0,0)[cc]{$+$}}
\put(31.00,10.00){\line(-1,0){25.00}} \put(18.50,10.00){\oval(19.00,10.00)[t]}
\put(36.00,10.00){\makebox(0,0)[cc]{$+$}} \put(41.00,10.00){\line(1,0){25.00}}
\put(48.00,10.00){\oval(8.00,6.00)[b]} \put(9.00,10.00){\line(2,-1){6.00}}
\put(15.00,7.00){\line(2,1){6.00}} \put(44.00,10.00){\line(2,1){10.00}}
\put(54.00,15.00){\line(2,-1){10.00}} \put(48.00,24.00){\vector(1,2){3.00}}
\put(6.00,0.00){\framebox(25.00,5.00)[cc]{\small $f_1'(1),f_1'(2),\mapsto$}}
\put(41.00,0.00){\framebox(25.00,5.00)[cc]{\small $f_1'(1),f_1'(2),\mapsto$}}
\put(3.00,10.00){\makebox(0,0)[cc]{$+$}}
\put(2.00,54.00){\makebox(0,0)[cc]{$+$}}
\put(2.00,30.00){\makebox(0,0)[cc]{$+$}}
\end{picture}
 \ .
\end{equation}

The last pictures in the second and the third lines of this formula denote one
and the same variety and annihilate, so we get only the sum of remaining twelve
varieties.

Now let us span this sum by a chain in $F_1$. As usual, any time as we have a
variety characterized by a picture with a zigzag (without arrows) we represent
it as a component of the boundary of a variety, whose picture is obtained from
this one by adding an arrow at one of the ends of the zigzag. Performing this
systematically, we find some ten varieties of codimension 1 in the greater cell
of $F_1$. These varieties are encoded by the left parts of the following
equations (\ref{d2-1})--(\ref{d2-91}), whose right-hand parts express the
boundaries of these varieties.
\begin{equation}
\label{d2-1} \unitlength=1.00mm \special{em:linewidth 0.4pt}
\linethickness{0.4pt}
\begin{picture}(124.00,29.00)
\put(3.00,24.00){\makebox(0,0)[cc]{$\partial$}}
\put(6.00,24.00){\line(1,0){25.00}} \put(16.00,24.00){\oval(16.00,10.00)[b]}
\put(29.00,24.00){\line(-1,1){5.00}} \put(24.00,29.00){\vector(-1,-1){5.00}}
\put(19.00,24.00){\line(-3,5){3.00}} \put(16.00,29.00){\vector(-1,-1){5.00}}
\put(34.00,24.00){\makebox(0,0)[cc]{$=$}} \put(37.00,24.00){\line(1,0){25.00}}
\put(46.50,24.00){\oval(15.00,10.00)[b]} \put(60.00,24.00){\line(-1,1){5.00}}
\put(55.00,29.00){\line(-1,-1){5.00}} \put(50.00,24.00){\line(-3,5){3.00}}
\put(47.00,29.00){\vector(-1,-1){5.00}}
\put(65.00,24.00){\makebox(0,0)[cc]{$+$}} \put(68.00,24.00){\line(1,0){25.00}}
\put(77.50,24.00){\oval(15.00,10.00)[b]} \put(91.00,24.00){\line(-1,1){5.00}}
\put(86.00,29.00){\vector(-1,-1){5.00}} \put(81.00,24.00){\line(-4,5){4.00}}
\put(77.00,29.00){\line(-4,-5){4.00}} \put(96.00,24.00){\makebox(0,0)[cc]{$+$}}
\put(99.00,24.00){\line(1,0){22.00}} \put(107.50,24.00){\oval(13.00,10.00)[b]}
\put(119.00,24.00){\line(-1,1){5.00}} \put(114.00,29.00){\vector(-1,-1){5.00}}
\put(109.00,24.00){\line(-3,5){3.00}} \put(106.00,29.00){\vector(-1,-1){5.00}}
\put(124.00,24.00){\makebox(0,0)[cc]{$+$}}
\put(3.00,8.00){\makebox(0,0)[cc]{$+$}} \put(36.00,8.00){\line(1,0){25.00}}
\put(45.50,8.00){\oval(15.00,10.00)[b]} \put(60.00,8.00){\line(-3,4){3.00}}
\put(57.00,12.00){\vector(-1,-1){4.00}} \put(38.00,8.00){\line(1,1){5.00}}
\put(43.00,13.00){\vector(1,-1){5.00}} \put(64.00,8.00){\makebox(0,0)[cc]{$+$}}
\put(67.00,8.00){\line(1,0){24.00}} \put(79.00,8.00){\oval(20.00,10.00)[b]}
\put(89.00,8.00){\line(-3,5){3.00}} \put(86.00,13.00){\vector(-1,-1){5.00}}
\put(81.00,8.00){\line(-3,5){3.00}} \put(78.00,13.00){\vector(-1,-1){5.00}}
\put(94.00,8.00){\makebox(0,0)[cc]{$+$}} \put(7.00,8.00){\line(1,0){24.00}}
\put(16.00,8.00){\oval(14.00,10.00)[b]} \put(28.00,8.00){\line(-1,1){5.00}}
\put(23.00,13.00){\vector(-1,-1){5.00}} \put(27.00,2.00){\makebox(0,0)[cc]{$2
\uparrow$}} \put(34.00,8.00){\makebox(0,0)[cc]{$+$}}
\put(97.00,8.00){\line(1,0){24.00}} \put(106.00,8.00){\oval(14.00,10.00)[b]}
\put(119.00,8.00){\line(-3,4){3.00}} \put(116.00,12.00){\vector(-3,-4){3.00}}
\put(119.00,6.00){\vector(1,-2){3.00}} \put(119.00,6.00){\vector(-1,-2){3.00}}
\put(115.00,3.00){\makebox(0,0)[cc]{1}} \put(123.00,3.00){\makebox(0,0)[cc]{2}}
\end{picture}
\end{equation}

\begin{equation}
\label{d2-2} \unitlength=1.00mm \special{em:linewidth 0.4pt}
\linethickness{0.4pt}
\begin{picture}(124.00,29.00)
\put(3.00,24.00){\makebox(0,0)[cc]{$\partial$}}
\put(6.00,24.00){\line(1,0){25.00}} \put(15.50,24.00){\oval(15.00,10.00)[b]}
\put(29.00,24.00){\line(-1,1){5.00}} \put(24.00,29.00){\line(-1,0){11.00}}
\put(13.00,29.00){\vector(0,-1){5.00}} \put(13.00,24.00){\line(1,2){2.00}}
\put(15.00,28.00){\vector(3,-4){3.00}}
\put(34.00,24.00){\makebox(0,0)[cc]{$=$}} \put(37.00,24.00){\line(1,0){25.00}}
\put(46.50,24.00){\oval(15.00,10.00)[b]} \put(60.00,24.00){\line(-1,1){5.00}}
\put(55.00,29.00){\vector(-1,-1){5.00}} \put(50.00,24.00){\line(-3,5){3.00}}
\put(47.00,29.00){\line(-3,-5){3.00}} \put(65.00,24.00){\makebox(0,0)[cc]{$+$}}
\put(68.00,24.00){\line(1,0){25.00}} \put(77.50,24.00){\oval(15.00,10.00)[b]}
\put(91.00,24.00){\line(-3,5){3.00}} \put(88.00,29.00){\line(-1,0){14.00}}
\put(74.00,29.00){\line(0,-1){5.00}} \put(74.00,24.00){\line(3,4){3.00}}
\put(77.00,28.00){\vector(3,-4){3.00}}
\put(96.00,24.00){\makebox(0,0)[cc]{$+$}} \put(99.00,24.00){\line(1,0){22.00}}
\put(107.50,24.00){\oval(13.00,10.00)[b]} \put(119.00,24.00){\line(-1,1){5.00}}
\put(114.00,29.00){\vector(0,-1){5.00}} \put(101.00,24.00){\line(3,5){3.00}}
\put(104.00,29.00){\vector(1,-1){5.00}}
\put(124.00,24.00){\makebox(0,0)[cc]{$+$}}
\put(27.00,9.00){\makebox(0,0)[cc]{$+$}} \put(30.00,9.00){\line(1,0){19.00}}
\put(38.00,9.00){\oval(12.00,10.00)[b]} \put(48.00,9.00){\line(-1,1){5.00}}
\put(43.00,14.00){\vector(-1,-1){5.00}} \put(38.00,9.00){\line(-3,4){3.00}}
\put(35.00,13.00){\vector(-3,-4){3.00}} \put(2.00,9.00){\makebox(0,0)[cc]{$+$}}
\put(5.00,9.00){\line(1,0){18.00}} \put(12.00,9.00){\oval(12.00,10.00)[b]}
\put(22.00,9.00){\line(-1,1){5.00}} \put(17.00,14.00){\vector(-1,-1){5.00}}
\put(21.00,2.00){\makebox(0,0)[cc]{$2 \downarrow$}}
\put(51.00,9.00){\makebox(0,0)[cc]{$+$}} \put(54.00,9.00){\line(1,0){17.00}}
\put(62.50,9.00){\oval(15.00,10.00)[b]} \put(60.00,9.00){\line(1,1){5.00}}
\put(65.00,14.00){\vector(0,-1){5.00}} \put(55.00,9.00){\line(1,1){5.00}}
\put(60.00,14.00){\vector(0,-1){5.00}} \put(74.00,9.00){\makebox(0,0)[cc]{$+$}}
\put(77.00,9.00){\line(1,0){17.00}} \put(82.00,9.00){\oval(8.00,10.00)[b]}
\put(93.00,9.00){\line(-2,5){2.00}} \put(91.00,14.00){\vector(-1,-1){5.00}}
\put(91.00,-1.00){\vector(-1,2){3.00}} \put(91.00,-1.00){\vector(1,2){3.00}}
\put(95.00,2.00){\makebox(0,0)[cc]{2}} \put(87.00,2.00){\makebox(0,0)[cc]{1}}
\put(97.00,9.00){\makebox(0,0)[cc]{$+$}} \put(100.00,9.00){\line(1,0){16.00}}
\put(108.00,9.00){\oval(14.00,10.00)[b]} \put(115.00,9.00){\line(-3,5){3.00}}
\put(112.00,14.00){\vector(-1,-1){5.00}} \put(117.00,3.00){\vector(2,1){6.00}}
\put(117.00,3.00){\vector(2,-1){6.00}}
\put(119.00,-1.00){\makebox(0,0)[cc]{$1$}}
\put(119.00,7.00){\makebox(0,0)[cc]{$2$}}
\end{picture}
\end{equation}

\begin{equation}
\label{d2-3} \unitlength=1.00mm \special{em:linewidth 0.4pt}
\linethickness{0.4pt}
\begin{picture}(125.00,29.00)
\put(2.00,24.00){\makebox(0,0)[cc]{$\partial$}}
\put(5.00,24.00){\line(1,0){25.00}} \put(14.50,24.00){\oval(15.00,8.00)[b]}
\put(28.00,24.00){\line(-6,5){6.00}} \put(22.00,29.00){\vector(-4,-3){6.67}}
\put(26.00,19.00){\makebox(0,0)[cc]{$2 \mapsto$}}
\put(33.00,24.00){\makebox(0,0)[cc]{$=$}} \put(36.00,24.00){\line(1,0){25.00}}
\put(45.00,24.00){\oval(14.00,8.00)[b]} \put(59.00,24.00){\line(-5,3){8.00}}
\put(51.00,28.67){\line(-5,-4){6.00}} \put(56.00,19.00){\makebox(0,0)[cc]{$2
\mapsto$}} \put(64.00,24.00){\makebox(0,0)[cc]{$+$}}
\put(67.00,24.00){\line(1,0){25.00}} \put(76.00,24.00){\oval(14.00,8.00)[b]}
\put(89.00,24.00){\line(-6,5){6.00}} \put(83.00,29.00){\vector(-4,-3){6.67}}
\put(87.00,19.00){\makebox(0,0)[cc]{$2 \updownarrow$}}
\put(95.00,24.00){\makebox(0,0)[cc]{$+$}} \put(98.00,24.00){\line(1,0){23.00}}
\put(106.50,24.00){\oval(13.00,8.00)[b]}
\put(117.00,19.00){\makebox(0,0)[cc]{$1 \mapsto$}}
\put(124.00,24.00){\makebox(0,0)[cc]{$+$}}
\put(1.00,8.00){\makebox(0,0)[cc]{$+$}} \put(4.00,8.00){\line(1,0){25.00}}
\put(12.00,8.00){\oval(12.00,8.00)[b]} \put(27.00,8.00){\line(-4,5){4.00}}
\put(23.00,13.00){\vector(-1,-1){5.00}} \put(18.00,0.00){\makebox(0,0)[cc]{$2
\mapsto$}} \put(32.00,8.00){\makebox(0,0)[cc]{$+$}}
\put(35.00,8.00){\line(1,0){24.00}} \put(47.00,8.00){\oval(20.00,8.00)[b]}
\put(57.00,8.00){\line(-1,1){5.00}} \put(52.00,13.00){\vector(-1,-1){5.00}}
\put(47.00,0.00){\makebox(0,0)[cc]{$2 \mapsto$}}
\put(62.00,8.00){\makebox(0,0)[cc]{$+$}} \put(65.00,8.00){\line(1,0){24.00}}
\put(76.00,8.00){\oval(16.00,8.00)[b]} \put(70.00,0.00){\makebox(0,0)[cc]{$1
\mapsto$}} \put(93.00,8.00){\makebox(0,0)[cc]{$+$}}
\put(96.00,8.00){\line(1,0){29.00}} \put(102.00,8.00){\makebox(0,0)[cc]{$*$}}
\put(122.00,8.00){\line(-2,1){10.00}} \put(112.00,13.00){\vector(-2,-1){10.00}}
\put(112.00,4.00){\makebox(0,0)[cc]{\small $f_1'''(1)/f_1''(1) \sim$}}
\put(112.00,0.00){\makebox(0,0)[cc]{\small $\sim \ \mapsto/ f_1''(1)$}}
\put(84.00,2.00){\vector(3,1){6.00}} \put(84.00,2.00){\vector(3,-1){6.00}}
\put(86.00,-1.00){\makebox(0,0)[cc]{$1$}}
\put(86.00,5.00){\makebox(0,0)[cc]{$2$}} \put(120.00,24.00){\line(-2,1){10.00}}
\put(110.00,29.00){\vector(-2,-1){10.00}}
\end{picture}
\end{equation}

\begin{equation}
\label{d2-4} \unitlength=1.00mm \special{em:linewidth 0.4pt}
\linethickness{0.4pt}
\begin{picture}(125.00,38.00)
\put(3.00,33.00){\makebox(0,0)[cc]{$\partial$}}
\put(6.00,33.00){\line(1,0){25.00}} \put(21.00,25.00){\vector(2,1){6.00}}
\put(21.00,25.00){\vector(2,-1){6.00}} \put(21.00,25.00){\vector(1,0){10.00}}
\put(21.00,23.00){\line(0,1){4.00}} \put(23.00,21.00){\makebox(0,0)[cc]{2}}
\put(23.00,29.00){\makebox(0,0)[cc]{1}}
\put(34.00,33.00){\makebox(0,0)[cc]{$=$}} \put(37.00,33.00){\line(1,0){25.00}}
\put(44.50,33.00){\oval(11.00,10.00)[b]} \put(60.00,33.00){\line(-1,1){5.00}}
\put(55.00,38.00){\line(-1,-1){5.00}} \put(13.50,33.00){\oval(11.00,10.00)[b]}
\put(29.00,33.00){\line(-1,1){5.00}} \put(24.00,38.00){\vector(-1,-1){5.00}}
\put(65.00,33.00){\makebox(0,0)[cc]{$+$}} \put(68.00,33.00){\line(1,0){25.00}}
\put(91.00,33.00){\line(-1,1){5.00}} \put(86.00,38.00){\vector(-1,-1){5.00}}
\put(75.50,33.00){\oval(11.00,10.00)[b]}
\put(96.00,33.00){\makebox(0,0)[cc]{$+$}} \put(52.00,25.00){\vector(2,1){6.00}}
\put(52.00,25.00){\vector(2,-1){6.00}} \put(52.00,25.00){\vector(1,0){10.00}}
\put(52.00,23.00){\line(0,1){4.00}} \put(54.00,21.00){\makebox(0,0)[cc]{2}}
\put(54.00,29.00){\makebox(0,0)[cc]{1}} \put(86.00,25.00){\makebox(0,0)[cc]{$1
\mapsto$}} \put(18.00,10.00){\makebox(0,0)[cc]{$+$}}
\put(98.00,33.00){\line(1,0){24.00}} \put(121.00,33.00){\line(-1,1){5.00}}
\put(116.00,38.00){\vector(-1,-1){5.00}}
\put(105.50,33.00){\oval(11.00,10.00)[b]}
\put(125.00,33.00){\makebox(0,0)[cc]{$+$}}
\put(116.00,25.00){\makebox(0,0)[cc]{$2 \mapsto$}}
\put(22.00,10.00){\line(1,0){25.00}} \put(29.50,10.00){\oval(11.00,10.00)[b]}
\put(45.00,10.00){\line(-1,1){5.00}} \put(40.00,15.00){\vector(-1,-1){5.00}}
\put(38.00,2.00){\makebox(0,0)[cc]{$1 \longleftrightarrow 2$}}
\put(38.00,0.50){\line(0,1){3.00}} \put(85.00,10.00){\makebox(0,0)[cc]{$+$}}
\put(88.00,10.00){\line(1,0){25.00}} \put(100.50,10.00){\oval(21.00,8.00)[t]}
\put(110.00,2.00){\makebox(0,0)[cc]{$2 \uparrow$}}
\put(51.00,10.00){\makebox(0,0)[cc]{$+$}} \put(54.00,10.00){\line(1,0){26.00}}
\put(77.00,10.00){\line(-2,1){10.00}} \put(67.00,15.00){\vector(-2,-1){10.00}}
\put(57.00,10.00){\makebox(0,0)[cc]{{\large $*$}}}
\put(67.00,5.00){\makebox(0,0)[cc]{\small $f_1'''(1)/f_1''(1) \sim$}}
\put(67.00,1.00){\makebox(0,0)[cc]{\small $\sim \ \mapsto/ f_1''(1)$}}
\put(90.00,5.00){\vector(1,0){12.00}} \put(90.00,5.00){\vector(3,1){9.00}}
\put(99.00,2.00){\vector(-3,1){6.00}} \put(96.00,3.00){\vector(-3,1){6.00}}
\put(90.00,2.00){\line(0,1){6.00}} \put(93.00,1.00){\makebox(0,0)[cc]{$2$}}
\put(93.00,8.00){\makebox(0,0)[cc]{$1$}}
\end{picture}
\end{equation}

\begin{equation}
\label{d2-5} \unitlength 1.00mm \linethickness{0.4pt}
\begin{picture}(126.00,42.00)
\put(3.00,37.00){\makebox(0,0)[cc]{$\partial$}}
\put(6.00,37.00){\line(1,0){25.00}} \put(18.50,37.00){\oval(21.00,8.00)[b]}
\put(29.00,37.00){\line(-6,5){6.00}} \put(23.00,42.00){\vector(-1,-1){5.00}}
\put(15.00,27.00){\vector(2,1){6.00}} \put(15.00,27.00){\vector(2,-1){6.00}}
\put(15.00,27.00){\vector(1,0){10.00}} \put(15.00,25.00){\line(0,1){4.00}}
\put(16.50,30.20){\makebox(0,0)[cc]{1}} \put(16.50,23.50){\makebox(0,0)[cc]{2}}
\put(34.00,37.00){\makebox(0,0)[cc]{$=$}} \put(37.00,37.00){\line(1,0){24.00}}
\put(49.00,37.00){\oval(20.00,8.00)[b]} \put(59.00,37.00){\line(-6,5){6.00}}
\put(53.00,42.00){\vector(-1,-1){5.00}} \put(48.00,29.00){\makebox(0,0)[cc]{$1
\mapsto$}} \put(68.00,37.00){\line(1,0){24.00}}
\put(80.00,37.00){\oval(20.00,8.00)[b]} \put(90.00,37.00){\line(-6,5){6.00}}
\put(84.00,42.00){\vector(-1,-1){5.00}} \put(79.00,29.00){\makebox(0,0)[cc]{$2
\mapsto$}} \put(95.00,37.00){\makebox(0,0)[cc]{$+$}}
\put(95.00,37.00){\makebox(0,0)[cc]{$+$}} \put(98.00,37.00){\line(1,0){25.00}}
\put(110.50,37.00){\oval(21.00,8.00)[b]} \put(121.00,37.00){\line(-6,5){6.00}}
\put(115.00,42.00){\vector(-1,-1){5.00}}
\put(110.00,29.00){\makebox(0,0)[cc]{$1 \longleftrightarrow 2$}}
\put(110.00,27.50){\line(0,1){3.00}} \put(11.00,14.00){\makebox(0,0)[cc]{$+$}}
\put(14.00,14.00){\line(1,0){24.00}} \put(26.00,14.00){\oval(20.00,8.00)[b]}
\put(36.00,14.00){\line(-6,5){6.00}} \put(30.00,19.00){\line(-1,-1){5.00}}
\put(22.00,4.00){\vector(2,-1){6.00}} \put(22.00,4.00){\vector(2,1){6.00}}
\put(24.00,7.00){\makebox(0,0)[cc]{1}} \put(24.00,1.00){\makebox(0,0)[cc]{2}}
\put(22.00,4.00){\vector(1,0){10.00}} \put(22.00,2.00){\line(0,1){4.00}}
\put(48.00,14.00){\line(1,0){25.00}} \put(60.50,14.00){\oval(21.00,8.00)[t]}
\put(69.00,8.00){\makebox(0,0)[cc]{$2 \uparrow$}}
\put(43.00,14.00){\makebox(0,0)[cc]{$+$}}
\put(65.00,37.00){\makebox(0,0)[cc]{$+$}}
\put(126.00,37.00){\makebox(0,0)[cc]{$+$}}
\put(78.00,14.00){\makebox(0,0)[cc]{$+$}} \put(83.00,14.00){\line(1,0){26.00}}
\put(96.00,14.00){\oval(19.00,8.00)[t]} \put(84.00,6.00){\makebox(0,0)[cc]{$1
\downarrow$}} \put(103.00,8.00){\makebox(0,0)[cc]{\small $f_1'''(1)/f_1''(1)
\sim$}} \put(103.00,4.00){\makebox(0,0)[cc]{\small $\sim \/ \mapsto/
f_1''(1)$}} \put(49.00,6.00){\vector(1,0){12.00}}
\put(49.00,6.00){\vector(3,1){9.00}} \put(58.00,3.00){\vector(-3,1){6.00}}
\put(55.00,4.00){\vector(-3,1){6.00}} \put(49.00,3.00){\line(0,1){6.00}}
\put(52.00,2.00){\makebox(0,0)[cc]{$2$}}
\put(52.00,9.00){\makebox(0,0)[cc]{$1$}}
\end{picture}
\end{equation}

\begin{equation}
\label{d2-6} \unitlength 1.00mm \special{em:linewidth 0.4pt}
\linethickness{0.4pt}
\begin{picture}(125.00,32.00)
\put(3.00,27.00){\makebox(0,0)[cc]{$\partial$}}
\put(6.00,27.00){\line(1,0){24.00}} \put(20.50,27.00){\oval(15.00,8.00)[b]}
\put(8.00,27.00){\line(3,5){3.00}} \put(20.00,32.00){\vector(1,-1){5.00}}
\put(11.00,32.00){\line(1,0){9.00}} \put(13.00,27.00){\line(3,4){3.00}}
\put(16.00,31.00){\vector(1,-2){2.00}}
\put(33.00,27.00){\makebox(0,0)[cc]{$=$}} \put(36.00,27.00){\line(1,0){24.00}}
\put(50.50,27.00){\oval(15.00,8.00)[b]} \put(38.00,27.00){\line(3,5){3.00}}
\put(41.00,32.00){\line(1,0){9.00}} \put(43.00,27.00){\line(3,4){3.00}}
\put(46.00,31.00){\vector(1,-2){2.00}} \put(50.00,32.00){\line(1,-1){5.00}}
\put(63.00,27.00){\makebox(0,0)[cc]{$+$}} \put(66.00,27.00){\line(1,0){24.00}}
\put(80.00,27.00){\oval(16.00,8.00)[b]} \put(68.00,27.00){\line(3,5){3.00}}
\put(71.00,32.00){\line(1,0){9.00}} \put(80.00,32.00){\vector(1,-1){5.00}}
\put(72.00,27.00){\line(3,4){3.00}} \put(75.00,31.00){\line(3,-4){3.00}}
\put(93.00,27.00){\makebox(0,0)[cc]{$+$}} \put(96.00,27.00){\line(1,0){24.00}}
\put(108.00,27.00){\oval(20.00,8.00)[b]} \put(98.00,27.00){\line(2,5){2.00}}
\put(100.00,32.00){\vector(1,-1){5.00}} \put(118.00,27.00){\line(-2,5){2.00}}
\put(116.00,32.00){\vector(-1,-1){5.00}}
\put(123.00,27.00){\makebox(0,0)[cc]{$+$}}
\put(2.00,11.00){\makebox(0,0)[cc]{$+$}} \put(5.00,11.00){\line(1,0){19.00}}
\put(18.50,11.00){\oval(9.00,8.00)[b]} \put(7.00,11.00){\line(6,5){6.00}}
\put(13.00,16.00){\vector(1,-1){5.00}} \put(11.00,5.00){\makebox(0,0)[cc]{$2
\downarrow$}} \put(27.00,11.00){\makebox(0,0)[cc]{$+$}}
\put(30.00,11.00){\line(1,0){19.00}} \put(44.00,11.00){\oval(8.00,8.00)[b]}
\put(32.00,11.00){\line(6,5){6.00}} \put(38.00,16.00){\vector(1,-1){5.00}}
\put(48.00,11.00){\line(0,1){5.00}} \put(48.00,16.00){\vector(-1,-1){5.00}}
\put(52.00,11.00){\makebox(0,0)[cc]{$+$}} \put(55.00,11.00){\line(1,0){19.00}}
\put(68.50,11.00){\oval(9.00,8.00)[b]} \put(57.00,11.00){\line(2,5){2.00}}
\put(59.00,16.00){\vector(1,-1){5.00}} \put(73.00,11.00){\line(-4,5){4.00}}
\put(69.00,16.00){\vector(0,-1){5.00}}
\put(77.00,11.00){\makebox(0,0)[cc]{$+$}} \put(105.00,11.00){\line(1,0){19.00}}
\put(114.50,11.00){\oval(15.00,8.00)[b]} \put(107.00,11.00){\line(3,5){3.00}}
\put(110.00,16.00){\vector(1,-1){5.00}}
\put(102.00,11.00){\makebox(0,0)[cc]{$+$}} \put(80.00,11.00){\line(1,0){20.00}}
\put(90.00,11.00){\oval(16.00,8.00)[b]} \put(98.00,11.00){\line(-3,5){3.00}}
\put(95.00,16.00){\vector(-1,-1){5.00}} \put(90.00,3.00){\makebox(0,0)[cc]{$1
\longleftrightarrow 2$}} \put(90.00,1.50){\line(0,1){3.00}}
\put(119.00,4.00){\vector(3,1){6.00}} \put(119.00,4.00){\vector(3,-1){6.00}}
\put(122.00,7.00){\makebox(0,0)[cc]{$3$}}
\put(122.00,1.00){\makebox(0,0)[cc]{$1$}}
\end{picture}
\end{equation}

\begin{equation}
\label{d2-7} \unitlength=1.00mm \special{em:linewidth 0.4pt}
\linethickness{0.4pt}
\begin{picture}(102.00,30.00)
\put(2.00,25.00){\makebox(0,0)[cc]{$\partial$}}
\put(5.00,25.00){\line(1,0){19.00}} \put(16.50,25.00){\oval(11.00,8.00)[b]}
\put(7.00,25.00){\line(1,1){5.00}} \put(12.00,30.00){\vector(1,-1){5.00}}
\put(9.00,19.00){\makebox(0,0)[cc]{$2 \mapsto$}}
\put(27.00,25.00){\makebox(0,0)[cc]{$=$}} \put(30.00,25.00){\line(1,0){19.00}}
\put(41.50,25.00){\oval(11.00,8.00)[b]} \put(32.00,25.00){\line(1,1){5.00}}
\put(37.00,30.00){\line(1,-1){5.00}} \put(34.00,19.00){\makebox(0,0)[cc]{$2
\mapsto$}} \put(52.00,25.00){\makebox(0,0)[cc]{$+$}}
\put(55.00,25.00){\line(1,0){19.00}} \put(66.50,25.00){\oval(11.00,8.00)[b]}
\put(57.00,25.00){\line(1,1){5.00}} \put(62.00,30.00){\vector(1,-1){5.00}}
\put(58.00,19.00){\makebox(0,0)[cc]{$2 \updownarrow$}}
\put(77.00,25.00){\makebox(0,0)[cc]{$+$}} \put(80.00,25.00){\line(1,0){19.00}}
\put(89.50,25.00){\oval(15.00,8.00)[b]} \put(82.00,25.00){\line(3,5){3.00}}
\put(85.00,30.00){\vector(1,-1){5.00}} \put(81.00,19.00){\makebox(0,0)[cc]{$1
\mapsto$}} \put(102.00,25.00){\makebox(0,0)[cc]{$+$}}
\put(24.00,9.00){\line(1,0){19.00}} \put(33.00,9.00){\oval(14.00,8.00)[t]}
\put(27.00,4.00){\makebox(0,0)[cc]{$1 \mapsto$}}
\put(37.00,4.00){\vector(3,1){6.00}} \put(37.00,4.00){\vector(3,-1){6.00}}
\put(40.00,7.00){\makebox(0,0)[cc]{$2$}}
\put(40.00,1.00){\makebox(0,0)[cc]{$1$}}
\put(21.00,9.00){\makebox(0,0)[cc]{$+$}}
\put(53.00,9.00){\makebox(0,0)[cc]{$+$}} \put(63.00,9.00){\line(1,0){24.00}}
\put(75.00,9.00){\oval(18.00,8.00)[t]} \put(63.00,4.00){\makebox(0,0)[cc]{$1
\downarrow$}} \put(81.00,6.00){\makebox(0,0)[cc]{\small $f_1'''(1)/f_1''(1)
\sim$}} \put(80.00,2.00){\makebox(0,0)[cc]{\small $\sim \/ \mapsto/ f_1''(1)$}}
\end{picture}
\end{equation}

\begin{equation}
\label{d2-8} \unitlength=1.00mm \special{em:linewidth 0.4pt}
\linethickness{0.4pt}
\begin{picture}(124.00,11.00)
\put(2.00,6.00){\makebox(0,0)[cc]{$\partial$}}
\put(5.00,6.00){\line(1,0){19.00}} \put(19.00,6.00){\oval(8.00,10.00)[b]}
\put(7.00,6.00){\line(6,5){6.00}} \put(13.00,11.00){\vector(1,-1){5.00}}
\put(15.00,6.00){\line(-4,-5){4.00}} \put(11.00,1.00){\vector(0,1){5.00}}
\put(27.00,6.00){\makebox(0,0)[cc]{$=$}} \put(30.00,6.00){\line(1,0){19.00}}
\put(44.00,6.00){\oval(8.00,10.00)[b]} \put(40.00,6.00){\line(-4,-5){4.00}}
\put(36.00,1.00){\vector(0,1){5.00}} \put(32.00,6.00){\line(6,5){6.00}}
\put(38.00,11.00){\line(6,-5){6.00}} \put(52.00,6.00){\makebox(0,0)[cc]{$+$}}
\put(55.00,6.00){\line(1,0){19.00}} \put(69.00,6.00){\oval(8.00,10.00)[b]}
\put(65.00,6.00){\line(-2,-5){2.00}} \put(63.00,1.00){\line(-2,5){2.00}}
\put(57.00,6.00){\line(6,5){6.00}} \put(63.00,11.00){\vector(1,-1){5.00}}
\put(77.00,6.00){\makebox(0,0)[cc]{$+$}} \put(80.00,6.00){\line(1,0){19.00}}
\put(93.50,6.00){\oval(9.00,10.00)[b]} \put(89.00,6.00){\line(-2,-5){2.00}}
\put(87.00,1.00){\vector(-1,1){5.00}} \put(82.00,6.00){\line(6,5){6.00}}
\put(88.00,11.00){\vector(1,-1){5.00}}
\put(102.00,6.00){\makebox(0,0)[cc]{$+$}} \put(105.00,6.00){\line(1,0){19.00}}
\put(118.50,6.00){\oval(9.00,10.00)[b]} \put(107.00,6.00){\line(6,5){6.00}}
\put(113.00,11.00){\vector(1,-1){5.00}} \put(109.00,1.00){\makebox(0,0)[cc]{$2
\uparrow$}}
\end{picture}
\end{equation}

\begin{equation}
\label{d2-9} \unitlength 1.00mm \linethickness{0.4pt}
\begin{picture}(125.00,16.00)
\put(2.00,11.00){\makebox(0,0)[cc]{$\partial$}}
\put(6.00,11.00){\line(1,0){24.00}} \put(17.00,11.00){\oval(20.00,10.00)[t]}
\put(33.00,11.00){\makebox(0,0)[cc]{$=$}} \put(36.00,11.00){\line(1,0){24.00}}
\put(48.00,11.00){\oval(20.00,10.00)[t]} \put(38.00,11.00){\line(3,-2){5.00}}
\put(43.00,7.67){\line(3,2){5.00}} \put(64.00,11.00){\makebox(0,0)[cc]{$+$}}
\put(68.00,11.00){\line(1,0){24.00}} \put(80.00,11.00){\oval(16.00,10.00)[t]}
\put(69.00,7.50){\makebox(0,0)[cc]{$1 \downarrow$}}
\put(99.00,11.00){\line(1,0){24.00}} \put(112.00,11.00){\oval(16.00,10.00)[t]}
\put(122.00,7.50){\makebox(0,0)[cc]{$2 \uparrow$}}
\put(96.00,11.00){\makebox(0,0)[cc]{$+$}}
\put(67.00,0.00){\framebox(27.00,5.00)[cc]{\small
$f_1''(1),f_1'''(1),\mapsto$}}
\put(98.00,0.00){\framebox(27.00,5.00)[cc]{\small
$f_1''(1),f_1'''(1),\mapsto$}}
\put(35.50,0.00){\framebox(25.00,5.00)[cc]{\small $f_1'(1),f_1'(2),\mapsto$}}
\put(5.50,0.00){\framebox(25.00,5.00)[cc]{\small $f_1'(1),f_1'(2),\mapsto$}}
\put(7.00,11.00){\line(3,-2){6.00}} \put(13.00,7.00){\vector(3,2){6.00}}
\end{picture}
\end{equation}

\begin{equation}
\label{d2-91} \unitlength 1.00mm \linethickness{0.4pt}
\begin{picture}(126.00,16.00)
\put(3.00,11.00){\makebox(0,0)[cc]{$\partial$}}
\put(6.00,11.00){\line(1,0){24.00}} \put(13.50,11.00){\oval(9.00,10.00)[b]}
\put(28.00,11.00){\line(-1,1){5.00}} \put(23.00,16.00){\vector(-1,-1){5.00}}
\put(34.00,11.00){\makebox(0,0)[cc]{$=$}} \put(69.00,11.00){\line(1,0){24.00}}
\put(91.00,11.00){\line(-2,1){10.00}} \put(81.00,16.00){\vector(-2,-1){10.00}}
\put(71.00,11.00){\makebox(0,0)[cc]{$*$}}
\put(66.00,11.00){\makebox(0,0)[cc]{$+$}} \put(37.00,11.00){\line(1,0){24.00}}
\put(44.50,11.00){\oval(9.00,10.00)[b]} \put(59.00,11.00){\line(-1,1){5.00}}
\put(98.00,11.00){\makebox(0,0)[cc]{$+$}} \put(101.00,11.00){\line(1,0){24.00}}
\put(113.00,11.00){\oval(18.00,10.00)[t]}
\put(119.00,7.50){\makebox(0,0)[cc]{$2 \uparrow$}}
\put(54.00,16.00){\line(-1,-1){5.00}}
\put(5.00,0.00){\framebox(26.00,5.00)[cc]{\small $f_1'(1),f_1'(2),\mapsto$}}
\put(67.00,0.00){\framebox(28.00,5.00)[cc]{\small
$f_1''(1),f_1'''(1),\mapsto$}}
\put(36.00,0.00){\framebox(26.00,5.00)[cc]{\small $f_1'(1),f_1'(2),\mapsto$}}
\put(100.00,0.00){\framebox(26.00,5.00)[cc]{\small $f_1'(1),f_1''(2),\mapsto$}}
\end{picture}
\end{equation}

Summing up the right-hand parts of these equations, we get the following
statement.

\begin{proposition}
The cycle $d^2(TT) \subset F_1$ presented by the linear combination
$($\ref{diftwo}$)$ is equal to the boundary $($mod 2$)$ of the sum of ten
varieties shown in the left parts of equations
$($\ref{d2-1}$)$--$($\ref{d2-91}$)$.
\end{proposition}

In this summation we use the following relations. Let us denote by (a;b) the
$b$-th summand in the right-hand part of the equation (a). Then (\ref{d2-1};6)
+ (\ref{d2-2};5) = (\ref{d2-3};2); \ (\ref{d2-1};7) + (\ref{d2-2};7) =
(\ref{d2-4};4); \ (\ref{d2-1};5) + (\ref{d2-2};6) = (\ref{d2-6};3); \
(\ref{d2-8};3) + (\ref{d2-6};6) = (\ref{d2-6};5); \ (\ref{d2-8};4) +
(\ref{d2-6};4) = (\ref{d2-7};2).

\subsection{The third differential}

Ten varieties described by left parts of equations (\ref{d2-1})--(\ref{d2-91})
form a chain in the term $F_1$ of the resolved discriminant, i.e. in the {\em
tautological resolution} of this discriminant, see \S \ref{method}.

Finally, we consider the image of this chain in the discriminant itself. The
image of any of ten components of this image is a subvariety in the space
${\mathcal K}_n$ of maps $\R^1 \to \R^n$, distinguished by conditions, whose
notation is obtained from the notation of the corresponding variety in $F_1$ by
replacing its unique chord by a zigzag with the same endpoints. It remains to
span the sum of these varieties by a chain in the space ${\mathcal K}_n$.
Proceeding as before, we find five varieties indicated in the left parts of the
following identities (\ref{d3-1})--(\ref{d3-5}).

\begin{equation}
\label{d3-1} \unitlength 1.00mm \linethickness{0.4pt}
\begin{picture}(123.00,46.00)
\put(3.00,40.00){\makebox(0,0)[cc]{$\partial$}}
\put(6.00,40.00){\line(1,0){34.00}} \put(39.00,40.00){\line(-4,3){8.00}}
\put(31.00,46.00){\vector(-4,-3){8.00}} \put(23.00,40.00){\line(-2,3){4.00}}
\put(19.00,46.00){\vector(-2,-3){4.00}}
\put(43.00,40.00){\makebox(0,0)[cc]{$=$}} \put(46.00,40.00){\line(1,0){34.00}}
\put(79.00,40.00){\line(-4,3){8.00}} \put(71.00,46.00){\vector(-4,-3){8.00}}
\put(63.00,40.00){\line(-2,3){4.00}} \put(59.00,46.00){\vector(-2,-3){4.00}}
\put(86.00,40.00){\line(1,0){34.00}} \put(119.00,40.00){\line(-4,3){8.00}}
\put(111.00,46.00){\vector(-4,-3){8.00}} \put(87.00,40.00){\line(2,-1){12.00}}
\put(99.00,34.00){\vector(2,1){12.00}} \put(6.00,22.00){\line(1,0){34.00}}
\put(23.00,22.00){\line(-2,3){4.00}} \put(19.00,28.00){\vector(-2,-3){4.00}}
\put(83.00,40.00){\makebox(0,0)[cc]{$+$}}
\put(123.00,40.00){\makebox(0,0)[cc]{$+$}}
\put(3.00,22.00){\makebox(0,0)[cc]{$+$}}
\put(43.00,22.00){\makebox(0,0)[cc]{$+$}} \put(46.00,22.00){\line(1,0){34.00}}
\put(58.00,22.00){\line(-2,3){4.00}} \put(54.00,28.00){\vector(-1,-1){6.00}}
\put(83.00,22.00){\makebox(0,0)[cc]{$+$}} \put(86.00,22.00){\line(1,0){34.00}}
\put(110.00,16.00){\makebox(0,0)[cc]{$2 \uparrow$}}
\put(19.00,6.00){\makebox(0,0)[cc]{$+$}} \put(56.00,6.00){\line(-1,0){34.00}}
\put(44.00,6.00){\line(-2,3){4.00}} \put(40.00,12.00){\vector(-1,-1){6.00}}
\put(54.00,6.00){\line(-2,3){4.00}} \put(50.00,12.00){\vector(-1,-1){6.00}}
\put(59.00,6.00){\makebox(0,0)[cc]{$+$}} \put(62.00,6.00){\line(1,0){34.00}}
\put(64.00,6.00){\line(2,-1){12.00}} \put(76.00,0.00){\line(1,0){6.00}}
\put(82.00,0.00){\vector(2,1){12.00}} \put(94.00,6.00){\line(-2,3){4.00}}
\put(90.00,12.00){\vector(-1,-1){6.00}} \put(84.00,6.00){\line(-2,3){4.00}}
\put(80.00,12.00){\vector(-1,-1){6.00}}
\put(123.00,22.00){\makebox(0,0)[cc]{$+$}}
\put(119.00,22.00){\line(-1,1){6.00}} \put(113.00,28.00){\line(-1,0){9.00}}
\put(104.00,28.00){\vector(-1,-1){6.00}} \put(24.00,6.00){\line(1,-1){6.00}}
\put(30.00,0.00){\line(1,0){8.00}} \put(38.00,0.00){\vector(1,1){6.00}}
\put(78.00,22.00){\line(-1,1){6.00}} \put(72.00,28.00){\line(-1,0){8.00}}
\put(64.00,28.00){\vector(-1,-1){6.00}} \put(48.00,22.00){\line(1,-1){6.00}}
\put(54.00,16.00){\line(1,0){8.00}} \put(62.00,16.00){\vector(1,1){6.00}}
\put(7.00,40.00){\line(1,-1){6.00}} \put(13.00,34.00){\line(1,0){12.00}}
\put(25.00,34.00){\vector(1,1){6.00}} \put(39.00,22.00){\line(-4,3){8.00}}
\put(31.00,28.00){\line(-4,-3){8.00}} \put(103.00,40.00){\line(-2,3){4.00}}
\put(99.00,46.00){\line(-2,-3){4.00}} \put(47.00,40.00){\line(1,-1){6.00}}
\put(53.00,34.00){\line(1,0){12.00}} \put(65.00,34.00){\line(1,1){6.00}}
\put(7.00,22.00){\line(1,-1){6.00}} \put(13.00,16.00){\line(1,0){12.00}}
\put(25.00,16.00){\vector(1,1){6.00}} \put(88.00,22.00){\line(1,-1){6.00}}
\put(94.00,16.00){\line(1,0){8.00}} \put(102.00,16.00){\vector(1,1){6.00}}
\end{picture}
\end{equation}

\begin{equation}
\label{d3-2} \unitlength=1.00mm \special{em:linewidth 0.4pt}
\linethickness{0.4pt}
\begin{picture}(123.00,51.00)
\put(3.00,45.00){\makebox(0,0)[cc]{$\partial$}}
\put(6.00,45.00){\line(1,0){34.00}} \put(46.00,45.00){\line(1,0){34.00}}
\put(83.00,45.00){\makebox(0,0)[cc]{$+$}} \put(86.00,45.00){\line(1,0){34.00}}
\put(6.00,7.00){\line(1,0){34.00}} \put(46.00,7.00){\line(1,0){34.00}}
\put(83.00,7.00){\makebox(0,0)[cc]{$+$}} \put(86.00,7.00){\line(1,0){34.00}}
\put(6.00,26.00){\line(1,0){34.00}} \put(46.00,26.00){\line(1,0){34.00}}
\put(83.00,26.00){\makebox(0,0)[cc]{$+$}} \put(86.00,26.00){\line(1,0){34.00}}
\put(7.00,45.00){\line(2,-1){12.00}} \put(19.00,39.00){\vector(2,1){12.00}}
\put(39.00,45.00){\line(-4,3){8.00}} \put(31.00,51.00){\line(-1,0){8.00}}
\put(23.00,51.00){\vector(-4,-3){8.00}} \put(15.00,45.00){\line(1,-1){4.00}}
\put(19.00,41.00){\vector(1,1){4.00}} \put(43.00,45.00){\makebox(0,0)[cc]{$=$}}
\put(79.00,45.00){\line(-4,3){8.00}} \put(71.00,51.00){\line(-1,0){8.00}}
\put(63.00,51.00){\vector(-4,-3){8.00}} \put(55.00,45.00){\line(1,-1){4.00}}
\put(59.00,41.00){\vector(1,1){4.00}} \put(86.00,45.00){\line(1,0){34.00}}
\put(87.00,45.00){\line(2,-1){12.00}} \put(99.00,39.00){\vector(2,1){12.00}}
\put(119.00,45.00){\line(-4,3){8.00}} \put(111.00,51.00){\line(-1,0){8.00}}
\put(103.00,51.00){\vector(-4,-3){8.00}}
\put(123.00,45.00){\makebox(0,0)[cc]{$+$}}
\put(123.00,26.00){\makebox(0,0)[cc]{$+$}} \put(88.00,26.00){\line(4,-3){8.00}}
\put(96.00,20.00){\vector(2,1){12.00}} \put(118.00,26.00){\line(-4,3){8.00}}
\put(110.00,32.00){\vector(-2,-1){12.00}}
\put(109.00,19.00){\makebox(0,0)[cc]{$2 \downarrow$}}
\put(78.00,26.00){\line(-1,1){6.00}} \put(72.00,32.00){\line(-1,0){18.00}}
\put(54.00,32.00){\vector(-1,-1){6.00}} \put(48.00,26.00){\line(2,-3){4.00}}
\put(52.00,20.00){\vector(1,1){6.00}} \put(48.00,26.00){\line(3,1){12.00}}
\put(60.00,30.00){\vector(2,-1){8.00}}
\put(43.00,26.00){\makebox(0,0)[cc]{$+$}}
\put(3.00,26.00){\makebox(0,0)[cc]{$+$}} \put(7.00,26.00){\line(2,-1){12.00}}
\put(19.00,20.00){\vector(2,1){12.00}} \put(15.00,26.00){\line(1,-1){4.00}}
\put(19.00,22.00){\vector(1,1){4.00}} \put(3.00,7.00){\makebox(0,0)[cc]{$+$}}
\put(8.00,7.00){\line(1,-1){6.00}} \put(14.00,1.00){\line(1,0){8.00}}
\put(22.00,1.00){\vector(1,1){6.00}} \put(38.00,7.00){\line(-1,1){6.00}}
\put(32.00,13.00){\line(-1,0){8.00}} \put(24.00,13.00){\vector(-1,-1){6.00}}
\put(18.00,7.00){\line(2,1){8.00}} \put(26.00,11.00){\vector(1,-2){2.00}}
\put(43.00,7.00){\makebox(0,0)[cc]{$+$}} \put(48.00,7.00){\line(1,-1){6.00}}
\put(54.00,1.00){\line(1,0){18.00}} \put(72.00,1.00){\vector(1,1){6.00}}
\put(78.00,7.00){\line(-1,1){6.00}} \put(72.00,13.00){\line(-1,0){8.00}}
\put(64.00,13.00){\vector(-1,-1){6.00}} \put(58.00,7.00){\line(3,-2){6.00}}
\put(64.00,3.00){\vector(1,1){4.00}} \put(118.00,7.00){\line(-1,1){6.00}}
\put(112.00,13.00){\line(-1,0){18.00}} \put(94.00,13.00){\vector(-1,-1){6.00}}
\put(88.00,7.00){\line(3,-2){9.00}} \put(97.00,1.00){\vector(1,1){6.00}}
\put(109.00,2.00){\makebox(0,0)[cc]{$1 \longleftrightarrow 2$}}
\put(109.00,0.50){\line(0,1){3.00}} \put(15.00,26.00){\line(1,1){6.00}}
\put(21.00,32.00){\line(1,0){12.00}} \put(33.00,32.00){\line(1,-1){6.00}}
\put(47.00,45.00){\line(1,-1){6.00}} \put(53.00,39.00){\line(1,0){12.00}}
\put(65.00,39.00){\line(1,1){6.00}} \put(95.00,45.00){\line(1,-1){4.00}}
\put(99.00,41.00){\line(1,1){4.00}}
\end{picture}
\end{equation}

\begin{equation}
\label{d3-3} \unitlength=1.00mm \special{em:linewidth 0.4pt}
\linethickness{0.4pt}
\begin{picture}(123.00,40.00)
\put(3.00,11.00){\makebox(0,0)[cc]{$+$}} \put(6.00,11.00){\line(1,0){22.00}}
\put(31.00,11.00){\makebox(0,0)[cc]{$+$}} \put(34.00,11.00){\line(1,0){22.00}}
\put(59.00,11.00){\makebox(0,0)[cc]{$+$}} \put(62.00,11.00){\line(1,0){22.00}}
\put(64.00,11.00){\line(3,2){9.00}} \put(73.00,17.00){\vector(3,-2){9.00}}
\put(82.00,11.00){\line(-3,-4){4.50}} \put(77.50,5.00){\vector(-3,4){4.50}}
\put(68.00,6.00){\makebox(0,0)[cc]{$2 \mapsto$}}
\put(36.00,11.00){\line(3,4){4.50}} \put(40.50,17.00){\vector(3,-4){4.50}}
\put(54.00,11.00){\line(-3,4){4.50}} \put(49.50,17.00){\vector(-3,-4){4.50}}
\put(42.00,6.00){\makebox(0,0)[cc]{$2 \mapsto$}}
\put(8.00,11.00){\line(3,-4){4.50}} \put(12.50,5.00){\vector(3,4){4.50}}
\put(26.00,11.00){\line(-3,2){9.00}} \put(17.00,17.00){\vector(-3,-2){9.00}}
\put(21.00,6.00){\makebox(0,0)[cc]{$1 \mapsto$}}
\put(3.00,35.00){\makebox(0,0)[cc]{$\partial$}}
\put(6.00,35.00){\line(1,0){24.00}} \put(28.00,35.00){\line(-3,2){7.50}}
\put(20.50,40.00){\vector(-3,-2){7.50}} \put(8.00,35.00){\line(3,-2){7.50}}
\put(15.50,30.00){\vector(3,2){7.50}} \put(25.00,30.00){\makebox(0,0)[cc]{$2
\mapsto$}} \put(36.00,35.00){\line(1,0){24.00}}
\put(58.00,35.00){\line(-3,2){7.50}} \put(50.50,40.00){\vector(-3,-2){7.50}}
\put(38.00,35.00){\line(3,-2){7.50}} \put(45.50,30.00){\line(3,2){7.50}}
\put(55.00,30.00){\makebox(0,0)[cc]{$2 \mapsto$}}
\put(66.00,35.00){\line(1,0){24.00}} \put(88.00,35.00){\line(-3,2){7.50}}
\put(80.50,40.00){\line(-3,-2){7.50}} \put(68.00,35.00){\line(3,-2){7.50}}
\put(75.50,30.00){\vector(3,2){7.50}} \put(85.00,30.00){\makebox(0,0)[cc]{$2
\mapsto$}} \put(96.00,35.00){\line(1,0){24.00}}
\put(118.00,35.00){\line(-3,2){7.50}} \put(110.50,40.00){\vector(-3,-2){7.50}}
\put(98.00,35.00){\line(3,-2){7.50}} \put(105.50,30.00){\vector(3,2){7.50}}
\put(33.00,35.00){\makebox(0,0)[cc]{$=$}}
\put(63.00,35.00){\makebox(0,0)[cc]{$+$}}
\put(93.00,35.00){\makebox(0,0)[cc]{$+$}}
\put(123.00,35.00){\makebox(0,0)[cc]{$+$}}
\put(115.00,30.00){\makebox(0,0)[cc]{$2 \updownarrow$}}
\put(87.00,11.00){\makebox(0,0)[cc]{$+$}} \put(90.00,11.00){\line(1,0){29.00}}
\put(117.00,11.00){\line(-2,1){10.00}} \put(107.00,16.00){\line(-1,0){5.00}}
\put(102.00,16.00){\vector(-2,-1){10.00}} \put(92.00,5.00){\makebox(0,0)[cc]{$1
\downarrow$}} \put(111.00,7.00){\makebox(0,0)[cc]{\small $-f_1'''(1)/f_1''(1)
\sim$}} \put(111.00,3.00){\makebox(0,0)[cc]{\small $\sim \mapsto / f_1''(1)$}}
\end{picture}
\end{equation}

\begin{equation}
\label{d3-4} \unitlength=1.00mm \special{em:linewidth 0.4pt}
\linethickness{0.4pt}
\begin{picture}(122.00,35.00)
\put(2.00,30.00){\makebox(0,0)[cc]{$\partial$}}
\put(19.00,23.00){\vector(2,1){6.00}} \put(19.00,23.00){\vector(2,-1){6.00}}
\put(19.00,23.00){\vector(1,0){10.00}} \put(19.00,21.00){\line(0,1){4.00}}
\put(21.00,19.00){\makebox(0,0)[cc]{2}} \put(21.00,27.00){\makebox(0,0)[cc]{1}}
\put(5.00,30.00){\line(1,0){24.00}} \put(27.00,30.00){\line(-2,1){10.00}}
\put(17.00,35.00){\vector(-2,-1){10.00}} \put(7.00,30.00){\line(1,-1){5.00}}
\put(12.00,25.00){\vector(1,1){5.00}} \put(49.00,23.00){\vector(2,1){6.00}}
\put(49.00,23.00){\vector(2,-1){6.00}} \put(49.00,23.00){\vector(1,0){10.00}}
\put(49.00,21.00){\line(0,1){4.00}} \put(51.00,19.00){\makebox(0,0)[cc]{2}}
\put(51.00,27.00){\makebox(0,0)[cc]{1}} \put(35.00,30.00){\line(1,0){24.00}}
\put(57.00,30.00){\line(-2,1){10.00}} \put(37.00,30.00){\line(1,-1){5.00}}
\put(42.00,25.00){\vector(1,1){5.00}} \put(79.00,23.00){\vector(2,1){6.00}}
\put(79.00,23.00){\vector(2,-1){6.00}} \put(79.00,23.00){\vector(1,0){10.00}}
\put(79.00,21.00){\line(0,1){4.00}} \put(81.00,19.00){\makebox(0,0)[cc]{2}}
\put(81.00,27.00){\makebox(0,0)[cc]{1}} \put(65.00,30.00){\line(1,0){24.00}}
\put(87.00,30.00){\line(-2,1){10.00}} \put(77.00,35.00){\vector(-2,-1){10.00}}
\put(67.00,30.00){\line(1,-1){5.00}} \put(95.00,30.00){\line(1,0){24.00}}
\put(117.00,30.00){\line(-2,1){10.00}}
\put(107.00,35.00){\vector(-2,-1){10.00}} \put(97.00,30.00){\line(1,-1){5.00}}
\put(102.00,25.00){\vector(1,1){5.00}}
\put(32.00,30.00){\makebox(0,0)[cc]{$=$}}
\put(62.00,30.00){\makebox(0,0)[cc]{$+$}}
\put(92.00,30.00){\makebox(0,0)[cc]{$+$}}
\put(122.00,30.00){\makebox(0,0)[cc]{$+$}}
\put(111.00,24.00){\makebox(0,0)[cc]{$1 \mapsto$}}
\put(72.00,25.00){\line(1,1){5.00}} \put(47.00,35.00){\line(-2,-1){10.00}}
\put(9.00,9.00){\line(1,0){24.00}} \put(31.00,9.00){\line(-2,1){10.00}}
\put(21.00,14.00){\vector(-2,-1){10.00}} \put(11.00,9.00){\line(1,-1){5.00}}
\put(16.00,4.00){\vector(1,1){5.00}} \put(6.00,9.00){\makebox(0,0)[cc]{$+$}}
\put(36.00,9.00){\makebox(0,0)[cc]{$+$}} \put(25.00,3.00){\makebox(0,0)[cc]{$2
\mapsto$}} \put(39.00,9.00){\line(1,0){24.00}}
\put(61.00,9.00){\line(-2,1){10.00}} \put(51.00,14.00){\vector(-2,-1){10.00}}
\put(41.00,9.00){\line(1,-1){5.00}} \put(46.00,4.00){\vector(1,1){5.00}}
\put(69.00,9.00){\makebox(0,0)[cc]{$+$}} \put(58.00,5.00){\vector(-1,0){5.00}}
\put(58.00,5.00){\vector(1,0){5.00}} \put(58.00,3.00){\line(0,1){4.00}}
\put(65.00,5.00){\makebox(0,0)[cc]{$2$}}
\put(51.00,5.00){\makebox(0,0)[cc]{$1$}} \put(75.00,9.00){\line(1,0){29.00}}
\put(102.00,9.00){\line(-2,1){10.00}} \put(92.00,14.00){\line(-1,0){5.00}}
\put(87.00,14.00){\vector(-2,-1){10.00}} \put(77.00,3.00){\makebox(0,0)[cc]{$1
\downarrow$}} \put(95.00,5.00){\makebox(0,0)[cc]{\small $f_1'''(1)/f_1''(1)
\sim$}} \put(95.00,1.00){\makebox(0,0)[cc]{\small $\sim \mapsto / f_1''(1)$}}
\end{picture}
\end{equation}

\begin{equation}
\label{d3-5} \unitlength 1.00mm \special{em:linewidth 0.4pt}
\linethickness{0.4pt}
\begin{picture}(125.00,17.00)
\put(3.00,12.00){\makebox(0,0)[cc]{$\partial$}}
\put(6.00,12.00){\line(1,0){24.00}} \put(28.00,12.00){\line(-2,1){10.00}}
\put(18.00,17.00){\vector(-2,-1){10.00}} \put(8.00,12.00){\line(1,-1){5.00}}
\put(13.00,7.00){\vector(1,1){5.00}} \put(36.00,12.00){\line(1,0){24.00}}
\put(58.00,12.00){\line(-2,1){10.00}} \put(38.00,12.00){\line(1,-1){5.00}}
\put(43.00,7.00){\vector(1,1){5.00}} \put(66.00,12.00){\line(1,0){24.00}}
\put(88.00,12.00){\line(-2,1){10.00}} \put(78.00,17.00){\vector(-2,-1){10.00}}
\put(68.00,12.00){\line(1,-1){5.00}} \put(96.00,12.00){\line(1,0){29.00}}
\put(123.00,12.00){\line(-2,1){10.00}} \put(113.00,17.00){\line(-1,0){5.00}}
\put(108.00,17.00){\vector(-2,-1){10.00}}
\put(93.00,12.00){\makebox(0,0)[cc]{$+$}}
\put(63.00,12.00){\makebox(0,0)[cc]{$+$}}
\put(48.00,17.00){\line(-2,-1){10.00}} \put(73.00,7.00){\line(1,1){5.00}}
\put(99.00,8.00){\makebox(0,0)[cc]{$1 \downarrow$}}
\put(5.00,0.00){\framebox(26.00,5.00)[cc]{\small $f_1'(1),f_1'(2),\mapsto$}}
\put(35.00,0.00){\framebox(26.00,5.00)[cc]{\small $f_1'(1),f_1'(2),\mapsto$}}
\put(65.00,0.00){\framebox(26.00,5.00)[cc]{\small $f_1'(1),f_1'(2),\mapsto$}}
\put(96.50,0.00){\framebox(28.00,5.00)[cc]{\small
$f_1''(1),f_1'''(2),\mapsto$}} \put(34.00,12.00){\makebox(0,0)[cc]{$=$}}
\end{picture}
\end{equation}

It is easy to check that the homological sum of the right-hand parts of these
identities is equal to the sum of our ten discriminant varieties obtained from
left parts of equalities (\ref{d2-1})--(\ref{d2-91}). (We use the following
relations: (\ref{d3-1};5) + (\ref{d3-2};5) = (\ref{d3-3};3); \ (\ref{d3-3};6) +
(\ref{d3-4};4) = (\ref{d3-3};5); \ (\ref{d3-3};4) + (\ref{d3-4};6) =
(\ref{d3-5},3); \ (\ref{d3-1};4) + (\ref{d3-1};6) + (\ref{d3-1};7) +
(\ref{d3-2};4) + (\ref{d3-2};6) + (\ref{d3-2};7) = 0.)

Therefore the sum of the five varieties indicated in left parts of equalities
(\ref{d3-1})--(\ref{d3-5}) is the desired relative cycle in ${\mathcal K}_n \
(\mbox{mod } \Sigma)$.

The sum of the first and the second of these five varieties (respectively, the
third variety, respectively, the difference of the fifth and the fourth
varieties) is exactly the variety indicated in item a) (respectively, b),
respectively, c)) of Theorem 1, which is thus completely proved.

\subsection{Problems}

{\bf 1. Algorithmization.} How to do algorithmically all the same for any other
finite-type cohomology class of the space of knots?

Namely, let us consider any homology class $\gamma$ of the discriminant of the
space of knots, having some finite filtration ("order") $p$ and presented by
its "principal part", i.e. by the corresponding homology class in the term $F_p
\setminus F_{p-1}$ of the resolved discriminant. This class always is described
by some linear combination of pictures (generalized chord diagrams) as in
(\ref{prinpart}), (\ref{prinpar}), see \S \ \ref{method}. To get the
combinatorial description of a cohomology class with this principal part, we
need to calculate all the steps of the spectral sequence starting from this
part. To do it, on any step we need to find the chains spanning the consecutive
boundaries $d^r(\gamma) \subset F_{p-r} \setminus F_{p-r-1}$. Above we have
used some obvious rules: if a piece of our cycle $d^r(\gamma)$ is described by
a picture like in (\ref{difone})--(\ref{d3-4}), then it is natural to kill it
by a piece of the spanning chain, described by almost the same picture, only
replacing some one zigzag without arrows by the same zigzag with arrow at one
its end, or replacing some condition of type \unitlength 1.00mm
\special{em:linewidth 0.4pt} \linethickness{0.4pt}
\begin{picture}(12.00,1.00)
\put(6.00,-1.00){\vector(1,0){6.00}} \put(6.00,-1.00){\vector(-1,0){6.00}}
\put(6.00,-1.00){\circle*{1.00}} \put(1.00,1.00){\makebox(0,0)[cc]{\small 1}}
\put(11.00,1.00){\makebox(0,0)[cc]{\small 2}}
\end{picture}
 \
by the condition of type \unitlength 1.00mm \special{em:linewidth 0.4pt}
\linethickness{0.4pt}
\begin{picture}(8.00,5.00)
\put(4.00,0.00){\vector(1,2){2.50}} \put(4.00,0.00){\vector(-1,2){2.50}}
\put(0.33,2.00){\makebox(0,0)[cc]{\small 1}}
\put(7.67,2.00){\makebox(0,0)[cc]{\small 2}}
\end{picture}
 \ ,
or replacing some condition \unitlength 1.00mm \special{em:linewidth 0.4pt}
\linethickness{0.4pt}
\begin{picture}(11.50,3.00)
\put(6.00,1.00){\makebox(0,0)[cc]{$1 \longleftrightarrow 2$}}
\put(6.00,-0.50){\line(0,1){3.00}}
\end{picture}
 \
by the condition \unitlength 1.00mm \special{em:linewidth 0.4pt}
\linethickness{0.4pt}
\begin{picture}(6.00,5.00)
\put(0.00,1.00){\vector(2,-1){4.00}} \put(0.00,1.00){\vector(2,1){4.00}}
\put(0.00,1.00){\vector(1,0){6.00}} \put(2.00,-2.00){\makebox(0,0)[cc]{$2$}}
\put(2.00,4.00){\makebox(0,0)[cc]{$1$}} \put(0.00,-0.50){\line(0,1){3.00}}
\end{picture}
 \ , etc.
But how to decide, which of these fragments (and for which piece of the cycle)
to replace first? At which endpoint to put the arrow? Is it possible to do it
always in such a way that all the other components of the boundary of this
spanning variety would be in some sense "of lesser complexity" than the killed
one, so that our algorithm converges inductively? Which other subvarieties in
the cells of $F_p \setminus F_{p-1}$ can occur in the process of performing
this algorithm ? What are the formal rules for calculating their boundaries ?

I presume that the main filtering degree should be the number of points in
$\R^1$ participating in the definition of the subvariety, and the orientation
of arrows is not important: say, the algorithm will work if we orient all of
them from the right to the left (although, of course, other choice can provide
somewhat easier formulas).

\bigskip

{\bf 2. Orientable case.}

To do all the same for homology with integer coefficients, i.e. taking into
account orientations of our varieties. In this problem, the answers for odd and
even $n$ will be different: already the chain (\ref{prinpart}) is a $\Z$-cycle
in $F_3 \sm F_2$ only for even $n$. If $n$ is even, is it correct, that all the
calculations of \S~\ref{proof1} remain valid after imposing appropriate signs
before the pictures ?

\subsection{Proof of Proposition \protect{\ref{realiz}}}
\label{proreal}

First we specify a loop in  the space of long knots as in this Proposition. We
can assume that the standard embedding $\R^1 \to \R^3$ (with which all long
knots should coincide close to the infinity) lies in the plane $\R^2$ and has
angle $\pi/4$ with the chosen direction "to the right". Let us consider the
standard long trefoil as shown in Fig.~\ref{cat}.

\begin{figure}
\begin{center}
\unitlength 1.00mm \linethickness{0.4pt}
\begin{picture}(69.00,77.00)
\bezier{132}(29.00,45.00)(28.00,28.00)(19.33,16.50)
\bezier{96}(29.00,45.00)(29.00,60.00)(38.00,59.00)
\bezier{32}(38.00,59.00)(42.00,58.00)(43.25,56.70)
\bezier{212}(46.67,54.33)(69.00,34.00)(56.00,16.00)
\bezier{212}(14.00,18.00)(39.00,1.00)(56.00,16.00)
\bezier{164}(27.00,39.00)(3.00,28.00)(14.00,18.00)
\bezier{156}(30.00,40.00)(42.00,46.00)(53.00,69.00)
\put(53.00,69.00){\vector(1,2){4.00}} \put(41.00,58.20){\circle{1.00}}
\put(10.70,25.00){\circle{1.00}} \put(46.94,63.04){\line(-2,-3){30.00}}
\put(47.06,62.96){\line(-2,-3){30.00}} \put(50.94,60.04){\line(-2,-3){30.00}}
\put(51.04,59.94){\line(-2,-3){30.00}} \put(18.94,12.04){\line(-2,-3){3.00}}
\put(19.06,11.94){\line(-2,-3){3.00}} \put(14.94,15.04){\line(-2,-3){3.00}}
\put(15.04,14.94){\line(-2,-3){3.00}}
\put(7.90,2.70){\framebox(1.00,1.00)[cc]{}}
\bezier{36}(47.06,62.96)(51.06,64.29)(51.06,59.96)
\bezier{36}(46.94,63.04)(50.94,64.37)(50.94,60.04)
\bezier{76}(17.00,14.00)(11.50,6.33)(5.50,0.00)
\bezier{20}(12.50,6.50)(14.00,6.00)(16.00,7.50)
\bezier{16}(16.00,7.50)(16.67,9.33)(15.30,10.90)
\bezier{20}(15.30,10.90)(13.50,11.50)(12.00,10.50)
\bezier{12}(12.00,10.50)(10.67,9.00)(11.00,8.00)
\end{picture}
\caption{} \label{cat}
\end{center}
\end{figure}

Namely, we assume that close to all crossing points the projections of tangent
directions to $\R^2$ are separated from the direction "to the right" or "to the
left": for certainty, let us make the angles between the direction "to the
right" and these projected tangent directions at consecutive 6 points to be
equal to $\pi/4$, $3\pi/10$, $-\pi/4$, $3\pi/4$, $2\pi/10$, and $\pi/4$
respectively.

We call this knot a "large" one, and tie a very small homothetic knot on its
initial segment indicated by a tiny square in Fig. \ref{cat}. Then we shrink
very much the large knot in the "vertical" direction (orthogonal to the plane
of our picture) so that it becomes almost flat and its derivative almost
horizontal, not changing the small knot. Then we move this small knot along the
large one in such a way as if it would be frozen in a small hard bead put on
this large knot. (On the same Fig.~\ref{cat} we show by the thick lines the
channel of the bead; in this case all the picture should be considered as that
of the small homothetic knot. In particular all the points of this knot where
the direction of its derivative is sufficiently far from the standard one, are
inside the bead.)

More precisely, we associate with this bead a orthonormal frame in $\R^3$ whose
first vector in the initial instant is vertical (i.e. orthogonal to the plane
of the picture) and the second vector is directed along the channel.

\begin{lemma}
Suppose that a) the ratio of the diameter of the channel to its length is equal
to a sufficiently small number $\varepsilon$, b) the coefficient of the
flattening of the large knot in the vertical direction is of order
$\varepsilon^2$, so that the absolute value of the "vertical" part of the
derivative of the large knot shown in Fig.~\ref{cat} is nowhere greater than
$\varepsilon^2$ times the length of its "horizontal" part, and c) the size of
the bead $($i.e. the homothety coefficient of two knots$)$ is equal to
$\varepsilon^3$.

Then we can move our bead along the entire large knot in such a way that

A) the first vector of its associated frame remains vertical all the time, and

B) there is a smooth one-parametric family of long knots in $\R^3$ such that at
any instant

i) they coincide with the large knot everywhere outside the convex hull of the
bead,

ii) their intersection with the bead itself remains fixed and is as shown in
Fig.~\ref{cat},

iii) in all the points of the knot inside the channel of the bead, the angle
between the derivative of the knot and the direction of the channel is less
than $\pi/4$. \quad $\square$
\end{lemma}

The first and last instants of this one-parameter family of knots obviously can
be joined by a homotopy not changing the topology of the knot diagram, and we
get a closed loop in the space of knots. Now let us calculate the intersection
number of this loop with the chain described in Theorem \ref{main}.

This loop can intersect the varieties indicated in statements a) and c) of this
Theorem only when triple intersections of the projection occur. This can happen
only if one of crossing points of the smaller knot moving along some branch of
the large knot passes above or below the other its branch: in total 18
suspicious instants. These instants should be counted with multiplicities. In
the case of variety described in statement a) the multiplicity is equal (mod 2)
to the number of other crossing points of the composite knot forming together
with this triple point a configuration satisfying all other conditions of this
statement; for variety described in c) the multiplicity is equal to 0 or 1
depending on the condition on the tangent frame.

It is easy to calculate that the desired configurations for the variety a)
exist only when our small knot passes the first time (i.e. along the lower
branch) the third crossing point of the large knot: moreover, all three
instants when one of crossing points of the small knot pass this point have
multiplicity 1. Therefore the total number of intersections of our path with
variety a) is equal to 3. Similarly, we meet the variety c) only once, when our
small knot (more precisely, its second crossing point) passes the first time
the first crossing point of the large knot. So, the intersection number with
variety c) is equal to 1.

The configurations of type b) can appear by two reasons. First, when the small
knot passes a crossing point of the large one (and namely an undercrossing)
then all its points go under the other branch of the large knot; at some
instant this happens with the point with the distinguished tangent direction.
Again, any such instant should be counted with multiplicities depending on the
order of other crossing points of the composite knot. It is easy to calculate
that only once this multiplicity can be not equal to zero. Namely, when our
small knot undercrosses the third crossing point, then at some instant this
situation appears with multiplicity 2. Further, when our small knot moves and
rotates together with the derivative of the large one, some of tangent lines at
its own crossing points can instantly become directed "to the right". (Namely,
only the tangent line at the undercrossing branch of the first or third
crossing point of the small knot is interesting for us.) There are exactly two
points of the large knot at which it happens: in Fig.~\ref{cat} they are
indicated by small circles. The multiplicity of the "lower" (in this picture)
point is equal to 1, and the multiplicity of the "upper" one is equal to 0.

Finally, the total number of intersection points of our path with the variety
indicated in Theorem \ref{main} is equal to $3+1+2+1=7$, and proposition
\ref{realiz} is completely proved.

\section{Comments on and proof of Theorem \protect{\ref{comain}}}
\label{proof2}

Four statements A, B, C, D of this theorem are discussed in corresponding parts
of this section.
\smallskip

{\bf A.} The variety in ${\mathcal K}_n$ given by the condition $f(0)=f(\pi)$
is a vector subspace of codimension $n.$ It is equal to the boundary of the
variety ${\mathcal A}$ described in statement A of Theorem \ref{comain}. The
map ${\mathcal K}_n \to \R^{n-1},$ sending any curve $f$ to the vector
$f_1(\pi)-f_1(0),$ defines an isomorphism between $\R^{n-1}$ and the normal
bundle of ${\mathcal A}$, in particular induces a coorientation of ${\mathcal
A}$ from any orientation of $\R^{n-1}$. Thus for any integral
$(n-1)$-dimensional cycle in ${\mathcal K}_n \setminus \Sigma$ its intersection
index with ${\mathcal A}$ is well defined and is equal to its linking number
with the subspace $\{f|f(0)=f(\pi)\}$. It follows from calculations in
\cite{arman}, \cite{fasis}, \cite{bjo} that such integral cycles exist only if
$n$ is even. For instance, let $S^{2k-1}$ be the unit sphere, and consider all
the fibers of the Hopf bundle $S^{2k-1} \to \CP^{k-1}$ supplied with natural
parametrizations respecting the natural orientations of these fibers. The set
of all these parametrized fibers is obviously homeomorphic to $S^{2k-1}$ (to
any parametrized fiber there corresponds the zero of the parameter) and has
exactly one intersection point with the variety ${\mathcal A}$.
\smallskip

{\bf B.} The algorithm of finding the spanning chain is as follows. The variety
${\mathcal L}$ described in Proposition \ref{ordone} is swept out by
1-parametric family of subspaces $L(\alpha) \subset {\mathcal K}_n$ of
codimension $n$: they are parametrized by points $\alpha$ of the half-circle
$S^1/\pm = \R^1/ \pi \Z$ and defined by conditions $f(\alpha)=f(\alpha+\pi)$.
Let us try to span all these spaces separately. Consider the trivial bundle
${\mathcal K}_n \times [0,\pi] \to [0,\pi]$ and subset in it consisting of
pairs $(\alpha,f)$ such that $f(\alpha)$ is above $f(\alpha+\pi)$ in $\R^n.$
This subset is a smooth submanifold with boundary, and its projection to
$[0,\pi]$ is a smooth fiber bundle. Forgetting the second coordinate $\alpha$
defines the projection of this manifold to ${\mathcal K}_n$. Its image is
exactly the variety ${\mathcal B}a$ described in statement Ba of Theorem
\ref{comain}. Its boundary consists of the variety ${\mathcal L}$ and images of
fibers of the above-described fiber bundle over the points $0$ and $\pi$. The
union of these two fibers is equal to the subspace distinguished by the
condition $f_1(0)=f_1(\pi)$, and is equal to the boundary of the half-space
${\mathcal B}b$ described in statement Bb of Theorem \ref{comain}.

Now we choose coorientations of these varieties. The variety ${\mathcal B}a$ is
singular. Any its regular point $f$ satisfies the condition
$f_1(\alpha)=f_1(\alpha+\pi)$ for exactly one $\alpha \in [0,\pi)$ and has
transverse self-intersection of the curve $f_1(S^1)$ at this point. Close to
such a point $f$ the coorientation of ${\mathcal B}a$ is defined as follows.
Fix our point $\alpha$ and define the map $({\mathcal K}_n,f) \to TS^{n-2}$
associating to any parametrized curve $g \approx f$ the point of $S^{n-2}$
equal to the direction of the vector $g'_1(\alpha+\pi)-g'_1(\alpha)$, and the
tangent vector at this point in $S^{n-2}$ equal to the projection of the vector
$g(\alpha+\pi)-g(\alpha)$ to the plane orthogonal to this vector
$g'_1(\alpha+\pi)-g'_1(\alpha)$. The preimage of the zero section of $TS^{n-2}$
under this map is tangent in ${\mathcal K}_n$ to the variety ${\mathcal B}a$,
in particular if we have a generic germ of a $(n-2)$-dimensional subvariety
(simplex) in ${\mathcal K}_n$ at the point $g$ then it is transversal to both
varieties and we can induce its desired orientation from (any fixed)
orientation of the bundle $TS^{n-2}$.

The coorientation of the variety ${\mathcal B}b$ is induced from a chosen
orientation of $\R^{n-2}$ by the map ${\mathcal K}_n \to \R^{n-2}$ by a map
sending any $f$ to the direction of the vector $f_2(\pi)-f_2(0) \in \R^{n-2}
\equiv \R^n/\{\uparrow, \mapsto\}$.
\smallskip

{\bf C}. Recall that the term $F_1$ of the simplicial resolution of $\Sigma$ is
the space of pairs
\begin{equation}
\label{f1} ((\alpha,\beta),f) \in \overline{B(S^1,2)} \times {\mathcal K}_n
\end{equation}
such that $f(\alpha)=f(\beta)$. In particular it is a vector bundle over
$\overline{B(S^1,2)}$. Let ${\mathcal M} \subset F_2 \setminus F_1$ be the
principal part of the $(2n-3)$-dimensional class of order 2 described in
Proposition \ref{ordtwo}. Its first differential $d^1(\mathcal M)$ is realized
by the subvariety in $F_1$ consisting of such pairs (\ref{f1}) that
$\beta=\alpha+\pi$ and $f$ satisfies not only the condition
$f(\alpha)=f(\alpha+\pi)$ but also the condition
$f(\alpha+\pi/2)=f(\alpha-\pi/2).$ The set of such pairs $(\alpha,\beta)$ is
the circle $\R^1 /\pi \Z,$ so our cycle $d^1({\mathcal M})$ is the space of a
vector bundle over the circle. To span it in $F_1$ consider the subvariety
${\mathcal M}'_1 \subset F_1$ consisting of such pairs (\ref{f1}) that again
$\beta=\alpha+\pi,$ $f(\alpha)=f(\beta),$ but the image of one of points
$f(\alpha \pm \pi/2)$ is {\em above} the other: namely, the image of those of
these two points which is separated from $0\in S^1$ by the points $\alpha,
\alpha+\pi$ is above the image of its antipode. This subvariety also forms a
fiber bundle over the circle $\R^1/\pi \Z$ of all such pairs
$(\alpha,\alpha+\pi)$. There is exactly one position of $\alpha$ over which
this fiber bundle fails to be locally trivial, namely $\alpha=0 (\mbox{mod }
\pi).$ The boundary of this subvariety is equal to the sum of the cycle
$d^1({\mathcal M})$ and the space of points (\ref{f1}) where $\alpha=
0(\mbox{mod } \pi)$, $f(0)=f(\pi)$ and $f_1(\pi/2)=f_1(-\pi/2)$. We span the
latter space by the similar {\em half}space ${\mathcal M}''_1$, defined by the
condition that $f(0)=f(\pi)$ and $f_1(\pi/2)$ lies to the right of
$f_1(-\pi/2)$. The sum ${\mathcal M}'_1 +{\mathcal M}''_1$ is the desired chain
in $F_1$ whose boundary is equal to $d^1({\mathcal M})$. Now we consider the
image $d_2({\mathcal M})$ of this chain in $\Sigma$ and try to represent it as
a boundary of some relative cycle in ${\mathcal K}_n (\mbox{mod }\Sigma).$ The
image of ${\mathcal M}''_2$ is obvious, and the image of ${\mathcal M}'_1$
consists of maps $f$  such that there exists $\alpha \in [0,\pi]$ such that
$f(\alpha)=f(\alpha+\pi)$, and the image of one of points $\alpha\pm \pi/2$
(namely, the one separated from 0 by $\alpha$ and $\alpha+\pi$) is above the
other.

It is natural to kill this variety  by the space of all maps $f$ described in
statement Ca of Theorem \ref{comain}. Its boundary consists of this image of
${\mathcal M}'_1$
 and the space of such maps $f$ that $f(\pi)$ is above $f(0)$
and $f_1(\pi/2)=f_1(-\pi/2).$ The boundary of the variety described in
statement Cb of Theorem \ref{comain} is equal to the sum of the latter space
and the image of ${\mathcal M}''_1$.

Statement C of Theorem \ref{comain} is thus proved for $\Z_2$-homology; the
proof of its integer version requires additionally only an accounting of
orientations.
\smallskip

{\bf D}. We shall use the pictures like in \S~\ref{proof1}, only the Wilson
loop will be shown not by a segment but by an oval with marked "zero" point on
its top. This point is referred to as $0$ in subscripts, and all the other
points participating in the definition of cells and their subvarieties are
numbered accordingly to the (counterclockwise) orientation of the Wilson loop.
All the calculus remains the same as in \S \ref{proof1}, only the boundary
operators will include the limit positions of our cells and their subvarieties
when some of defining them points tend to $0$.

As we are interested in integral homology classes, we shall take care of
orientations of all our varieties in the cells of the standard cell
decompositions of terms $F_i \setminus F_{i-1}.$ This orientation consists of
the orientation of the cell and the (co)orientation of the subvariety in it.
The choice of these orientations will follow the guidelines indicated in \S 3
of \cite{V1} or \S V.3.3 of \cite{fasis}. Namely, they consist of the following
orientations (taken in that order): a) the orientation of the simplex
participating in the construction of the simplicial resolution (i.e. the
simplex $\tilde \Delta(J)$ or some its non-marginal face); b) the coorientation
of the subspace $L(J)$ of the space  $ {\mathcal K}_n$; c) the orientation of
the space of equivalent point configurations $J \subset S^1$; d) the
(co)orientation of the subvariety in the cell. The first three orientations are
specified exactly as in \S 3 of \cite{V1} (but now in c) we can move only {\em
nonzero} points). Often the subvariety in the cell is given by several
conditions of the form: "there are additional points in $\R^1$ whose images
$f(\cdot) \in \R^n$ (or their projections to some fixed subspace) coincide with
one another or with images of some points participating in the definition of
the cell", or at least our subvariety forms an open subset in a subvariety
defined in such a way. In this case the orientation d) also is defined by the
sequence consisting of $\alpha$) the (co)orientation of the vector subspace
defined by these conditions in the vector spaces counted in the step b) above,
and $\beta$) the orientation of the space of configurations of additional
points. These orientations also are specified as in \cite{V1}, \cite{fasis}; to
define the coorientations of subspaces we assume that the direction "up" in
$\R^n$ is the first vector of the canonical frame, and the direction "to the
right" is the second. All the forthcoming calculations refer to exactly this
choice of orientations.
\medskip

The principal part of the considered class in $F_2\setminus F_1$ is as follows:
\begin{equation}
\label{cprinpar} \unitlength 1.00mm \special{em:linewidth 0.4pt}
\linethickness{0.4pt}
\begin{picture}(80.00,14.50)
\put(25.00,9.00){\oval(30.00,10.00)[]}
\put(3.00,9.00){\makebox(0,0)[cc]{$V_2=$}}
\put(45.00,9.00){\makebox(0,0)[cc]{$+$}} \put(65.00,9.00){\oval(30.00,10.00)[]}
\put(65.00,14.00){\circle*{1.00}} \put(25.00,14.00){\circle*{1.00}}
\put(22.00,4.00){\oval(12.00,8.00)[t]} \put(28.00,4.00){\oval(12.00,8.00)[b]}
\put(65.00,4.00){\oval(16.00,8.00)[t]} \put(61.00,4.00){\oval(8.00,4.00)[t]}
\put(69.00,4.00){\oval(8.00,4.00)[t]}
\end{picture}
 \ ,
\end{equation}
see \cite{V1}. The second summand has no boundary in $F_1$, the boundary of the
first is as follows:

\begin{equation}
\label{cdone} \unitlength 1.00mm \special{em:linewidth 0.4pt}
\linethickness{0.4pt}
\begin{picture}(79.00,15.50)
\put(24.00,10.00){\oval(30.00,10.00)[]} \put(64.00,10.00){\oval(30.00,10.00)[]}
\put(64.00,15.00){\circle*{1.00}} \put(24.00,15.00){\circle*{1.00}}
\put(24.00,5.00){\oval(12.00,8.00)[t]} \put(64.00,5.00){\oval(12.00,8.00)[t]}
\put(64.00,5.00){\line(1,-1){5.00}} \put(69.00,0.00){\line(1,1){5.00}}
\put(24.00,5.00){\line(-1,-1){5.00}} \put(19.00,0.00){\line(-1,1){5.00}}
\put(44.00,10.00){\makebox(0,0)[cc]{$-$}}
\put(3.00,10.00){\makebox(0,0)[cc]{$(-1)^n$}}
\end{picture}
 \ .
\end{equation}

Arguing as previously, we try to span the two terms of this chain by the
varieties encoded in left parts of the next two equations, respectively:

\begin{equation}
\label{cspan1} \unitlength 1.00mm \special{em:linewidth 0.4pt}
\linethickness{0.4pt}
\begin{picture}(123.00,95.00)
\put(26.00,14.00){\oval(26.00,10.00)[]} \put(26.00,19.00){\circle*{1.00}}
\put(63.00,14.00){\oval(26.00,10.00)[]} \put(63.00,19.00){\circle*{1.00}}
\put(105.00,14.00){\oval(26.00,10.00)[]} \put(105.00,19.00){\circle*{1.00}}
\put(20.00,35.00){\oval(26.00,10.00)[]} \put(20.00,40.00){\circle*{1.00}}
\put(63.00,35.00){\oval(26.00,10.00)[]} \put(63.00,40.00){\circle*{1.00}}
\put(105.00,35.00){\oval(26.00,10.00)[]} \put(105.00,40.00){\circle*{1.00}}
\put(20.00,63.00){\oval(26.00,10.00)[]} \put(20.00,68.00){\circle*{1.00}}
\put(63.00,63.00){\oval(26.00,10.00)[]} \put(63.00,68.00){\circle*{1.00}}
\put(105.00,63.00){\oval(26.00,10.00)[]} \put(105.00,68.00){\circle*{1.00}}
\put(20.00,85.00){\oval(26.00,10.00)[]} \put(20.00,90.00){\circle*{1.00}}
\put(66.00,85.00){\oval(26.00,10.00)[]} \put(66.00,90.00){\circle*{1.00}}
\put(105.00,85.00){\oval(26.00,10.00)[]} \put(105.00,90.00){\circle*{1.00}}
\put(3.00,85.00){\makebox(0,0)[cc]{$\partial$}}
\put(44.00,85.00){\makebox(0,0)[cc]{$=(-1)^{n-1}$}}
\put(85.00,85.00){\makebox(0,0)[cc]{$+$}}
\put(123.00,85.00){\makebox(0,0)[cc]{$+$}}
\put(41.00,63.00){\makebox(0,0)[cc]{$-$}}
\put(84.00,63.00){\makebox(0,0)[cc]{$-$}}
\put(3.00,35.00){\makebox(0,0)[cc]{$\partial$}}
\put(41.00,35.00){\makebox(0,0)[cc]{$=(-1)^n$}}
\put(84.00,35.00){\makebox(0,0)[cc]{$-(-1)^n$}}
\put(123.00,35.00){\makebox(0,0)[cc]{$+$}}
\put(44.00,14.00){\makebox(0,0)[cc]{$-$}}
\put(84.00,14.00){\makebox(0,0)[cc]{$-(-1)^n$}}
\put(22.00,80.00){\oval(10.00,8.00)[t]} \put(22.00,80.00){\line(-1,-1){5.00}}
\put(17.00,75.00){\vector(-1,1){5.00}} \put(68.00,80.00){\oval(10.00,8.00)[t]}
\put(68.00,80.00){\line(-1,-1){5.00}} \put(63.00,75.00){\line(-1,1){5.00}}
\put(105.00,80.00){\line(0,-1){5.00}} \put(105.00,75.00){\line(-1,0){14.00}}
\put(91.00,75.00){\line(0,1){20.00}} \put(91.00,95.00){\line(1,0){14.00}}
\put(105.00,95.00){\vector(0,-1){4.50}} \put(105.00,58.00){\line(0,1){10.00}}
\put(110.00,58.00){\line(-1,-1){5.00}} \put(105.00,53.00){\vector(-1,1){5.00}}
\put(63.00,58.00){\oval(14.00,8.00)[t]} \put(66.00,58.00){\line(-1,-1){5.00}}
\put(61.00,53.00){\vector(-1,1){5.00}} \put(20.00,58.00){\oval(14.00,8.00)[t]}
\put(16.00,52.00){\vector(4,1){8.00}} \put(16.00,52.00){\vector(4,-1){8.00}}
\put(20.00,48.00){\makebox(0,0)[cc]{$2$}}
\put(20.00,56.00){\makebox(0,0)[cc]{$1$}}
\put(18.00,30.00){\oval(10.00,8.00)[t]} \put(18.00,30.00){\line(1,-1){5.00}}
\put(23.00,25.00){\vector(1,1){5.00}} \put(61.00,30.00){\oval(10.00,8.00)[t]}
\put(61.00,30.00){\line(1,-1){5.00}} \put(66.00,25.00){\line(1,1){5.00}}
\put(105.00,30.00){\line(0,-1){5.00}} \put(105.00,25.00){\line(1,0){14.00}}
\put(119.00,25.00){\line(0,1){20.00}} \put(119.00,45.00){\line(-1,0){14.00}}
\put(105.00,45.00){\vector(0,-1){4.50}}
\put(105.00,30.00){\oval(14.00,8.00)[t]}
\put(105.00,80.00){\oval(14.00,8.00)[t]} \put(105.00,9.00){\line(0,1){10.00}}
\put(100.00,9.00){\line(1,-1){5.00}} \put(105.00,4.00){\vector(1,1){5.00}}
\put(63.00,9.00){\oval(14.00,8.00)[t]} \put(60.00,9.00){\line(1,-1){5.00}}
\put(65.00,4.00){\vector(1,1){5.00}} \put(26.00,9.00){\oval(14.00,8.00)[t]}
\put(22.00,3.00){\vector(4,1){8.00}} \put(22.00,3.00){\vector(4,-1){8.00}}
\put(26.00,0.00){\makebox(0,0)[cc]{$2$}}
\put(26.00,7.00){\makebox(0,0)[cc]{$1$}}
\put(3.00,63.00){\makebox(0,0)[cc]{$+$}}
\put(6.00,14.00){\makebox(0,0)[cc]{$+(-1)^n$}}
\end{picture}
\end{equation}

The varieties shown by the second from the end pictures in both these equations
coincide geometrically, and their canonical orientations differ by the factor
$(-1)^{n-1}$. Therefore the linear combination of left parts of these equations
taken with coefficients $-1$ and $(-1)^{n-1}$ respectively is equal in $F_1$ to
the sum of the expression (\ref{cdone}) and two last varieties in these
equations. If $n=3$ then the sum of last two varieties is equal to zero.
Indeed, any of these varieties consists of pairs (\ref{f1}) with $\alpha=0,$
$f(\alpha)=f(\beta),$ taken with some multiplicities. These multiplicities
always are opposite, because they are equal (up to signs) to different
combinatorial expressions for the linking numbers of two ``smoothened'' loops
into which the point $f(0)=f(\beta)$ breaks the curve $f(S^1).$

However, if $n>3$ then the sum of these two varieties is only homologous to
zero, but not equal to it. We shall encode this sum by the picture in the
right-hand part of the following equation:

\begin{equation}
\label{cspan2} \unitlength 1.00mm \linethickness{0.4pt}
\begin{picture}(91.00,10.50)
\put(3.00,5.00){\makebox(0,0)[cc]{$\partial$}}
\put(24.00,5.00){\oval(32.00,10.00)[]} \put(24.00,0.00){\line(0,1){10.00}}
\put(24.00,10.00){\circle*{1.00}} \put(6.00,2.00){\line(0,1){6.00}}
\put(6.00,5.00){\vector(1,0){34.00}} \put(50.00,5.00){\makebox(0,0)[cc]{$=
(-1)^n$}} \put(75.00,5.00){\oval(32.00,10.00)[]}
\put(75.00,10.00){\circle*{1.00}} \put(75.00,10.00){\line(0,-1){10.00}}
\put(75.00,5.00){\vector(1,0){16.00}} \put(75.00,5.00){\vector(-1,0){16.00}}
\end{picture}
 \ .
\end{equation}

The variety assumed in the left part of this equation consists of all pairs
(\ref{f1}) in $F_1$ such that $\alpha=0$ and additionally there are points
$\gamma \in (0,\beta)$ and $\delta \in (\beta,2\pi)$ such that the projection
of $f(\delta)$ to $\R^{n-1}$ lies "to the right" from that of $f(\beta)$.

So, the desired chain in $F_1$ spanning the cycle (\ref{cdone}) is equal to

\begin{equation}
\label{cspanlst} - \left( \unitlength 1.00mm \special{em:linewidth 0.4pt}
\linethickness{0.4pt}
\begin{picture}(108.00,8.50)
\put(12.00,3.00){\oval(24.00,10.00)[]}
\put(32.00,3.00){\makebox(0,0)[cc]{$+(-1)^n$}}
\put(52.00,3.00){\oval(24.00,10.00)[]}
\put(74.00,3.00){\makebox(0,0)[cc]{+$(-1)^n$}}
\put(96.00,3.00){\oval(24.00,10.00)[]} \put(82.00,0.00){\line(0,1){6.00}}
\put(82.00,3.00){\vector(1,0){26.00}} \put(50.00,-2.00){\oval(10.00,8.00)[t]}
\put(50.00,-2.00){\line(1,-1){5.00}} \put(55.00,-7.00){\vector(1,1){5.00}}
\put(52.00,8.00){\circle*{1.00}} \put(96.00,8.00){\circle*{1.00}}
\put(96.00,8.00){\line(0,-1){10.00}} \put(12.00,8.00){\circle*{1.00}}
\put(14.00,-2.00){\oval(10.00,8.00)[t]} \put(15.00,-2.00){\line(-1,-1){5.00}}
\put(10.00,-7.00){\vector(-1,1){5.00}}
\end{picture} \right) .
\end{equation}

Its image in $\Sigma$ is expressed by the formula

\begin{equation}
\label{cd-2} -\left( \unitlength 1.00mm \special{em:linewidth 0.4pt}
\linethickness{0.4pt}
\begin{picture}(101.00,8.50)
\put(12.00,3.00){\oval(24.00,10.00)[]} \put(34.00,3.00){\makebox(0,0)[cc]{$+
(-1)^n$}} \put(53.00,3.00){\oval(24.00,10.00)[]}
\put(70.00,3.00){\makebox(0,0)[cc]{$+$}} \put(89.00,3.00){\oval(24.00,10.00)[]}
\put(75.00,0.00){\line(0,1){6.00}} \put(75.00,3.00){\vector(1,0){26.00}}
\put(51.00,-2.00){\line(1,-1){5.00}} \put(56.00,-7.00){\vector(1,1){5.00}}
\put(53.00,8.00){\circle*{1.00}} \put(89.00,8.00){\circle*{1.00}}
\put(12.00,8.00){\circle*{1.00}} \put(15.00,-2.00){\line(-1,-1){5.00}}
\put(10.00,-7.00){\vector(-1,1){5.00}} \put(10.00,-2.00){\line(1,1){5.00}}
\put(15.00,3.00){\line(1,-1){5.00}} \put(56.00,-2.00){\line(-1,1){5.00}}
\put(51.00,3.00){\line(-1,-1){5.00}} \put(89.00,-2.00){\line(1,1){4.00}}
\put(93.00,2.00){\line(-2,3){4.00}}
\end{picture}
 \right) .
\end{equation}

We need to span this chain by a relative cycle in ${\mathcal K}_n$ (mod
$\Sigma$). For this we have

\begin{equation}
\label{ppvv} \unitlength 1.00mm \special{em:linewidth 0.4pt}
\linethickness{0.4pt}
\begin{picture}(124.00,40.50)
\put(2.00,35.00){\makebox(0,0)[cc]{$\partial$}}
\put(15.00,35.00){\oval(22.00,10.00)[]} \put(15.00,40.00){\circle*{1.00}}
\put(57.00,35.00){\oval(22.00,10.00)[]} \put(57.00,40.00){\circle*{1.00}}
\put(85.00,35.00){\oval(22.00,10.00)[]} \put(85.00,40.00){\circle*{1.00}}
\put(111.00,35.00){\oval(20.00,10.00)[]} \put(111.00,40.00){\circle*{1.00}}
\put(23.00,10.00){\oval(24.00,10.00)[]} \put(23.00,15.00){\circle*{1.00}}
\put(53.00,10.00){\oval(24.00,10.00)[]} \put(53.00,15.00){\circle*{1.00}}
\put(83.00,10.00){\oval(24.00,10.00)[]} \put(83.00,15.00){\circle*{1.00}}
\put(112.00,10.00){\oval(20.00,10.00)[]} \put(112.00,15.00){\circle*{1.00}}
\put(18.00,30.00){\line(-1,1){5.00}} \put(13.00,35.00){\vector(-1,-1){5.00}}
\put(12.00,30.00){\line(1,-1){5.00}} \put(17.00,25.00){\vector(1,1){5.00}}
\put(36.00,35.00){\makebox(0,0)[cc]{$=(-1)^{n-1}$}}
\put(64.00,30.00){\line(-1,-1){5.00}} \put(59.00,25.00){\line(-1,1){5.00}}
\put(61.00,30.00){\line(-1,1){5.00}} \put(56.00,35.00){\vector(-1,-1){5.00}}
\put(78.00,30.00){\line(1,1){5.00}} \put(83.00,35.00){\line(1,-1){5.00}}
\put(81.00,30.00){\line(1,-1){5.00}} \put(86.00,25.00){\vector(1,1){5.00}}
\put(71.00,35.00){\makebox(0,0)[cc]{$-$}}
\put(99.00,35.00){\makebox(0,0)[cc]{$-$}}
\put(124.00,35.00){\makebox(0,0)[cc]{$+$}}
\put(106.00,30.00){\line(1,-1){5.00}} \put(111.00,25.00){\vector(1,1){5.00}}
\put(111.00,30.00){\vector(0,1){9.50}} \put(112.00,5.00){\vector(0,1){9.50}}
\put(117.00,5.00){\line(-1,-1){5.00}} \put(112.00,0.00){\vector(-1,1){5.00}}
\put(81.00,5.00){\line(1,-1){5.00}} \put(86.00,0.00){\vector(1,1){5.00}}
\put(90.50,5.00){\line(-3,2){7.50}} \put(83.00,10.00){\vector(-3,-2){7.50}}
\put(53.00,5.00){\line(3,-4){3.75}} \put(56.75,0.00){\vector(3,4){3.75}}
\put(53.00,5.00){\line(-3,4){3.75}} \put(49.25,10.00){\vector(-3,-4){3.75}}
\put(25.00,5.00){\line(-1,1){5.00}} \put(20.00,10.00){\vector(-1,-1){5.00}}
\put(15.00,5.00){\line(3,-2){7.50}} \put(22.50,0.00){\vector(3,2){7.50}}
\put(38.00,10.00){\makebox(0,0)[cc]{$+$}}
\put(8.00,10.00){\makebox(0,0)[cc]{$+$}}
\put(68.00,10.00){\makebox(0,0)[cc]{$+$}}
\put(98.00,10.00){\makebox(0,0)[cc]{$-$}}
\end{picture} ,
\end{equation}

\begin{equation}
\label{cnev1} \unitlength 1.00mm \special{em:linewidth 0.4pt}
\linethickness{0.4pt}
\begin{picture}(119.00,10.50)
\put(1.00,5.00){\makebox(0,0)[cc]{$\partial$}}
\put(19.00,5.00){\oval(26.00,10.00)[]} \put(4.00,8.00){\line(0,-1){6.00}}
\put(4.00,5.00){\vector(1,0){28.00}} \put(19.00,10.00){\circle*{1.00}}
\put(19.00,0.00){\vector(0,1){9.50}} \put(38.00,5.00){\makebox(0,0)[cc]{$=-$}}
\put(44.00,2.00){\line(0,1){6.00}} \put(44.00,5.00){\vector(1,0){28.00}}
\put(59.00,5.00){\oval(26.00,10.00)[]} \put(59.00,10.00){\circle*{1.00}}
\put(59.00,0.00){\line(1,1){4.00}} \put(63.00,4.00){\line(-2,3){4.00}}
\put(79.00,5.00){\makebox(0,0)[cc]{$+$}}
\put(87.50,5.00){\makebox(0,0)[cc]{$(-1)^n$}}
\put(106.00,5.00){\oval(26.00,10.00)[]} \put(106.00,10.00){\circle*{1.00}}
\put(106.00,5.00){\vector(1,0){13.00}} \put(106.00,5.00){\vector(-1,0){13.00}}
\put(106.00,0.00){\vector(0,1){9.50}}
\end{picture}
 .
\end{equation}

The sum of the fourth, fifth and sixth terms in the right-hand part of
(\ref{ppvv}) is equal to zero. The sum of varieties encoded by the pictures in
the third and seventh terms is equal to the variety shown by the last picture
in (\ref{cnev1}). Therefore the desired relative cycle is equal to the linear
combination of left parts of (\ref{ppvv}) and (\ref{cnev1}) taken with
coefficients $1$ and $(-1)^n$ respectively. Theorem \ref{comain} is completely
proved.


\begin{thebibliography}{99}

\bibitem{arnold} V.~I.~Arnold, {\it On some topological
invariants of algebraic functions,} Trudy Moskov. Mat. Obshch. 1970, 21,
27--46. Engl. transl.: Transact. Moscow Math. Soc., 1970, 21, 30--52.

\bibitem{BN} D.~Bar-Natan, {\em On the Vassiliev Knot Invariants},
Topology {\bf 34}, No. 2, 1995, 423-472.

\bibitem{BL} J.~S.~Birman, X.-S.~Lin, {\em Knot polynomials and
Vassiliev's invariants}, Invent. Math. 1993, {\bf 111}(2), 225--270.

\bibitem{cart} P.~Cartier, {\it Construction combinatoire des invariants
de Vassiliev,} C.R.Acad. Sci. Paris, S\'erie I, {\bf 316} (1993), 1205--1210.

\bibitem{cotta}
A.~S.~Cattaneo, P.~Cotta-Ramusino, R.~Longoni, {\it Configuration spaces and
Vassiliev classes in any dimension,} preprint, 1999-2000.

\bibitem{GPV} M.~Goussarov, M.~Polyak, and O.~Viro,
{\it Finite type invariants of classical and virtual knots,} Topology 2000;
{\tt http://www.math.uu.se/\,$\tilde{}$\,viro/pw.ps}

\bibitem{lannes} J.~Lannes, {\it Sur l'invariants de Vassiliev de degr\'e
inf\'erieur ou \'egal \`a 3.} L'Enseignement Math\'ematique, 1993, {\bf 39}
(3-4), 295--316.

\bibitem{MM} A.~B.~Merkov,
{\it Vassiliev invariants classify plane curves and doodles,} preprint 1998,
{\tt http://www.botik.ru/\,$\tilde{}$\,duzhin/as-papers/ficd-dvi.zip}

\bibitem{M} A.~B.~Merkov,
{\it Segment--arrow diagrams and invariants of ornaments,} to appear in Mat.
Sbornik.

\bibitem{PV} M.~Polyak and O.~Viro,
{\it Gauss diagram formulas for Vassiliev invariants,} {\em Internat. Math.
Res. Notes} {\bf 11}, 1994, 445--453.

\bibitem{tetu} D.~M.~Teiblum, V.~E.~Turchin, {\it Private communication,}
1995.

\bibitem{Turchin} V.~Tourtchine,
{\it Sur l'homologie des espaces des n\oe uds non-compacts,} Submitted to Publ.
IHES, 2000.

\bibitem{sveta} S.~D.~Tyurina, {\it On the Lannes and Viro-Polyak type
formulas for finite type invarians}, Matem. Zametki (Math. Notes), {\bf 65:}6
(1999).

\bibitem{V1} V.~A.~Vassiliev,
{\it Cohomology of knot spaces}, in: Adv. in Sov. Math.; Theory of
Singularities and its Appl. (V.~I.~Arnold, ed.). AMS, Providence, R.I., 1990,
p.~23--69.

\bibitem{congr} V.~A.~Vassiliev, {\it Complexes of connected graphs,}
in: The I.~M.~Gel'fand's mathematical seminars 1990--1992 (L.~Corvin,
I.~Gel'fand, J.~Lepovsky, eds.); Birkh\"auser, Basel, 1993, 223--235.

\bibitem{fasis} V.~A.~Vassiliev, {\it Topology of complements
of discriminants,} Moscow, Phasis, 1997 (in Russian).

\bibitem{arman} V.~A.~Vassiliev,
{\it On invariants and homology of spaces of knots in arbitrary manifolds,} in:
Topics in Quantum Groups and Finite-Type Invariants. Mathematics at the
Independent University of Moscow (B.~Feigin and V.~Vassiliev, eds.). AMS
Translations. Ser.~2. {\bf 185}, 1998, AMS, Providence RI, 155--182.

\bibitem{bjo} V.~A.~Vassiliev, {\em Topology of two-connected graphs and
homology of spaces of knots,} in: Differential and symplectic topology of knots
and curves (S.~L.~Tabachnikov, ed.); AMS Transl., Ser. 2, vol. 190, AMS,
Providence RI, 1999, 253--286.

\bibitem{kotor} V.~A.~Vassiliev,
{\em Topological order complexes and resolutions of discriminant sets,}
Publications de l'Institut Math\'ematique Belgrade, Nouvelle s\'erie {\bf
66(80)} (1999), 165--185.

\bibitem{ZZ} Ziegler, G.~M. and \v{Z}ivaljevi\'c, R.~T.
{\em Homotopy type of arrangements via diagrams of spaces,} Math. Ann. {\bf
295} (1993), 527--548.

\end{thebibliography}
\end{document}